\documentclass[]{article}
\usepackage[letterpaper,margin=01.in]{geometry}
\usepackage{amssymb, amsmath}
\usepackage{graphicx}
\usepackage[utf8]{inputenc}
\usepackage{subcaption}
\usepackage[shortlabels]{enumitem}
\title{Myocardial ischemic effects on cardiac electro-mechanical activity}
\author{B.V. Rathish Kumar \\ Department of Mathematics and Statistics, Indian Institue of Technology Kanpur, India \\ Meena Pargaei \\ Department of Mathematics and Statistics, Indian Institue of Technology Kanpur, India \\ Govt. P. G. College Champawat, Uttarakhand, India, \\ Luca F. Pavarino \\ Department of Mathematics, University of Pavia, \\ Simone Scacchi \\ Department of Mathematics, University of Milan}

\begin{document}

\maketitle

\begin{abstract}
In this work, we investigated the effect of varying strength of Hyperkalemia and Hypoxia, in a human cardiac tissue with a local ischemic subregion, on the electrical and mechanical activity of healthy and ischemic zones of the cardiac muscle. The Monodomain model in a deforming domain is taken with the addition of mechanical feedback and stretch activated channel current coupled with the ten Tusscher human ventricular membrane model. The equations of finite elasticity are used to describe the deformation of the cardiac tissue. The resulting coupled electro-mechanical PDEs-ODEs non-linear system is solved numerically using finite elements in space and finite difference method in time. We examined the effect of local ischemia on the cardiac electrical and mechanical activity in different cases. We concluded that the spread of Hyperkalemic or Hypoxic region alters the electro-mechanical coupling in terms of the action potential ($v$), intracellular calcium ion concentration $[Ca^{+2}]_i$, active tension, ($T_A$), stretch ($\lambda$), stretch rate ($ \frac{d \lambda}{dt}$). With the increase in the size of the ischemic region by factor of five, approximately $45\%$ variation in the stretch rate $\frac{d \lambda}{dt}$ is noticed. It is also shown that ischemia affects the deformation (expansion and contraction) of the heart.
\end{abstract}

\section{Introduction}
The heart is considered as a muscular pump which is connected to the systemic and pulmonary vascular systems. The main purpose of
the heart and vascular system is to sustain an appropriate supply of nutrients as oxygenated blood and metabolic substrates to the 
body. The left atrium and ventricle chambers pump the blood from pulmonary veins to the aorta, whereas the right atrium and 
ventricle chambers pump the blood from vena cava to pulmonary arteries. The cardiac cycle consists of two phases: diastole, 
during which the ventricles relax and the heart fills with blood; and systole, during which the ventricles contract and pump
blood out of the heart to arteries. 
To this end, the heart needs three types of cells: 
1. SA node or pacemaker cells, which produce an electrical signal; 2. Conductors to spread the pacemaker signals; and 3. 
Contractile cells, to mechanically pump blood. Pacemaker cells start the electrical sequence of depolarization and repolarization. 
The electrical signal generated from the SA node travels to the ventricles via the atrioventricular node(AV node), the bundle of His, the  right and  left bundle branches and Purkinje fibers. As the depolarization signal reaches the contractile cells 
they start to contract and as the repolarization signal reaches the myocardial cells, they start to relax. In this way the
electrical signals and the mechanical pumping action of the heart are connected via excitation-contraction coupling.
Heart's electrical and mechanical activity can be deduced from the ECG pattern. P wave represents the atrial depolarization, 
QRS complex due to the ventricular depolarization also start of the ventricular contraction, T wave indicate the ventricular 
repolarization and the beginning of ventricular relaxation.

There exist one factor, called stretch, which is responsible for the activation and inactivation of the mechanically controlled 
ion channels. In 1997, Hu \& Sachs \cite{Hu} revealed the stretch activated channels (SACs) in cardiac cells. Positive response
of the SACs contribute for heart to stretch. Cardiac muscle stretch also leads to the increase in active tension. 

For the finite elasticity material, there are various models for the strain energy function based on the different
laws. For example, in \cite{elastic1} Guccione et al. used the transversely isotropic exponential law strain energy function and
the orthotropic and isotropic laws have been used in \cite{elastic2} and 
\cite{elastic3}, respectively. Various models for the development of the active tension and the description
of the dynamics of calcium ions and the cross bridges bindings have been proposed in \cite{VRPM,HSME,TAGG,RMA,RJAL,MA,STRA,JDRM,SJMJ,PJ,SSJEN,XVN,VTJHN,VJJN,BBJR,JYAN,IJH}. In \cite{hunter},
Hunter et al. presented the model based on the review of experimental data for the mechanics of active and
passive cardiac muscle. They concluded that it could be interpreted as the four state variable model. (i) myocardial tissue passive elasticity, (ii) tension dependence on the troponin C and $Ca^{+2}$ bindings, (iii)
tropomyosin movement kinetics and the actively cycling cross-bridges binding sites and the length dependence
of this process, and, (iv) cross-bridge tension development kinetics under variations of the myofilament length.
In \cite{smith}, developed a biophysically and dynamically stable model for the cardiac myocyte of rat at the room temperature. This model is capable of
determining the changes in the concentration of intracellular calcium ion in response to the stretch. In \cite{land}, 
presented a multiscale electro-mechanical model for the murince heart which is capable of investigating the effects of excitation-contraction
coupling. In this paper,they show a strong relation between the tension and the velocity to explain the differences between the 
single cell tension and the whole organ pressure transient.

For the cardiac electrical activity models at the tissue level, Bidomain and Monodomain models have been considered in various 
researches \cite{JYAN,BBJR,VJJN,XVN,SJMJ,HSME,MA,STRA,TAGG,RMA,RJAL}. For the ionic membrane models related to the different species has been used in
\cite{JYAN,BBJR,VJJN,XVN,SJMJ,HSME,RMA,RJAL,MA,STRA,TAGG} etc. 

Following the electro-mechanical coupling approach presented in \cite{IJH,BBJR,VJJN,JMD,STRA,SPDJ,FAR,PJ,DNP,RMA,SE,HSME,DGFA} we need to impose the electrical models on the deformed domain to propose the coupling between the electrical and mechanical model. Then by the 
Lagrangian framework again this model is reformulated into the reference domain. 

There exist coupled cardiac electro-mechanical models which explain the traveling of electric signal into the cardiac muscle and 
the contraction-relaxation process. These models consists : 1) the electrical  model at the tissue, which consists, the system of
non-linear parabolic reaction-diffusion equations of the degenerate type, 2) the cell level model to describe the flow of ionic 
currents via cellular membrane  consist, system of ordinary differential equations (ODEs), 3) deformation of cardiac tissue is 
modeled by the quasi-static finite elasticity model, 4)intracellular calcium dynamics and the cross bridges bindings are described
by the active tension model which consists of a system of non-linear ODEs.
Existence and uniqueness for the coupled electro-mechanical model is given in \cite{existence1,existence2}.

Cardiac arrhythmia involves the abnormality in the pattern of the heartbeat, if the electrical impulse fails to start from 
the SA node or if there is some abnormality in the impulse initiated from the SA node. Myocardial ischemia takes
place when there is abnormality in the blood flow to heart and in oxygen supply to the heart. It is one of the
main causes of sudden death. Cardiac arrhythmia and heart attack can also occur due to the myocardial
ischemia. Due to this myocardial ischemia metabolism, electrophysiological and mechanical changes appear which results
in the reduction of action potential duration (APD), action-potential upstroke, resting-potential shift, reduction in conduction
velocity (CV), active tension, stretch along the fiber and the stretch rate in the ischemic cells or tissue [2, 3, 4].

Various experimental and numerical simulation work have been carried out to analyze the electrophysiological and metabolic properties of the healthy and ischemic cells and tissue \cite{Dutta,Janse,Pogwizd,Tice,Lucaischemia,Rodriguez}.
Mathematical and computational studies
have been focused on investigating various physiological effects in order to understand the relationship between the 
electrophysiological and metabolic parameters. It can be analyzed from the past studies that these variations 
are basically due to, a) Hyperkalemia which affects the cardiac cell resting potential due to the 
increase in the concentration of extracellular potassium \cite{Schaapher,Pandit}, b) Hypoxia defined as the reduced oxygen supply, 
it reduces the cell metabolism and ADP to ATP concentration ratio in the cell get changed. Successively, 
this change affects the opening of the particular ATP-dependent potassium channels \cite{Van,Weiss}. This ATP-dependent potassium channels is also modulated by the mechanical environment. 

There are several research papers in which they have studied about ischemia in the animals. 
\cite{Fiolet,Carmeliet,Furukawa,Coronel,MA,Pandit1,Wilensky,Schaapher}. Since it is difficult to get the human data, and hence it is challenging to extend the properties from animals to human. But it can be achieved via computational modeling. Even though many human models have been built 
and analyzed via healthy cell data their appropriateness for the ischemia computation is not known. So, it is essential to examine the effect of varying ischemic parameters in the ischemic region. In \cite{panfilov,Sutton,Taggart}, studies have been done for the global ischemia in human under different ischemic conditions. In \cite{pargaei2019}, the authors has discussed the effect of local cardiac ischemia on the electrical activity of the ischemic region and the neighboring healthy region in a human cardiac tissue. 

For the electro-mechanical coupling models, the bio-electrical activity experiences three main feedback from the mechanical deformation:
\begin{enumerate}
	\item 
	conductivity feedback: the influence of the deformation gradient on the conductivity coefficients of the electric current flow model;
	\item 
	convection feedback: the influence of deformation gradient and deformation
	rate on the electric current flow model;
	\item 
	ionic feedback: the influence of stretch-activated membrane channels on the ionic current.
\end{enumerate}

In \cite{lucamech} authors considered these three cases together and study their effects in a strongly coupled anisotropic cardiac electro-mechanical model. X. Jie et al. \cite{XVN} used a electromechanical model for the ventricle of rabbit to analyze the arrhythmia originates from the ischemia induced electrophysiological and the mechanical changes. 
As per our literature survey no one has considered the electro-mechanical model of the human cardiac tissue to analyze the local ischemia effects on the electro-mechanical activity.

In this study, we will discuss the influence of local ischemic region on electrical and mechanical activity in the deforming human heart. For the electrical activity at cell level, we will use the ten Tusscher human ventricular membrane model \cite{ten}. 
We will consider the above three cases (i), (ii), (iii) and study the effect in deformed Monodomain model. We will also show the effect of ischemia on this model.
We will discuss the two types of ischemia namely Hyperkalemia and Hypoxia. We will discuss the influence of these two types of local ischemia, by changing the corresponding parameters in the 2D 
tissue level model, on the electrical and mechanical activity in terms of the action potential (AP), activation time (AT), repolarization time (RT), action potential duration (APD) and intracellular calcium ion concentration $[Ca^{+2}]_i$, active tension, ($T_A$), stretch ($\lambda$), stretch rate ($ \frac{d \lambda}{dt}$). We will also see the effect of the spread of ischemic zone in both the ischemic subregion and the healthy region on the electrical and mechanical properties. 

In the next section, we will describe the cardiac electro-mechanical model and also discuss about the ischemic parameters. We will also present the modeling of local ischemic zone. In the next section, we will present the fully discrete system using finite elements in space and backward Euler method in time. In section 4, we will discuss the numerical results. 

\section{Electro-mechanical models of cardiac tissue}
Cardiac electro-mechanical model is a combination of electrical and mechanical models. In this section, we will describe the 
mechanical and electrical models of cardiac tissue.

\subsection{Cardiac tissue mechanical model}
Let $\hat{\Omega}$ is the undeformed cardiac domain and $\Omega (t)$ is the deformed cardiac domain at time $'t'$ with 
$X=(X_1,X_2)$ and $x=(x_1,x_2)$ are the spatial coordinates respectively.

Define the deformation map, 
$\psi_t(X) = \psi(X,t) = x $ from $\hat{\Omega}$ to $\Omega (t)$ and $u(X,t)=x-X$ as the displacement vector.


In this work, cardiac tissue is considered as a nonlinear elastic material.
\newline
The deformation gradient tensor$F(X,t)=\{F_{i,j}\} = \{\frac{\partial x_i }{\partial X_j}, i,j=1,2,3 \}$,
\newline
$J(X,t)= detF(X,t)$
\newline
Cauchy-Green deformation tensor $C=F^TF$,
\newline
Lagrange-Green strain tensor $\epsilon=\frac{1}{2}(C-1)$.
\newline
Also the deformed body satisfy the following steady state force equilibrium equation described as 
\begin{align}
	\label{mech}
	div(FS(u,X))=0 , X \in \hat \Omega,
\end{align}
where $S$ is the second Piola-Kirchoff stress tensor. 
Prescribed displacement is imposed on the Dirichlet boundary $x(X,t) = \hat x (X), X \in \hat{ \partial \Omega_d}$ and Neumann 
boundary with no traction force $n^t FS(u(X,t),X) =0 , X \in \hat{ \partial \Omega_n}$.

The second Piola-Kirchoff stress tensor $S$, as taken in the studies \cite{hunter,kerckhoffs,vetter,lucamech}, is the combination
of three components namely, 
\begin{align*}
	S= S^{A} + S^{P} + S^{V},
\end{align*}
where active component $S^{A}$ generated biochemically, passive elastic component $S^{P}$ and volumetric component $S^{V}$. 

Next, $S^{P,V}$ is defined as follows
\begin{align}
	S_{ij}^{P,V} = \frac{1}{2} \bigg(\frac{\partial W^{P,V}}{\partial \epsilon_{ij}} + \frac{\partial W^{P,V}}{\partial \epsilon_{ji}} \bigg) , i,j =1,2,
\end{align}

where $W^{P,V}$ are the passive and volumetric strain energy functions, variety of these has been proposed in 
\cite{costa,gucci,holz,kuhl,nash,reme,schmid,usyk} and $\epsilon$ is the Green-Lagrange strain.

In this work, cardiac tissue is modeled as an orthotropic hyperelastic material. The passive strain energy function for a hyperelastic material is given as in \cite{TAGG}:
\begin{align*}
	&W^{P} = \frac{1}{2}(e^{Q}-1),\\ \nonumber
	&Q= \sum_{i=l,t} b_i (I_{4i}-1)^2 + b_{lt}I_{8lt}^2,
\end{align*}

where $I_{4l}=\hat b_l^T C \hat b_l$, $I_{4t}=\hat b_t^T C \hat b_t$ and
$I_{8lt}=\hat b_l^T C \hat b_t$,

and volumetric strain energy function for a hyperelastic material is given as
\begin{align*}
	W^V = K (J-1)^2,
\end{align*}
where $K$ is a positive bulk modulus. 

Now, the active component of second Piola-Kirchoff stress tensor $S^A$ is defined in the form of active tension originated along 
the myocytes.

An isometric contraction of the cardiac tissue, generates an active tension without changing length of the cardiac tissue. 
Mechanical model of active tension is modeled by the myofilaments dynamics activated by calcium. We consider, as suggested
in \cite{lucamech}, that the active force moves only in the fiber direction. Therefore, active Cauchy stress is given as
\begin{align}
	\sigma^A(x,t) = J^{-1}T_A a_l(x) \otimes a_l(x),
\end{align}

where $a_l=\frac{F \hat a_l}{|F \hat a_l|}$ is a unit vector parallel to the local fiber direction and $T_A$ is the active  tension in the deformed domain.

The active component of second Piola-Kirchoff stress tensor $S^A$ is defined as
\begin{align}
	S^{A}= JF^{-1} \sigma^A F^{-T},
\end{align}
and the stretch rate along the fiber is defined as
\begin{align}
	\lambda = \sqrt {\hat a_l^T C \hat a_l}.
\end{align}

The active tension is dependent on calcium dynamics, stretch and stretch rate i.e. $T_A =T_A(Ca_i, \lambda , \frac{d \lambda}{dt})$, as discussed in \cite{land} is given by the system of ODEs:

\begin{align}
	\label{tension}
	\nonumber
	& \frac{d(trpn)}{dt} = k_{trpn} \bigg( \bigg (\frac{Ca_i}{Ca_{50}(1+ \beta (\lambda -1))} \bigg )^{n_{trpn}} \bigg), \\ \nonumber
	& \frac{dXb}{dt} = k_{Xb} \bigg( (trpn_{50}) (trpn)^{n_{Xb}} (1-Xb) - \frac{1}{(trpn_{50}) (trpn)^{n_{Xb}}} Xb  \bigg), \\ \nonumber
	&\frac{dq_i}{dt} = A_i \frac{d\lambda}{dt} - \alpha_i q_i, i=1,2,\\
	& T_A=g(q) h(\lambda) Xb , q=q_1+q_2.
\end{align}
where, function $g(q)$ is the effect on tension is defined as :

\[  g(q)= \left\{ 
\begin{array}{ll}
	\frac{aq+1}{1-q} & q\leq 0 \\
	\frac{1+(a+2)Q}{1+q} & q> 0 \\
\end{array} 
\right. \]. 

For more details refer to \cite{land}.

\subsection{Cardiac tissue electrical model on a deforming domain: Monodomain model}
We assume that the cardiac tissue is insulated. The influence of the cardiac tissue deformation on the Monodomain model in 
the strong coupling framework is due to three different mechano-electrical feedback: a) presence of deformation gradient
$F$ in the conductivity coefficients structure, b) presence of deformation gradient $F$ and the deformation rate in the
conductive term, and c) presence of stretch $\lambda$ in the $I_{ion}$.

Monodomain model on the deformed domain $\Omega(t)$ is defined as follows:
\begin{align}
	\label{mvdefor}
	& \frac{\partial v}{\partial t}- div(D(x)\nabla v) + I_{ion}(v,w,c,\lambda) = I_{app},    x \in \Omega (t),  0\leq t \leq T\\
	\label{mwdefor}
	&\frac{\partial w}{\partial t}-g(v,w)= 0, x \in \Omega (t),  0\leq t \leq T\\
	\label{mcdefor}
	&\frac{\partial c}{\partial t}-f(v,w,c)= 0,  x \in \Omega (t),  0\leq t \leq T\\
	&v(x,0)= v_0(x,0), \hspace{5mm} w(x,0)=w_0(x,0), \hspace{5mm} c(x,0)=c_0(x,0), x \in \Omega\\
	&n^T D(x) \nabla v =0,  x \in \partial \Omega, 0\leq t \leq T
\end{align}
where,
\begin{align*}
	D(x)= {\sigma_l} a_l(x) {a_l}^T(x)+{\sigma_t} a_t(x) {a_t}^T(x),
\end{align*}
where, $a_l(x)$ is a unit vector along the local fiber direction and $a_t(x)$ is a unit vector transverse to the fiber direction and they both are orthogonal to each other. $\sigma_l$ and  $\sigma_t$ are the conductivity coefficients in the corresponding directions.
\begin{align*}
	D(x)= {\sigma_t} I +({\sigma_l} -{\sigma_t} )a_l(x) {a_l}^T(x). 
\end{align*}

Let $\hat v(X,t), \hat w(X,t), \hat c(X,t)$ and $\hat D(X,t)$ are the action potential, gating, concentration 
variables on the undeformed domain $\hat \Omega$. Then

\begin{align*}
	&\frac{\partial \hat v}{\partial t}= \frac{\partial v} {\partial t} + \nabla v . \frac{\partial \psi_t }
	{\partial t}, (\text{convective term})\\
	& \nabla v(x,t)= F^{-T} Grad \hat v(X,t), \\
	&dx=JdX,
\end{align*}
Thus the monodomain system on the reference domain $\hat{\Omega}$ is given as follows:
\begin{align}
	\label{mvr}
	\nonumber
	J \bigg( \frac{\partial \hat v}{\partial t}- F^{-T} \nabla \hat v . \frac{d\psi_t}{dt} \bigg )- div(J F^{-1} D F^{-T}\nabla \hat v) + J I_{ion}(\hat v, \hat w, \hat c,\lambda)  = J I_{app}, \hspace{2em}  X \in \hat{\Omega},
\end{align}
coupled with the system of ODEs for the gating and concentration variables $w(x,t)$ and $c(x,t)$ respectively, given by
\begin{align}
	&\frac{\partial w}{\partial t}-g(v,w)= 0, &  X \in \hat \Omega, T\\
	&\frac{\partial c}{\partial t}-f(v,w,c)= 0, &  X \in \hat \Omega.
\end{align}
Initial conditions:
\begin{align}
	\nonumber
	& \hat v(X,0)= \hat v_0(X), \hspace{2mm} X \in \hat \Omega, \\&  w(X,0)=w_0(X), \hspace{5mm} c(X,0)=c_0(X), \hspace{1mm}  X \in \hat \Omega.
\end{align}
Boundary condition:
\begin{align}
	n^T F^{-1} D(x) F^{-T} \nabla \hat v =0, \hspace{5mm}  \text{on} \hspace{5mm} \partial \hat \Omega \times (0,T).
\end{align}
Now, we compute the conductivity tensor $\hat D(X,t)$ on the reference domain is given as:
\begin{align*}
	\hat D(X,t)= F^{-1} D(x(X,t)) F^{-T}(X,t), \hspace{8mm} X \in \hat{\Omega}
\end{align*}
We know the conductivity tensor in the deformed domain $\Omega (t)$ is defined as
\begin{align*}
	D(x)= {\sigma_t} I +({\sigma_l} -{\sigma_t} )a_l(x) {a_l}^T(x).
\end{align*}
Now, let $\hat a_l(X)$ is a unit vector parallel to the local fiber direction in the reference domain, then in the deformed domain the unit vector $a_l(x)$ along the local fiber direction is obtained by 
\begin{align*}
	a_l &= \frac{F \hat a_l}{\|F \hat a_l\|},\\
	&= \frac{F \hat a_l}{\sqrt{\hat a_l^T C \hat a_l}}
\end{align*}
Thus,
\begin{align*}
	a_l a_l^T=  \frac{F \hat a_l \hat a_l^T F^T}{{\hat a_l^T C \hat a_l}}
\end{align*}

\subsection{Ionic model}
$I_{ion}(v, w, c,\lambda)$, $g(v,w)$ and $f(v,w,c)$ are defined by the cell level model. In this work we are using the human
ventricular cell level model, ten Tusscher model \cite{ten} (TT06). Dynamics of the ionic currents is modeled by the 17 ODEs. The ionic
current flow $I_{sac}$ due to the stretch activated channels and $ K_{ATP}$ current is also taken into account. Thus the total ionic current is given by
\begin{align}
	\nonumber
	& I_{ion}= I_{Na} + I_{K1} + I_{to}+ I_{Kr}+ I_{Ks}+ I_{CaL}+ I_{NaCa}+ I_{NaK} +\\
	& \hspace{0.6cm} I_{pCa}+ I_{pK}+ I_{bCa}+ I_{bNa} + I_{KATP}+ I_{sac},
\end{align}

where, the fast sodium $(I_{Na})$, L-type calcium$(I_{CaL})$, 
transient outward$(I_{to})$ , rapid delayed rectifier$(I_{Kr})$ and slow delayed rectifier $(I_{Ks})$, and inward rectifier 
currents $(I_{K1})$, $I_{NaCa}$ is $Na/Ca^{2+}$ exchange current, $I_{NaK}$ is $Na/K^{+}$ pump current, $I_{pCa}$ and $I_{pK}$ are
plateau $Ca^{2+}$ and $K^+$ currents, $I_{bCa}$ and $I_{bK}$ are background
$Ca^{2+}$ and $K^+$ currents and $I_{KATP}$ is the$ K_{ATP}$ current, all are modeled as in \cite{ten}, $I_{KATP}$ which accounts for
Hypoxia is defined as in \cite{panfilov}:
\begin{align}
	I_{KATP}= G_{KATP} f_{ATP} \bigg (\frac{[K^+]_o}{5.4}\bigg )^{0.3} \frac{1}{40+3.5e^{0.025v}}(v-E_k),
\end{align}
where, $G_{KATP}$ is the conductance and $E_k$ is the resting potential of potassium.

$I_{sac}$, the stretch activated current \cite{hancox} is given as follows : 
\begin{align*}
	&I_{sac}= I_{sac, n} + I_{Ko},\\
	&I_{sac,n} = I_{sac, Na} + I_{sac, K},
\end{align*}
where, $I_{sac,n}$ is the non-specific current defined as follows:
\begin{align*}
	&I_{sac, Na}= G_{sac} \gamma_{slsac}(v-v_{Na}) \bigg( - \frac{v_R+85}{v_R-65} \bigg), \\
	&	I_{sac, K}= G_{sac} \gamma_{slsac}(v-v_{K}).
\end{align*}

The specific $K^+$ stretch-dependent current is given by
\begin{align*}
	I_{Ko}= G_{Ko} \frac{\gamma_{sl,Ko}}{1+exp(-(10+v)/45)} (v-v_K),
\end{align*}
where, $\gamma_{sl,Ko}=3(\lambda-1)+0.7$.

The two types of ischemia, Hyperkalemia and Hypoxia will be incorporated into the model by varying the parameters corresponding to each of the two types of ischemia.

\begin{figure}
	\vspace{-8em}
	\centering
	\includegraphics[width=0.8\textwidth]{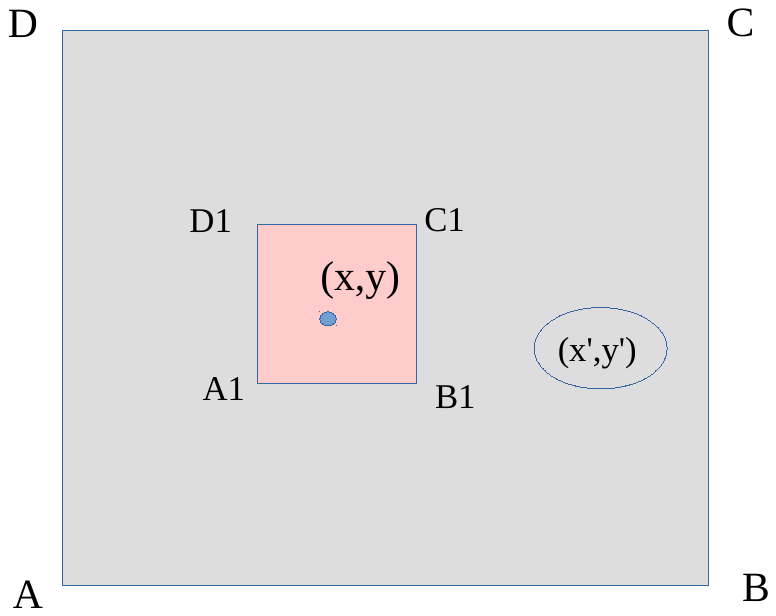}
	\vspace{-35em}
	\caption{: \textbf{Cardiac tissue $ABCD$ with single ischemic subregion $A1B1C1D1$, $(x,y)$ and $(x',y')$ are the points inside and outside the ischemic subregion.}}
	\label{1region_name}
\end{figure}

Human Cardiac tissue is considered as a square $A B C D$, and the ischemic region is taken as a small square $A_1 B_1 C_1 D_1$.
In this study, local ischemia is modeled by varying the ischemic parameter value in the sub-domain $A_1 B_1 C_1 D_1$ only and using
the normal value of the parameter in the remaining part of the cardiac tissue $ABCD$. Also considering the ischemia as mild, moderate
and severe by taking the $[K^+] _o$ parameter value in the corresponding range in the domain $A_1 B_1 C_1 D_1$ only. Then, we compare
the electrophysiological properties of the ischemic and non-ischemic region points say $(x,y)$ and $(x',y')$ respectively
(see Fig.\eqref{1region_name}).

\section{Weak Formulation}
First of all we define the strong form the model.
Monodomain model on the deformed domain $\Omega(t)$ is defined as follows:
\begin{align}
	\label{elast}
	\text{(elasticity equation)} \hspace{2mm} &- \nabla . (\nabla u \sigma(x, \gamma))=f,&  x \in \Omega(t), a.e. t \in (0,T)\\
	\label{mv}
	\text{(Tissue level model)} \hspace{2mm} & \frac{\partial v}{\partial t}- div(D(x)\nabla v) + I_{ion}(v,w,c,\lambda) = I_{app},&  x \in \Omega (t),  0\leq t \leq T\\
	\label{mw}
	\text{(cell level model)} \hspace{2mm} &\frac{\partial w}{\partial t}=G(v,w),& x \in \Omega(t),  0\leq t \leq T\\
	\label{mc}
	&\frac{\partial c}{\partial t}=h(v,w,c),& x \in \Omega(t),  0\leq t \leq T\\
	&\frac{\partial \gamma}{\partial t}=S(\gamma,w),&  x \in \Omega(t),  0\leq t \leq T\\
	&v(x,0)= v_0(x,0), \hspace{2mm} w(x,0)=w_0(x,0), \hspace{2mm} c(x,0)=c_0(x,0),& x \in \Omega(t)\\
	\label{bndrydeform}
	&n^T D(x) \nabla v =0,&  x \in \partial \Omega, 0\leq t \leq T.
\end{align}
Now, we will define the general form of the ionic model so that the weak formulation of the above model equations can be written.

The general form of ionic current $I_{ion} : \mathbb{R} \times \mathbb{R}^d \times (0, \infty )^n \to \mathbb{R}$ is defined as follows:
\begin{align}
	\label{Ion1mech}
	I_{ion}(v,w,c) := \sum_{i=1}^n Q_i(v,w,logc_i) + R(v,w,c), 
\end{align}
where, 
\begin{align}
	\label{Ionpropertymech}
	\nonumber
	& Q_i \in C^1(\mathbb{R} \times \mathbb{R}^d \times \mathbb{R}),\\ \nonumber
	& 0 < \underline{\rm H}(w) \leq \frac{\partial}{\partial \zeta} Q_i(v,w,\zeta) \leq \overline{\rm H}(w),\\
	& |\frac{\partial}{\partial v} Q_i(v,w,0)| \leq  L(w),
\end{align}
$\underline{\rm H}, \overline{\rm H}, L \in C^0(\mathbb{R}^d, \mathbb{R}_{+})$ and  
\begin{align}
	\label{Ionproperty1mech}
	R \in C^0(\mathbb{R} \times \mathbb{R}^d \times (0, \infty )^n ) \cap Lip(\mathbb{R} \times [0,1]^d \times
	(0,+ \infty)^m).
\end{align}
The following system of ODEs describing the dynamics of gating variables,
\begin{align}
	\label{wODEmech}
	\frac{\partial w_i}{\partial t} = G_i(v,w_i), \hspace{5mm} i=1,2,...,d.
\end{align}
where, $G_i$ has the particular form as
\begin{align*}
	G_i(v,w) := \alpha_i(v)(1-w_i)-\beta_i(v)w_i, \hspace{5mm} i=1,2,...,d,
\end{align*}
and $\alpha_i$ and $\beta_i$ are positive rational exponential functions in $v$.

Assumption, $\forall \hspace{2mm} i=1,2,...,d$,
\begin{align}
	\label{wpropertymech}
	\nonumber
	&G_i : \mathbb{R}^2 \to \mathbb{R}\hspace{5mm}  \textbf{is locally Lipschitz continuous},\\ \nonumber
	& G_i(v,0) \geq 0 \hspace{5mm} \forall v \in \mathbb{R},\\
	&G_i(v,1) \leq 0 \hspace{5mm} \forall v \in \mathbb{R}.
\end{align}
ODEs describing the dynamics of concentration variables is given by
\begin{align}
	\label{cODEmech}
	\frac{\partial c_j}{\partial t} = h_j(v,w,c_j):= Q_j(v,w,c_j)+ R_j(v,w,c_j)  , \hspace{5mm} j=1,2,...,n,
\end{align}
where ,
\begin{align}
	\label{cpropertymech}
	R_j \in C^0(\mathbb{R} \times \mathbb{R}^d \times (0, \infty )^n ) \cap Lip(\mathbb{R} \times [0,1]^k \times
	(0,+ \infty)^m), \hspace{5mm} j=1,2,...,n.
\end{align}
The function $S$ is given by $S(\gamma , w) = \beta (\sum_{j=1}^{k} \eta_j w_j -\eta_0 \gamma)$, where $\beta$ and $\eta_j$ are the positive pathological parameters.
Now, we will define the weak formulation of the above defined model.

\subsection{Weak formulation}
Definition: A weak solution of the problem \eqref{elast}-\eqref{bndrydeform} is a set $(u,v,w,c,\gamma)$ such that $u \in L^2 (0,T; H^1(\Omega)^3)$, 
$v \in  L^2 (0,T; H^1(\Omega))$, $v_t \in  L^2 (0,T; (H^1(\Omega))')$, $w \in C(0,T; L^2(\Omega)^m)$, $c \in C(0,T; L^2(\Omega)^n)$,
$z \in C(0,T; L^2(\Omega))$, and it satisfy the following:

$\forall \phi_1 \in  H^1(\Omega)^3, \phi _2 \in H^1(\Omega), \phi _3 \in L^2(\Omega)^m, \phi _4 \in L^2(\Omega)^n, \phi _5 \in L^2(\Omega)$ a.e. $t \in (0,T)$ such that:
\begin{align}
	\label{weak_u}
	&\int_{\Omega} \nabla u \sigma(x, \gamma):\nabla \phi_1 = \int_{\Omega} f \phi_1 dx + \int_{\partial \Omega} u \phi _1 dx, \\
	\label{weak_v}
	&\int_{\Omega} \frac{\partial v}{\partial t} \phi _2 + \int_{\Omega} D(x,\nabla u) \nabla v  \nabla \phi_2 
	+ \int_{\Omega} I_{ion}(v,w,c) \phi _2  = \int_{\Omega} I_{app} \phi _2,\\
	\label{weak_w}
	&\int_{\Omega} \frac{\partial w}{\partial t} \phi _4 = \int_{\Omega} G(v,w) \phi_3, \\
	\label{weak_c}
	&\int_{\Omega} \frac{\partial c}{\partial t} \phi _4 = \int_{\Omega} h(v,w,c) \phi_4, \\
	\label{weak_gamma}
	&\int_{\Omega} \frac{\partial \gamma}{\partial t} \phi _5 = \int_{\Omega} S(\gamma,w) \phi_5. 
\end{align}
\textbf{Note}: $H^1(\varOmega)= \{u:u\in L^2(\varOmega), \frac{\partial u}{\partial x_i} \in L^2(\varOmega), i=1,2 \}$,\\
\textbf{Note}: $L^2(\Omega)=\{u : \int_{\varOmega} u ^2 < \infty \}$.

\section{Numerical Scheme}
The above electro-mechanical model is solved via finite element method for space discretization and implicit-explicit (IE) Euler method for the time discretization.

Cardiac tissue is discretized with rectangular grid $T_{h_e}$ (grid size $h_e$) for the electrical monodomain model \eqref{mvr} and  $T_{h_m}$ (grid size $h_m$) for the mechanical model \eqref{mech}, where $T_{h_e}$ is a refinement of $T_{h_m}$. We will then use quadrilateral finite elements in space for both the electrical and mechanical models. 
For time discretization of both electrical and mechanical models we will use equal time steps $dt$. 

We will solve this coupled electro-mechanical model in two steps: 

1) Time discretization of the ODEs to find $w$, $c$ then solve the mechanical model to find the deformed system.

2) Solve the electrical monodomain model to compute $v$.

1) a) Compute $w^{n+1}, c^{n+1}$, for a given $v^n, w^{n}, c^{n}$, by solving 
\eqref{mwdefor} and \eqref{mcdefor} via Implicit- explicit Euler method:
\begin{align}
	&w^{n+1}= w^n+ dt G(v^n, w^{n+1}),\\
	&c^{n+1}= c^n+ dt f(v^n, w^{n+1}, c^n);
\end{align}

b) Now, use the obtained $w^{n+1}, c^{n+1}$ to solve the mechanical model \eqref{mech}- \eqref{tension} and compute the new deformed 
coordinates $x^{n+1}$ in terms of the deformation gradient tensor $F^{n+1}$.

2)Solve the electrical monodomain model using the linear finite elements in space and the IE Euler method in time, for a given $w^{n+1}, c^{n+1}$ and $F^{n+1}$ and compute the action potential $v^{n+1}$.

Define 
\begin{align*}
	&V_h = \{ v \in H^1(\Omega): v \hspace{0.2mm} \text{is continuous in} \hspace{0.2mm} \Omega : {v}_{\mid E} \in Q_1(E) , \forall E \in T_{h_e} \},\\ 
	&W_h = \{ u \in H^1(\Omega)^3: u \hspace{0.2mm} \text{is continuous in} \hspace{0.2mm} \Omega : {u}_{\mid E} \in Q_1(E) ,\forall E \in T_{h_m}\}.
\end{align*}

After solving the ODEs, solve the mechanical equation using the FEM (mechanical discretization, $T_{h_m}$) and calculate the deformation gradient at the next level.
%

Choose $\{\xi_i\}_{i=1}^N$ and $\{\psi_i\}_{i=1}^N$ as the basis functions in the space $V_h$ and $W_h$ respectively. Then we can write the following:

\begin{align*}
	&u_n=\sum_{i=0}^n a_{u,i} \boldsymbol{ \psi}_n,\\
	&v_n=\sum_{i=0}^n a_{v,i} \xi_n, 
\end{align*}

The mechanical equation using the finite element method (mechanical discretization, $T_{h_m}$) becomes as follows:
\begin{align*}
	\int_{T_{h_m}} S \nabla \Big (\sum_{i=1}^{N}u_{i,t} \psi_i(x) \Big ) :\nabla \psi_j = 0.
\end{align*} 
and then solve it using the Newton's method. 

Also, solve the monodomain equation using the finite element method (electrical discretization, $T_{h_e}$) in space, the weak form of the monodomain equation becomes as follows:
\begin{align*}
	\int_{T_{h_e}}\Big (\sum_{i=1}^{N}v_{i,t} \xi_i(x) \Big ) \xi_j + \int_{T_{h_e}} D(x,\nabla u)\nabla \Big (\sum_{i=1}^{N}v_i\xi_i(x)\Big ) \nabla \xi_j= \int_{T_{h_e}} I_{ion}^h \xi_j,
\end{align*}

\begin{align*}
	\sum_{i=0}^n a'_{v,i} \int_{\Omega} \xi_i \xi_n + \sum_{i=0}^n a_{v,i} \int_{\Omega} D(x,\nabla u_{i,n}) \nabla \xi_i \nabla \xi_n 
	+ \int_{\Omega} I_{ion}(a_{v,i},a_{w,i},a_{c,i}) \xi_n  = \int_{\Omega} I_{app} \xi_n.
\end{align*}
This is an ODE in time which is now solved using the IE Euler method.

\section{Results and Discussion}
We will present the results computed from the numerical simulation performed via Matlab Ra2015 library. 
Cardiac tissue is taken as a two dimensional slab of size $[0,1]^2$. The conductivity tensor in two-dimension we are using is  as follows:
\begin{align*}
	D=\begin{bmatrix} \sigma_l  & 0 &  \\ 0 & \sigma_t \end{bmatrix}
\end{align*} 
We need to choose an appropriate grid system so that the numerical solution is computed up-to acceptable accuracy. Such a goal is achieved through grid validation
test. Prior to the grid validation we will briefly discuss the discontinuity treatment of
$I_{ion}$ expression due to ischemic subregions inside a healthy cardiac tissue. While
Hyperkalemic ischemic zones results in the discontinuity of $[K^+]_o$ parameter at the
interface of healthy and ischemic parts of a cardiac tissue the Hypoxic ischemic zones
results in the discontinuity of $f_{ATP}$ parameter at the interface. Here these discontinuities
are treated to a desired order of accuracy by a simple regularization effect through
appropriate interpolating polynomials. These local interpolating polynomials of high
order accuracy are constructed following the procedure discussed by Alsoudani and
Bogle in \cite{bogel} using refined localized data around the discontinuity. Different grid systems
consisting of (a) 441 (b) 625 (c) 841 (d) 1089 and (e) 2601 degrees of freedom
(dofs) have been considered. We compare the action potential obtained using these grid systems (as in \cite{pargaei2019}). In all these cases only a marginal variation (less than 2–3\% for 441 and 625 grids and 0.5\% for others) is found.

The parameters value in the model is taken from the reference \cite{lucamech}. Applied stimulus $I_{app}$ is taken as 200 $mA/cm^3$ to start the depolarization process \cite{lucamech}. The initial condition for the action potential and all the gating variables of the TT06 ionic model are taken as the corresponding resting values while the boundary of the domain is assumed to be insulated.
\newline
\textbf{Activation time (AT)} is defined as the unique time $\tau_a(x)$ during the upstroke phase of the action potential 
when $v(x,\tau_a(x))=-50$. \textbf{Repolarization time (RT)} is defined as the unique time $\tau_r(x)$ during the repolarization
phase of the action potential when $v(x,\tau_r(x))=0.9v_R$. \textbf{Action potential duration (APD)} is defined as $\tau_r(x)-\tau_a(x)$.

First of all, we will present the results for the simulation of the four situations S1, S2, S3 and S4 in the cardiac tissue which are,
(S1) without the mechanical feedback, (S2) consider the mechanical feedback into the diffusion term of the Monodomain model but
without $I_{sac}$ current and the convective term, (S3) Take the mechanical feedback in the diffusion term and the $I_{sac}$ current,
and, (S4)  Take the mechanical feedback in the diffusion term with $I_{sac}$ current and the convective term. We will also examine 
the impact on the temporal variation of the action potential $(v)$, intracellular calcium concentration 
($[Ca^{+2}]_i$), active tension $(T_A)$, stretch along the fiber $(\lambda)$ and stretch rate $(\frac{d\lambda}{dt})$ at the points (describe later) in the cardiac tissue with the change of Hyperkalemia and Hypoxia both. We will also examine the increased ischemic region size effect onto the non-ischemic region in terms of its electro-mechanical properties. 

\subsection{(Case 1) Cardiac tissue without ischemic subregion} 
\subsubsection{Effect of S2, S3 and S4 situations on electrical and mechanical waveforms}
\begin{figure*}
	\vspace{-0.5em}
	\centering
	\begin{subfigure}[t]{0.5\textwidth}
		\centering
		\includegraphics[width=0.7\textwidth]{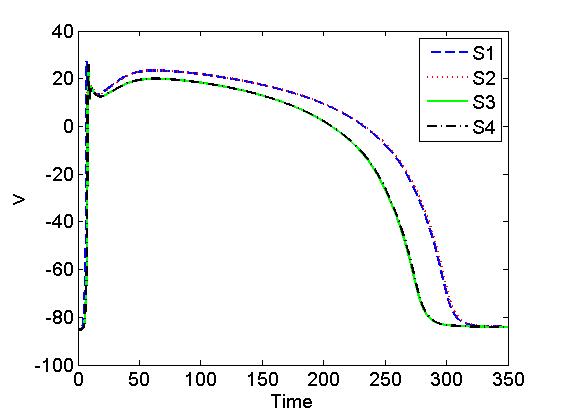}
		\vspace{-0.6em}
		\caption{}
		\label{AP_Ko_norm_diff_case_Pt1}
	\end{subfigure}\hfill	
	\hspace{-10em}
	\begin{subfigure}[t]{0.5\textwidth}
		\centering
		\includegraphics[width=0.7\textwidth]{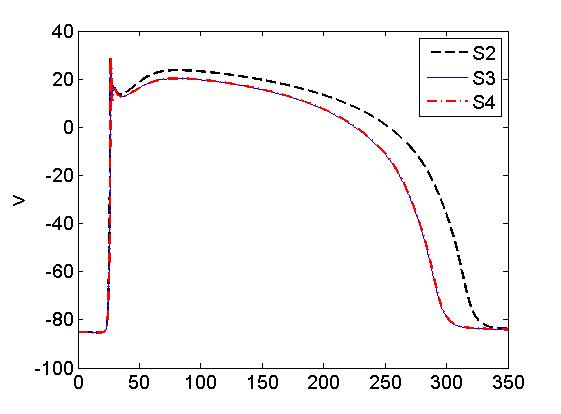}
		\vspace{-0.6em}
		\caption{}
		\label{AP_Ko_norm_diff_case_Pt3}
	\end{subfigure}\\
	\begin{subfigure}[t]{0.5\textwidth}
		\centering
		\includegraphics[width=0.7\textwidth]{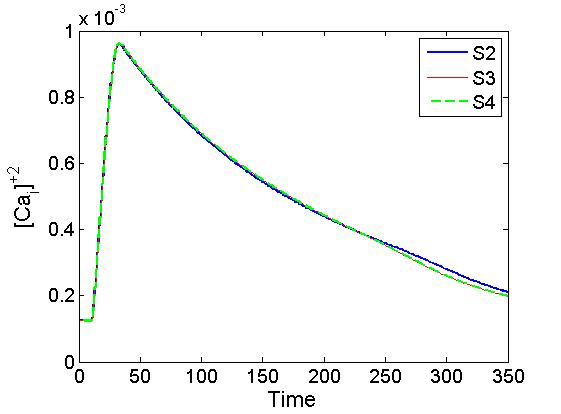}
		\vspace{-0.6em}
		\caption{}
		\label{Cai_Ko_norm_diff_case_Pt1}
	\end{subfigure}\hfill
	\hspace{-10em}
	\begin{subfigure}[t]{0.5\textwidth}
		\centering
		\includegraphics[width=0.7\textwidth]{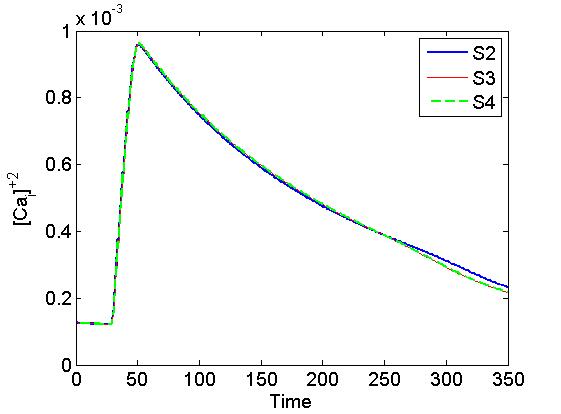}
		\vspace{-0.6em}
		\caption{}
	\end{subfigure}\\
	\begin{subfigure}[t]{0.5\textwidth}
		\centering
		\includegraphics[width=0.7\textwidth]{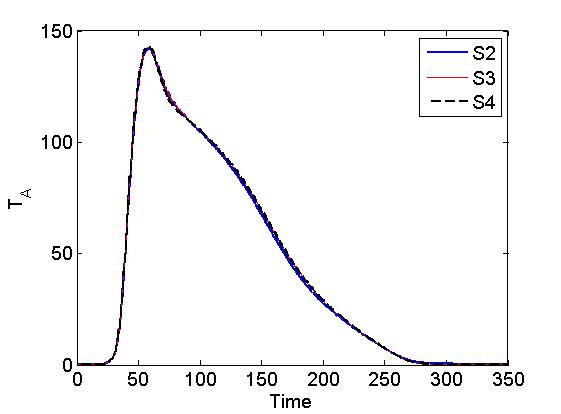}
		\vspace{-0.6em}
		\caption{}
	\end{subfigure}\hfill
	\hspace{-10em}
	\begin{subfigure}[t]{0.5\textwidth}
		\centering
		\includegraphics[width=0.7\textwidth]{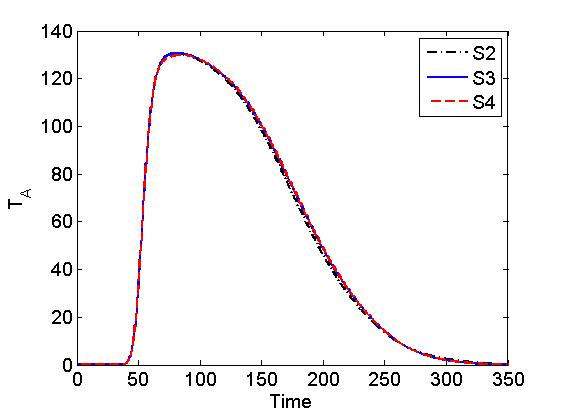}
		\vspace{-0.6em}
		\caption{}
	\end{subfigure}\\
	\vspace{-1em}
	\begin{subfigure}[t]{0.5\textwidth}
		\centering
		\includegraphics[width=0.7\textwidth]{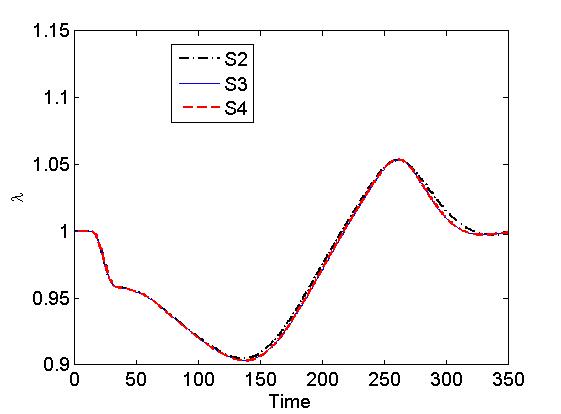}
		\vspace{-0.6em}
		\caption{}
		\label{lambda_Ko_norm_diff_case_Pt1}
	\end{subfigure}\hfill
	\hspace{-10em}
	\begin{subfigure}[t]{0.5\textwidth}
		\centering
		\includegraphics[width=0.7\textwidth]{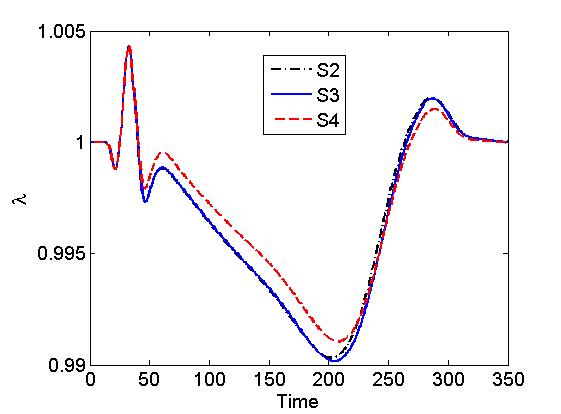}
		\vspace{-0.6em}
		\caption{}
		\label{lambda_Ko_norm_diff_case_Pt3}
	\end{subfigure}\\
	\begin{subfigure}[t]{0.5\textwidth}
		\centering
		\includegraphics[width=0.7\textwidth]{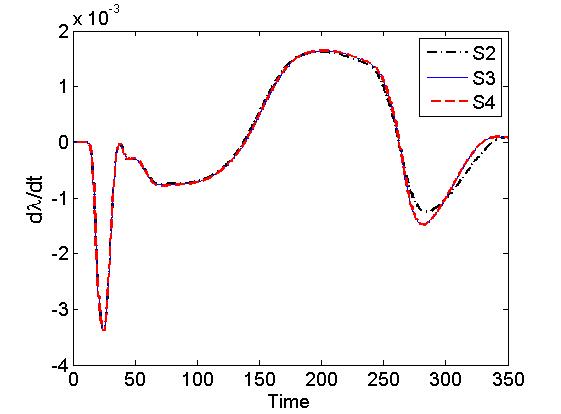}
		\vspace{-0.6em}
		\caption{}
		\label{dlambda_Ko_norm_diff_case_Pt1}	
	\end{subfigure}\hfill
	\hspace{-10em}
	\begin{subfigure}[t]{0.5\textwidth}
		\centering
		\includegraphics[width=0.7\textwidth]{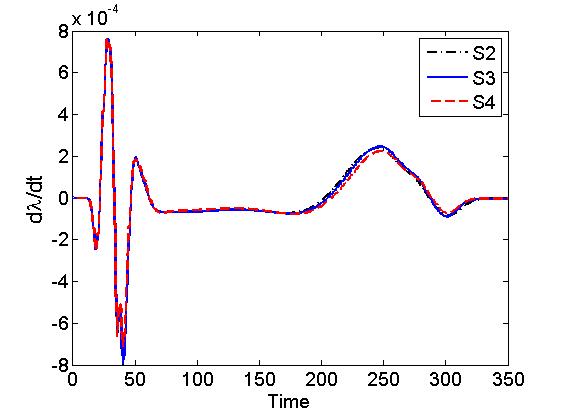}
		\vspace{-0.6em}
		\caption{}
		\label{dlambda_Ko_norm_diff_case_Pt3}
	\end{subfigure}\\
	\caption{: \textbf{$v$ (a), $[Ca^{+2}]_i$ (b), $T_A$ (c), $\lambda$ (d), $\frac{d \lambda}{dt} (e)$ for the cases, without MEF, without $I_{sac}$, 
			with $I_{sac}$ and with $I_{sac}+$ CONV at M1(column1) and M2(column2).}}
	\label{AP_MP_diff_case_Pt1}
\end{figure*}

We presented the action potential $(v)$, intracellular calcium concentration ($[Ca^{+2}]_i$), active tension $(T)$, stretch along the fiber $(\lambda)$ and stretch rate $(\frac{d\lambda}{dt})$ at the points (0.1875, 0.1875) (M1) and (0.5, 0.5) (M2) of the human cardiac tissue (see Fig. \eqref{AP_MP_diff_case_Pt1}, in the subsequent simulations M1 will be in the ischemic zone). The simulations are plotted in Fig \eqref{AP_MP_diff_case_Pt1}. From Fig. \eqref{AP_Ko_norm_diff_case_Pt1} and \eqref{AP_Ko_norm_diff_case_Pt3}, it can be examined that action potential$(v)$ for the situations S1 and S2 i.e. without mechanical feedback and the addition of mechanical feedback without $I_{sac}$ respectively, had a little difference. This implies that the addition of mechanical feedback had a little effect on the cardiac electrical activity in terms of ventricular action potential profile and hence on QRS complex of ECG.

Further, as we add the current generated due to the stretch activated channels i.e. $I_{sac}$, the plateau phase, repolarization phase and hence APD of action potential gets affected. Fast upstroke, early repolarization and hence $10\%$ decrease in APD is noticed. Next, with the addition of convective term i.e. S4, negligible difference in action potential profile in comparison to S2 is noticed.

Since, S1 is without mechanical case, therefore, we will consider the mechanical parameters for the situations S2, S3 and S4. Now, from Fig. \eqref{Cai_Ko_norm_diff_case_Pt1}- \eqref{dlambda_Ko_norm_diff_case_Pt3}, an alteration in the waveforms of the mechanical parameters is noticed with the addition of $I_{sac}$ at the points M1 and M2 of the human cardiac tissue. With the addition of the $I_{sac}$, waveform of $[Ca^{+2}]_i$ at the points M1 and M2 get affected during the resting phase and hence the contractile force and the contractility. While addition of convective term does not affect the contractility. 

Stretch along the fiber $(\lambda)$ and hence the stretch rate $\frac{d \lambda}{dt}$ changes with the addition of $I_{sac}$. From Fig. \eqref{lambda_Ko_norm_diff_case_Pt1} and \eqref{lambda_Ko_norm_diff_case_Pt3}, it is to be analyzed that during the systolic contraction there is shortening of fiber i.e. $\lambda < 1$ and after this, myocyte stretch $\lambda > 1$.  There is almost $1\%$ change in the length of myocytes at point M1 and  $3\%$ change in the length of myocytes at point M2 relative to the resting length.

Thus, it is concluded that ventricular contraction without $I_{sac}$ had little effect on the action potential and hence on ECG. With the addition of  $I_{sac}$ i.e. stretch regulated electrical activities alter the AP phase as repolarization phase and APD and hence alter the QRS complex and QT interval of ECG. The mechanical stretch induced current  $I_{sac}$ also changes the myocardial resting $[Ca^{+2}]_i$ which is related to the changes in contractile force and hence contractility. This is consistent with the experimental studies ().
Fiber length also changes due to the stretch activated channels, which is responsible for the contraction and expansion of the heart. In short, it is to be verified that the stretch activated channels plays an important role in the electro-mechanical activity of HCT or in the functioning of the heart.

\subsection{(Case 2) Cardiac tissue with ischemic subregion} 
\begin{figure*}
	\centering
	\vspace{-5em}
	\hspace{2em}
	\includegraphics[width=1.8\textwidth]{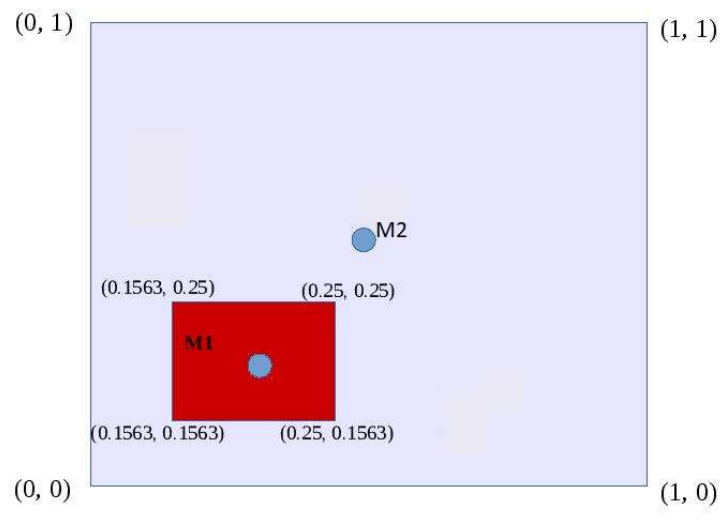}
	\vspace{-40em}
	\caption{: \textbf{Cardiac tissue with one ischemic subregion, $M1$ is the point inside the ischemic subregion and $M2$ is outside the ischemic subregion.}}
	\label{1region}
\end{figure*}

We consider the cardiac tissue of size $[0,1]^2$ with ischemic subregion as $[0.1563, 0.25]^2$. The ischemic parameters $[K^+]_o$ and $f_{ATP}$ are varied in the subregion only. The impact of the
change for these parameter values are analyzed on the ischemic region and on the healthy neighborhood of the ischemic region. In this case, we will examine the effect at points $M1 (0.1875, 0.1875)$ and $M2 (0.5, 0.5)$. $M1$ lies inside the ischemic zone of the cardiac tissue while $M2$ lies in the healthy neighborhood of this ischemic region (see Fig. \eqref{1region}). Now, we will examine the
electrophysiological and the mechanical properties of the cardiac tissue via action potential $(v)$, intracellular calcium concentration $[Ca^{+2}]_i$, active tension $T_A$, stretch along fiber $(\lambda)$ and stretch rate $(\frac{d \lambda}{dt})$ at these two points. 

\subsubsection{{(a) Hyperkalemia (Effect of ischemia when only $[K^+]_o$ varies)}}
\begin{figure*}
	\hspace{-5em}
	\begin{subfigure}[t]{0.5\textwidth}
		\centering
		\includegraphics[width=0.6\textwidth]{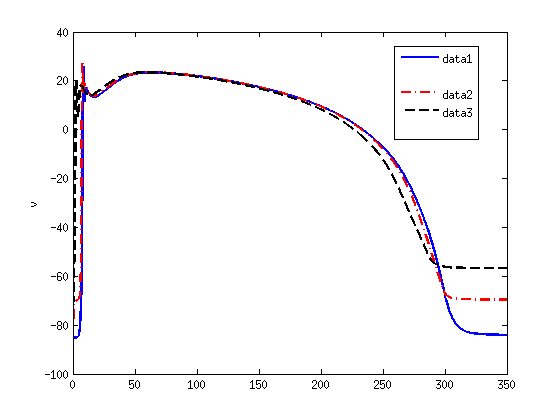}
		\caption{}
	\end{subfigure} %
	\hspace{-10em}
	\begin{subfigure}[t]{0.5\textwidth}
		\centering
		\includegraphics[width=0.6\textwidth]{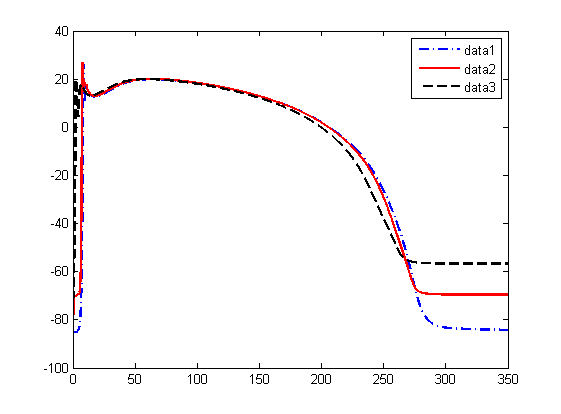}
		\caption{}
	\end{subfigure}\hspace{-10em}
	\begin{subfigure}[t]{0.5\textwidth}
		\centering
		\includegraphics[width=0.6\textwidth]{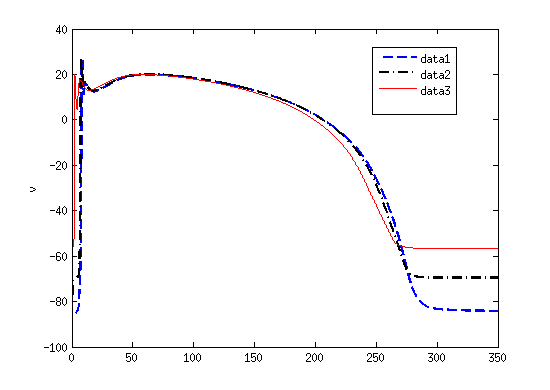}
		\caption{}
	\end{subfigure} \\
	\hspace{1em}
	\begin{subfigure}[t]{0.5\textwidth}
		\includegraphics[width=0.6\textwidth]{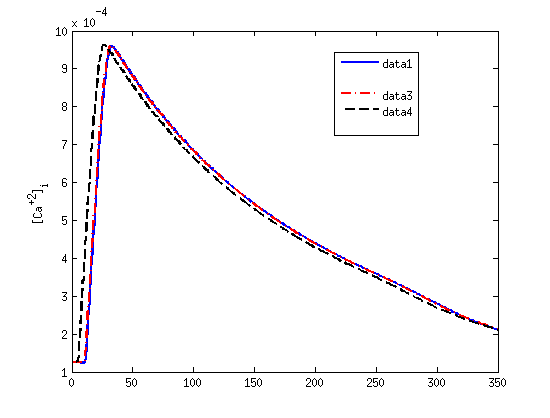}
		\caption{}
		\label{Cai_diff_Ko_no_ISAC_Pt1}
	\end{subfigure}\hspace{-15em}
	\begin{subfigure}[t]{0.5\textwidth}
		\centering
		\includegraphics[width=0.6\textwidth]{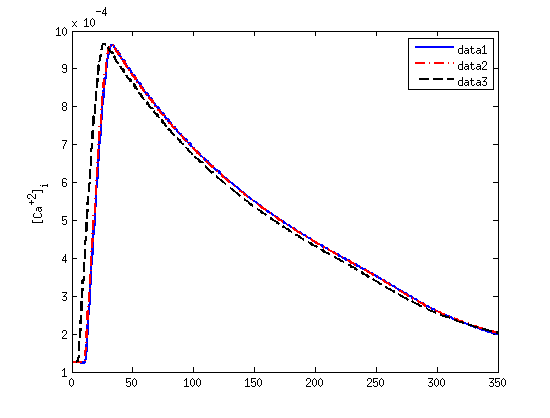}
		\caption{}
		\label{Cai_diff_Ko_with_ISAC_Pt1}
	\end{subfigure}\hspace{-10em}
	\begin{subfigure}[t]{0.5\textwidth}
		\centering
		\includegraphics[width=0.6\textwidth]{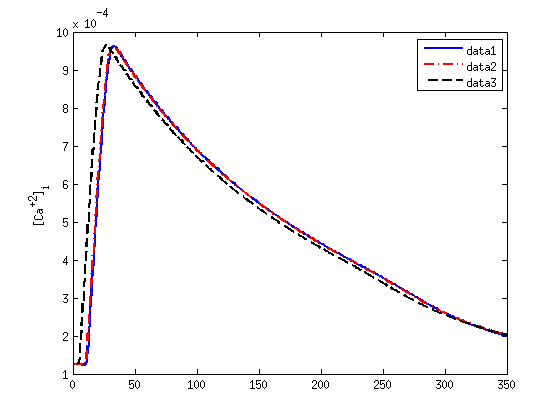}
		\caption{}
		\label{Cai_diff_Ko_with_ISAC_CONV_Pt1}
	\end{subfigure} \\
	\begin{subfigure}[t]{0.5\textwidth}
		\includegraphics[width=0.6\textwidth]{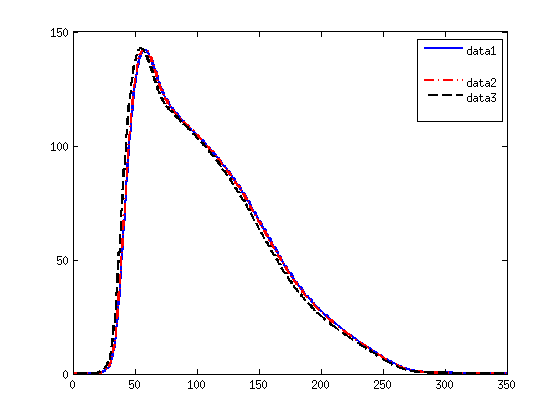}
		\caption{}
	\end{subfigure}\hspace{-15em}
	\begin{subfigure}[t]{0.5\textwidth}
		\centering
		\includegraphics[width=0.6\textwidth]{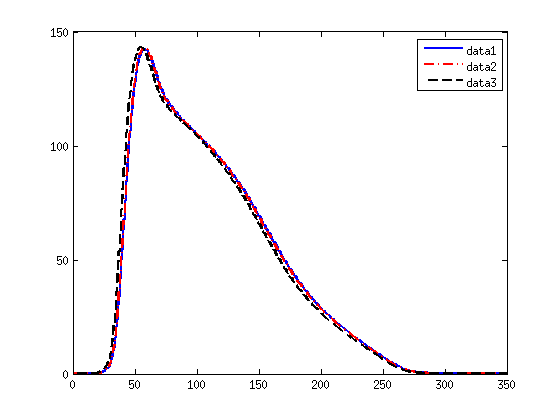}
		\caption{}
	\end{subfigure}\hspace{-10em}
	\begin{subfigure}[t]{0.5\textwidth}
		\centering
		\includegraphics[width=0.6\textwidth]{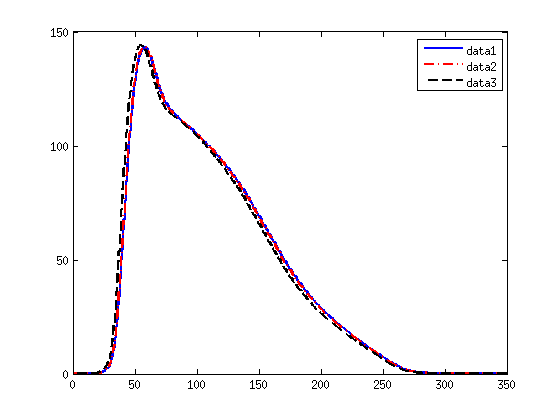}
		\caption{}
	\end{subfigure}\\
	\begin{subfigure}[t]{0.5\textwidth}
		\includegraphics[width=0.6\textwidth]{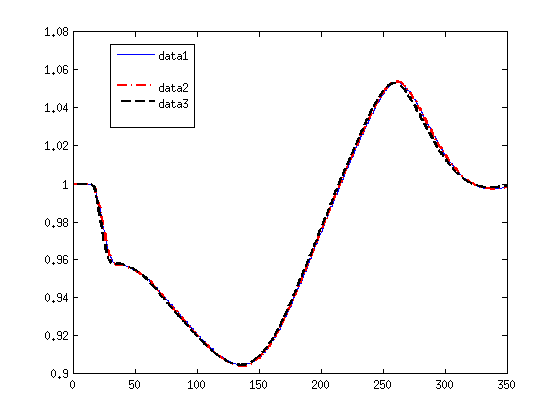}
		\caption{}
		\label{lambda_diff_Ko_no_ISAC_Pt1}
	\end{subfigure}\hspace{-15em}
	\begin{subfigure}[t]{0.5\textwidth}
		\centering
		\includegraphics[width=0.6\textwidth]{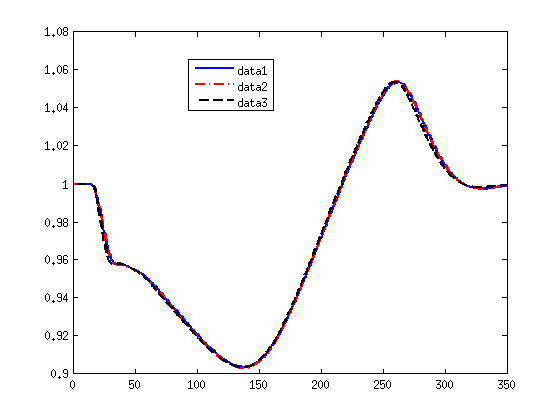}
		\caption{}
		\label{lambda_diff_Ko_with_ISAC_Pt1}
	\end{subfigure}\hspace{-10em}
	\begin{subfigure}[t]{0.5\textwidth}
		\centering
		\includegraphics[width=0.6\textwidth]{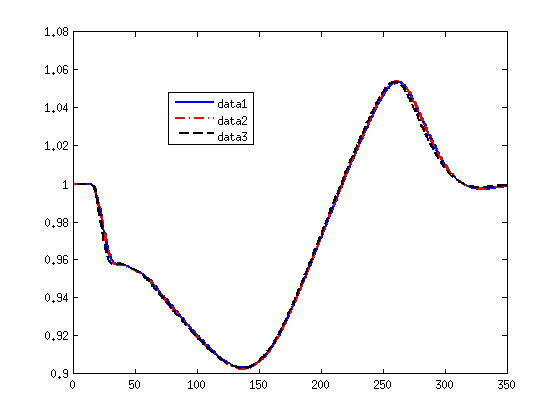}
		\caption{}
		\label{lambda_diff_Ko_with_ISAC_CONV_Pt1}
	\end{subfigure}
	\\
	\begin{subfigure}[t]{0.5\textwidth}
		\includegraphics[width=0.6\textwidth]{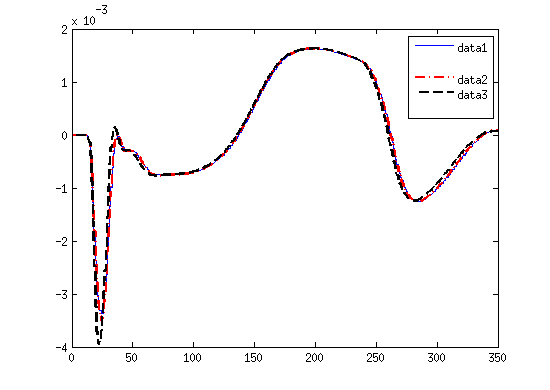}
		\caption{}
		\label{dlambda_diff_Ko_no_ISAC_Pt1}
	\end{subfigure}\hspace{-15em}
	\begin{subfigure}[t]{0.5\textwidth}
		\centering
		\includegraphics[width=0.6\textwidth]{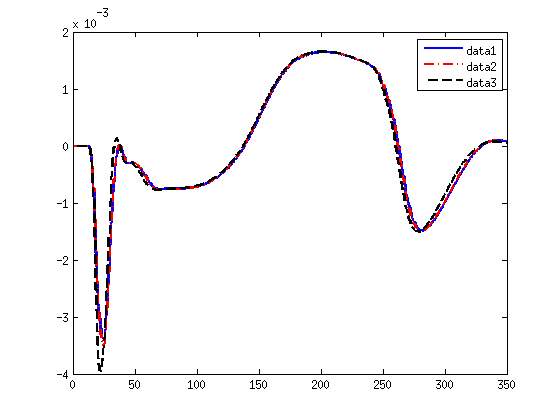}
		\caption{}
		\label{dlambda_diff_Ko_with_ISAC_Pt1}
	\end{subfigure}\hspace{-10em}
	\begin{subfigure}[t]{0.5\textwidth}
		\centering
		\includegraphics[width=0.6\textwidth]{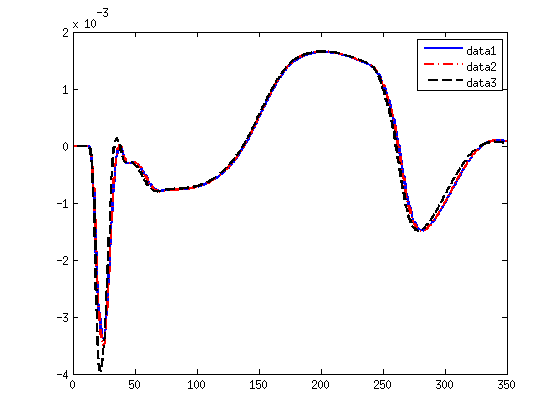}
		\caption{}
		\label{dlambda_diff_Ko_with_ISAC_CONV_Pt1}
	\end{subfigure}
	\caption{: \textbf{$v$ (first row), $[Ca^{+2}]_i$ (second row), $T_A$ (third row), $\lambda$ (fourth row), $\frac{d \lambda}{dt}$ (fifth row) for the cases without $I_{sac}$ (first column), with $I_{sac}$ (second column) and with $I_{sac}+$ CONV (fourth column) at M1, for different values of $[K^+]_o $, data1 ($[K^+]_o =5.4$),  data2 ($[K^+]_o =12$),  data3 ($[K^+]_o =20$)).}}
	\label{diffKo_diff_case_Pt1}
\end{figure*}

\begin{figure*}
	\hspace{-5em}
	\begin{subfigure}[t]{0.5\textwidth}
		\centering
		\includegraphics[width=0.6\textwidth]{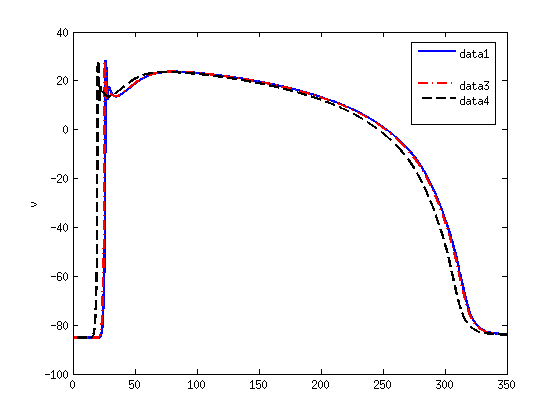}
		\caption{}
	\end{subfigure} %
	\hspace{-10em}
	\begin{subfigure}[t]{0.5\textwidth}
		\centering
		\includegraphics[width=0.6\textwidth]{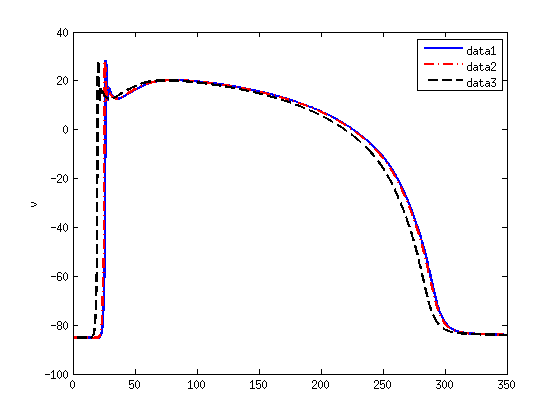}
		\caption{}
	\end{subfigure}\hspace{-10em}
	\begin{subfigure}[t]{0.5\textwidth}
		\centering
		\includegraphics[width=0.6\textwidth]{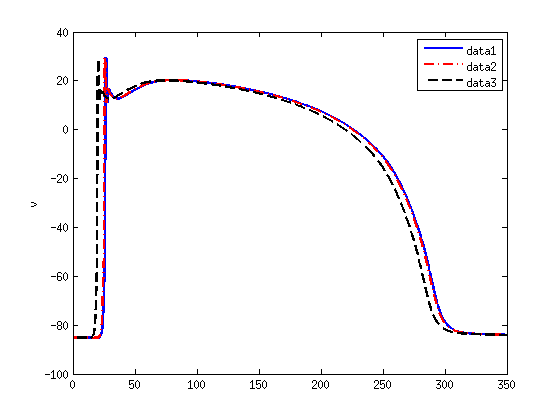}
		\caption{}
	\end{subfigure} \\
	\hspace{1em}
	\begin{subfigure}[t]{0.5\textwidth}
		\includegraphics[width=0.6\textwidth]{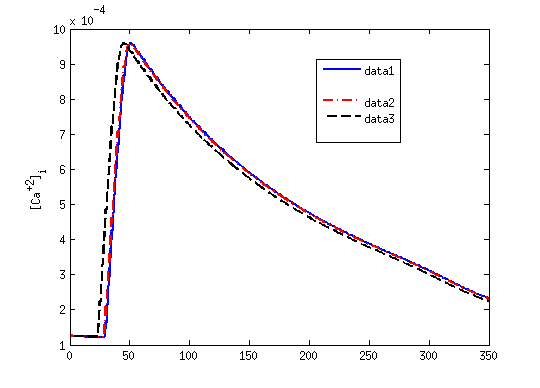}
		\caption{}
		\label{Cai_diff_Ko_no_ISAC_Pt3}
	\end{subfigure}\hspace{-15em}
	\begin{subfigure}[t]{0.5\textwidth}
		\centering
		\includegraphics[width=0.6\textwidth]{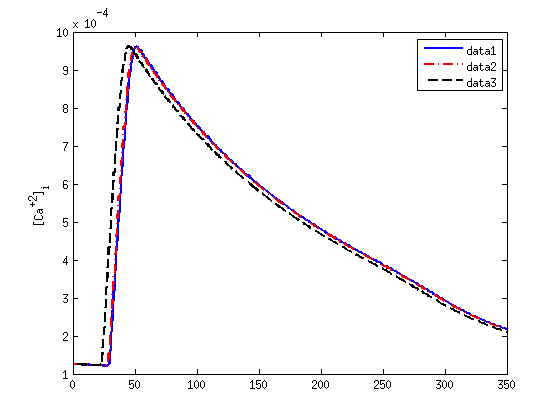}
		\caption{}
		\label{Cai_diff_Ko_with_ISAC_Pt3}
	\end{subfigure}\hspace{-10em}
	\begin{subfigure}[t]{0.5\textwidth}
		\centering
		\includegraphics[width=0.6\textwidth]{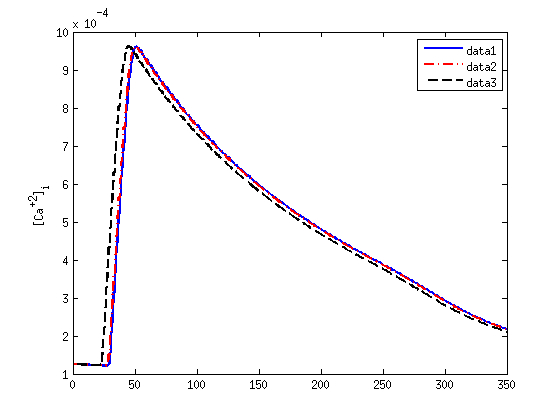}
		\caption{}
		\label{Cai_diff_Ko_with_ISAC_CONV_Pt3}
	\end{subfigure} \\
	\begin{subfigure}[t]{0.5\textwidth}
		\includegraphics[width=0.6\textwidth]{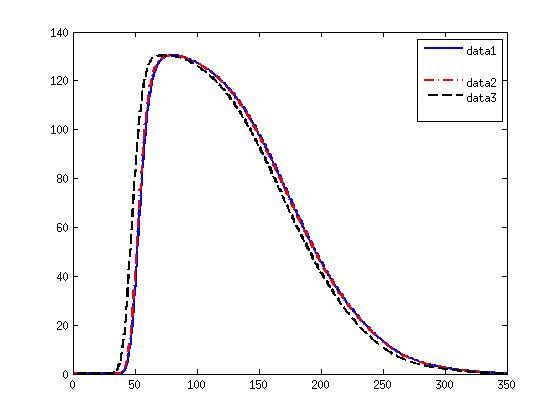}
		\caption{}
		\label{Ta_diff_Ko_no_ISAC_Pt3}
	\end{subfigure}\hspace{-15em}
	\begin{subfigure}[t]{0.5\textwidth}
		\centering
		\includegraphics[width=0.6\textwidth]{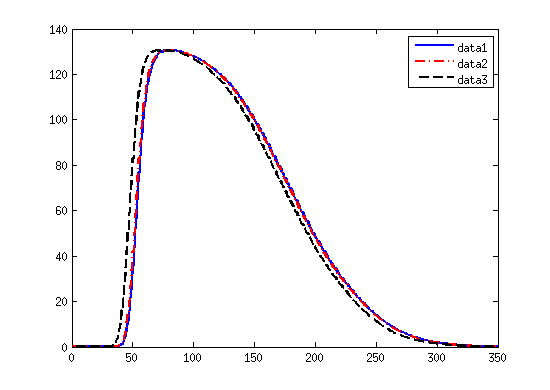}
		\caption{}
		\label{Ta_diff_Ko_with_ISAC_Pt3}
	\end{subfigure}\hspace{-10em}
	\begin{subfigure}[t]{0.5\textwidth}
		\centering
		\includegraphics[width=0.6\textwidth]{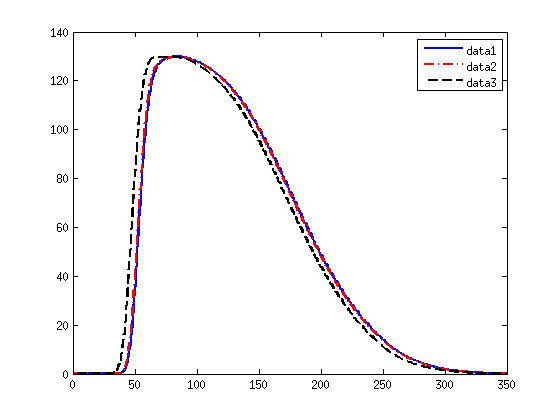}
		\caption{}
		\label{Ta_diff_Ko_with_ISAC-CONV_Pt3}
	\end{subfigure}\\
	\begin{subfigure}[t]{0.5\textwidth}
		\includegraphics[width=0.6\textwidth]{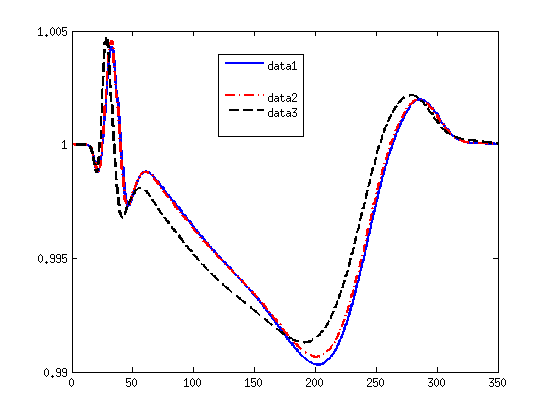}
		\caption{}
		\label{lambda_diff_Ko_no_ISAC_Pt3}
	\end{subfigure}\hspace{-15em}
	\begin{subfigure}[t]{0.5\textwidth}
		\centering
		\includegraphics[width=0.6\textwidth]{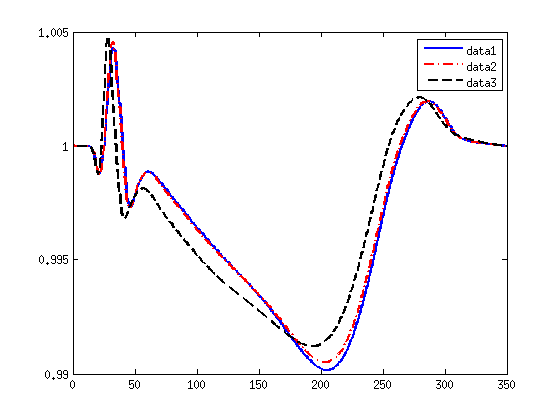}
		\caption{}
		\label{lambda_diff_Ko_with_ISAC_Pt3}
	\end{subfigure}\hspace{-10em}
	\begin{subfigure}[t]{0.5\textwidth}
		\centering
		\includegraphics[width=0.6\textwidth]{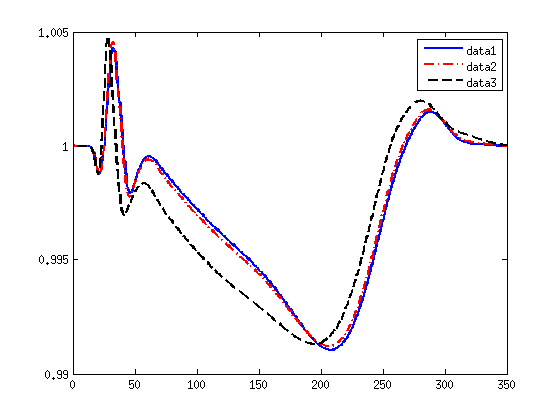}
		\caption{}
		\label{lambda_diff_Ko_with_ISAC_CONV_Pt3}
	\end{subfigure}
	\\
	\begin{subfigure}[t]{0.5\textwidth}
		\includegraphics[width=0.6\textwidth]{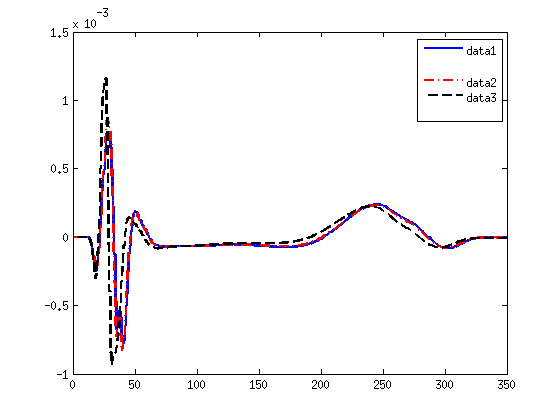}
		\caption{}
		\label{dlambda_diff_Ko_no_ISAC_Pt3}
	\end{subfigure}\hspace{-15em}
	\begin{subfigure}[t]{0.5\textwidth}
		\centering
		\includegraphics[width=0.6\textwidth]{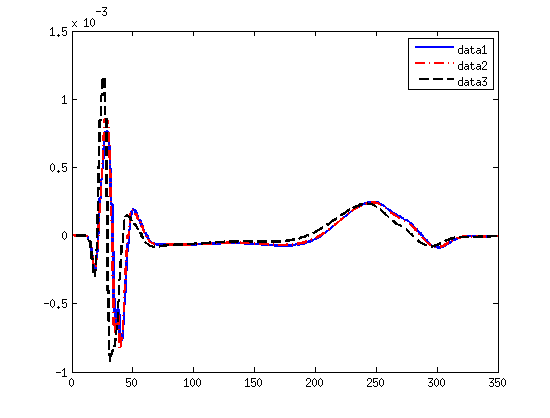}
		\caption{}
		\label{dlambda_diff_Ko_with_ISAC_Pt3}
	\end{subfigure}\hspace{-10em}
	\begin{subfigure}[t]{0.5\textwidth}
		\centering
		\includegraphics[width=0.6\textwidth]{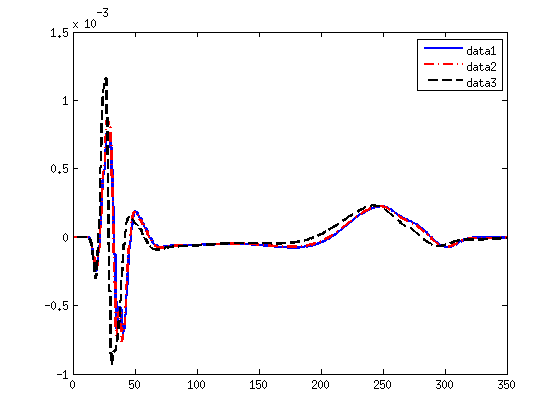}
		\caption{}
		\label{dlambda_diff_Ko_with_ISAC_CONV_Pt3}
	\end{subfigure}
	\caption{: \textbf{$v$ (first row), $[Ca^{+2}]_i$ (second row), $T_A$ (third row), $\lambda$ (fourth row), $\frac{d \lambda}{dt}$ (fifth row) for the cases without $I_{sac}$ (first column), with $I_{sac}$ (second column) and with $I_{sac}+$ CONV (fourth column) at M2, for different values of $[K^+]_o $, data1 ($[K^+]_o =5.4$),  data2 ($[K^+]_o =12$,  data3 ($[K^+]_o =20$).)}}
	\label{diffKo_diff_case_Pt3}
\end{figure*}

$[K^+]_o$ is taken in the range $5.4-20 mM$ \cite{shaw}. $[K^+]_o=5.4mM$ is considered as the normal value. In Fig. \eqref{diffKo_diff_case_Pt1} and \eqref{diffKo_diff_case_Pt3}, action potential $(v)$, intracellular calcium concentration 
($[Ca^{+2}]_i$), active tension $(T)$, stretch along the fiber $(\lambda)$ and stretch rate $(\frac{d\lambda}{dt})$, at the points $M1$ and $M2$ corresponding to S2, S3 and S4 for different values of $[K^+]_o$ is drawn. 
In Fig. \eqref{diffKo_diff_case_Pt1}, effect of Hyperkalemia at the point $M1$ in the ischemic 
region corresponding to S2, S3 and S4 can be easily visualized. With the increase in the strength of Hyperkalemia, resting potential (RP) increasingly reaches to 
the higher level of potential values than that of the normal RP value. RP corresponding to S2, S3 and S4 for $[K^+]_o=20mM$ turns out to be -50mV while the RP for the normal value of $[K^+]_o$ is -85mV. Thus, there is approximately $40\%$ variation in the RP when Hyperkalemia strength increases up to 20mM. For S2, the ischemic cell repolarize faster than that of the normal cell and hence the APD also reduces with the increase 
in the Hyperkalemia strength also the cell gets activated earlier as $[K^+]_o$ reaches to 20mM. This faster repolarization and
reduction in APD at M1 is more evident with the addition of the ionic current generated by the stretch activated channels i.e. $(I_{sac})$. While the addition of convective term (i.e. S4) shows negligible change in AP.  Changes in the values of APD at $M1$ is given as follows: \\

\begin{tabular}{ |c | c |c | c|}
	\hline
	& S2 & S3 & S4 \\ 
	\hline
	APD (at M1 when $[K^+]_o =20$mM) & 272 & 265 & 264 \\
	\hline
\end{tabular}\\\\

From Fig. \eqref{diffKo_diff_case_Pt3}, we can see that at M2 there are small changes in all the phases of the action potential for all the cases S2, S3, and, S4, when $[K^+]_o$ reaches to 20mM. Early depolarization, faster repolarization is visible. With the addition of $(I_{sac})$ nearby healthy cell M2 goes faster to the resting state while the addition of convective term shows negligible effect.

From Fig. \eqref{Cai_diff_Ko_no_ISAC_Pt1}- \eqref{Cai_diff_Ko_with_ISAC_CONV_Pt1}and \eqref{Cai_diff_Ko_no_ISAC_Pt3}-\eqref{Cai_diff_Ko_with_ISAC_CONV_Pt3}, at M1 and M2, changes in the waveform of the concentration of the intracellular calcium ions $[Ca^{+2}]_i$ for the cases S2, S3, and S4 are negligible when $[K^+]_o \leq 12mM$. Changes in the mechanical waveform of $[Ca^{+2}]_i$ are visible
when $[K^+]_o $ reaches 20mM and this change in $[Ca^{+2}]_i$ waveform leads to change in contraction force and hence contractility of the cardiac tissue.  
It is to be noticed that at $M1$ and $M2$ the effect in the waveform of $T_A$ for the cases S2, S3, and S4 is negligible when $[K^+]_o \leq 12mM$ and the effect in the change in mechanical waveform of $T_A$  is visible when $[K^+]_o$ 
reaches to 20mM.  
It is also clearly visible that the influence of Hyperkalemia on the stretch along fiber $\lambda$ and stretch rate, $\frac{d \lambda}{dt}$ i.e. the change in the waveform of $\lambda$  and hence $\frac{d \lambda}{dt}$ is negligible at M1 corresponding to S2, S3 and S4. The change in the waveform of $\lambda$ and hence $\frac{d \lambda}{dt}$ is more at the point M2 i.e. stretch along the fiber and stretch rate is more with the increase in Hyperkalemia strength. 
There is more than 3\% change in the length of the myocytes relative to the resting length at point M2.

\begin{figure*}
	\hspace{-5em}
	\begin{subfigure}[t]{0.5\textwidth}
		\centering
		\includegraphics[width=0.7\textwidth]{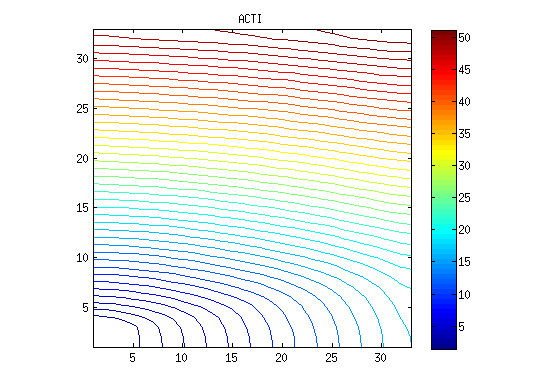}
		\caption{}
	\end{subfigure} %
	\hspace{-10em}
	\begin{subfigure}[t]{0.5\textwidth}
		\centering
		\includegraphics[width=0.7\textwidth]{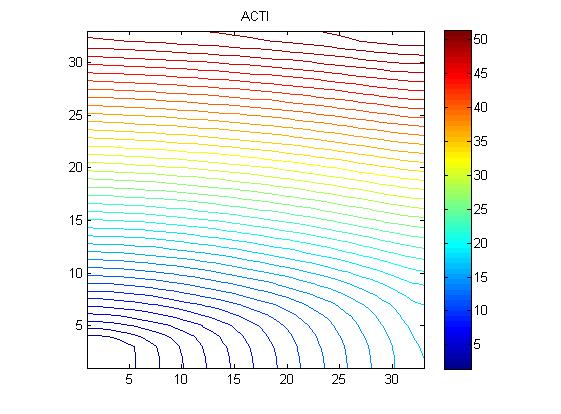}
		\caption{}
	\end{subfigure}\hspace{-10em}
	\begin{subfigure}[t]{0.5\textwidth}
		\centering
		\includegraphics[width=0.7\textwidth]{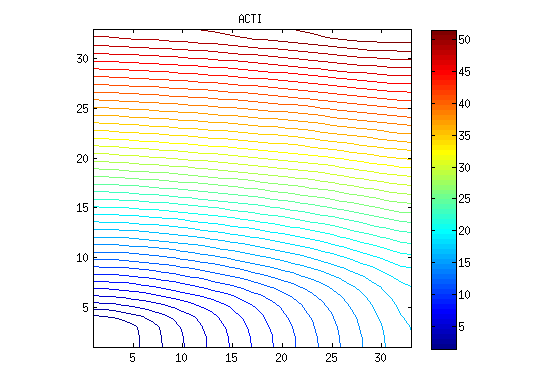}
		\caption{}
	\end{subfigure} \\
	\begin{subfigure}[t]{0.5\textwidth}
		\includegraphics[width=0.7\textwidth]{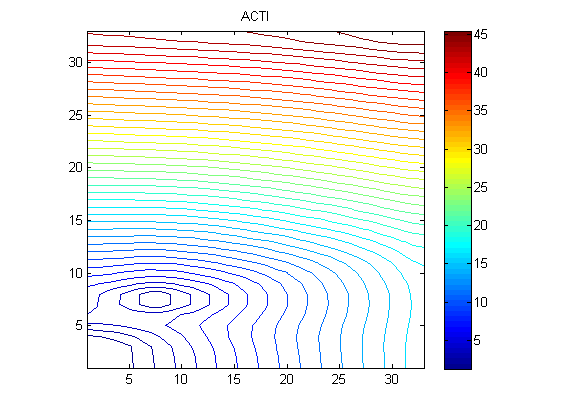}
		\caption{}
	\end{subfigure}\hspace{-14em}
	\begin{subfigure}[t]{0.5\textwidth}
		\centering
		\includegraphics[width=0.7\textwidth]{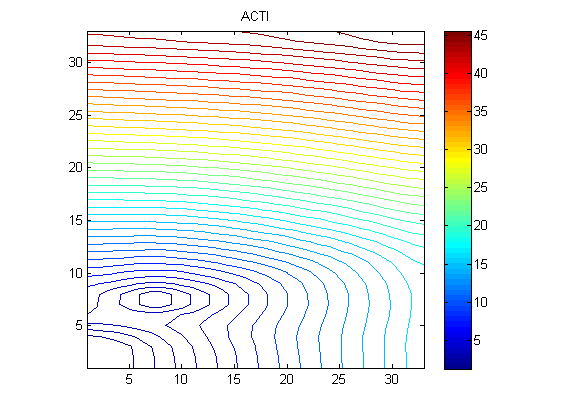}
		\caption{}
	\end{subfigure}\hspace{-10em}
	\begin{subfigure}[t]{0.5\textwidth}
		\centering
		\includegraphics[width=0.7\textwidth]{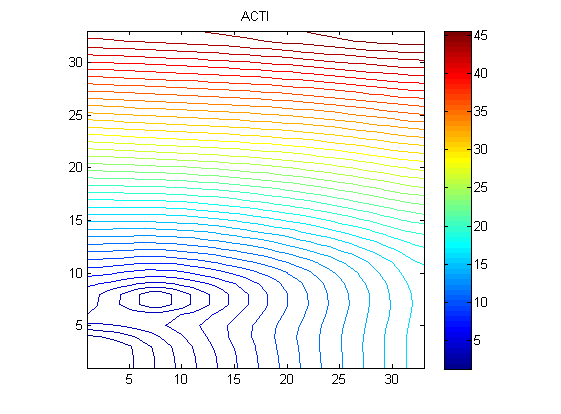}
		\caption{}
	\end{subfigure}\\
	\caption{: \textbf{AT for the cases without $I_{sac}$ (first column), with $I_{sac}$ (second column) and with $I_{sac}+$ CONV (third column) for $[K^+]_o =5.4$mM (first row) and $[K^+]_o =20$mM (second row).}}
	\label{AT_Ko_20}
\end{figure*}

\begin{figure*}
	\hspace{-5em}
	\begin{subfigure}[t]{0.5\textwidth}
		\centering
		\includegraphics[width=0.7\textwidth]{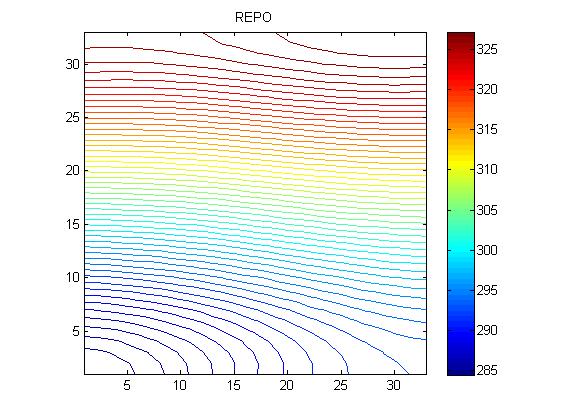}
		\caption{}
	\end{subfigure} %
	\hspace{-10em}
	\begin{subfigure}[t]{0.5\textwidth}
		\centering
		\includegraphics[width=0.75\textwidth]{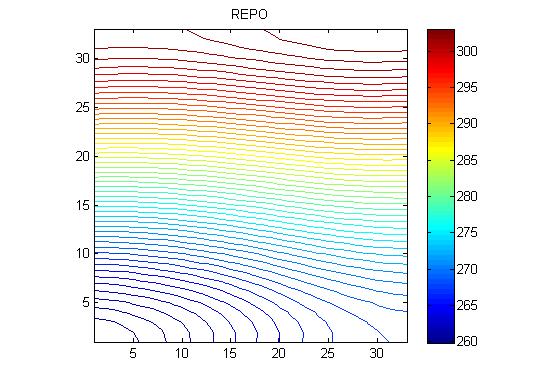}
		\caption{}
	\end{subfigure}\hspace{-10em}
	\begin{subfigure}[t]{0.5\textwidth}
		\centering
		\includegraphics[width=0.7\textwidth]{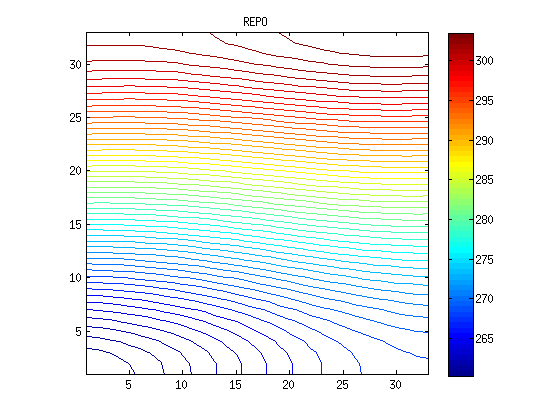}
		\caption{}
	\end{subfigure} \\
	\begin{subfigure}[t]{0.5\textwidth}
		\includegraphics[width=0.7\textwidth]{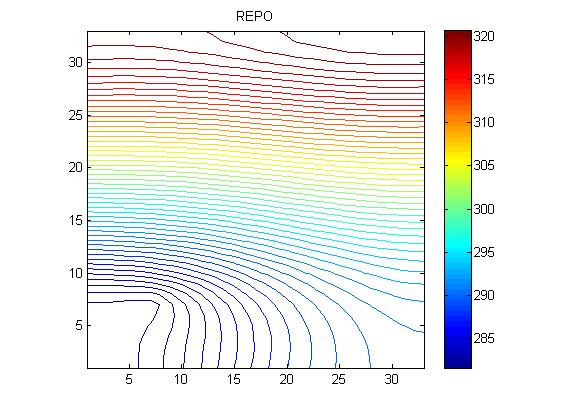}
		\caption{}
	\end{subfigure}\hspace{-14em}
	\begin{subfigure}[t]{0.5\textwidth}
		\centering
		\includegraphics[width=0.7\textwidth]{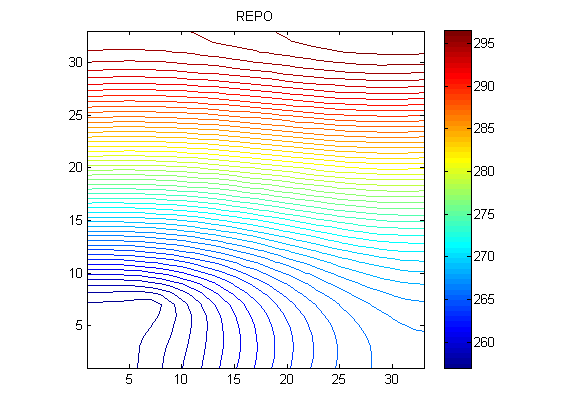}
		\caption{}
	\end{subfigure}\hspace{-10em}
	\begin{subfigure}[t]{0.5\textwidth}
		\centering
		\includegraphics[width=0.7\textwidth]{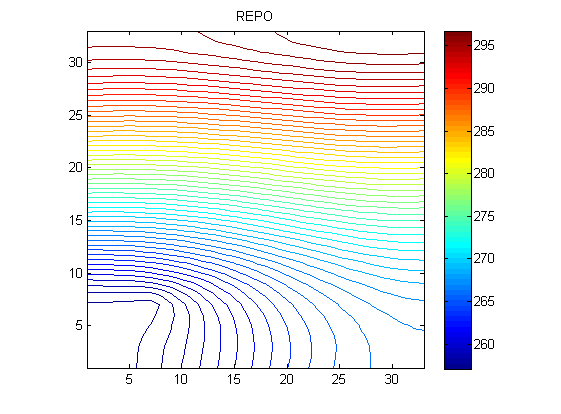}
		\caption{}
	\end{subfigure}\\
	\caption{: \textbf{RT for the cases without $I_{sac}$ (first column), with $I_{sac}$ (second column) and with $I_{sac}+$ CONV (third column) for $[K^+]_o =5.4$mM (first row) and $[K^+]_o =20$mM (second row).}}
	\label{RT_Ko_20}
\end{figure*}

\begin{figure*}
	\hspace{-5em}
	\begin{subfigure}[t]{0.5\textwidth}
		\centering
		\includegraphics[width=0.7\textwidth]{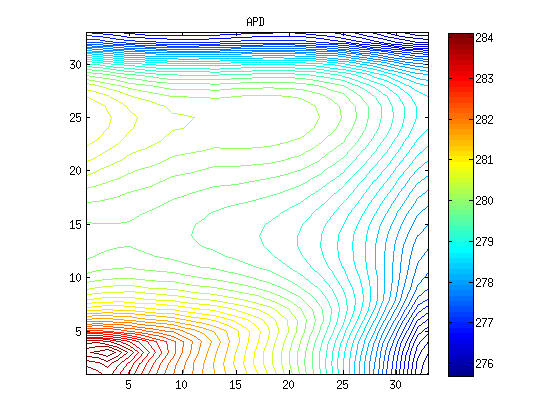}
		\caption{}
	\end{subfigure} %
	\hspace{-10em}
	\begin{subfigure}[t]{0.5\textwidth}
		\centering
		\includegraphics[width=0.75\textwidth]{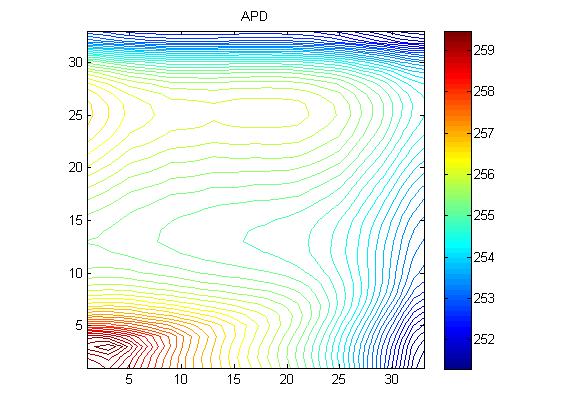}
		\caption{}
	\end{subfigure}\hspace{-10em}
	\begin{subfigure}[t]{0.5\textwidth}
		\centering
		\includegraphics[width=0.7\textwidth]{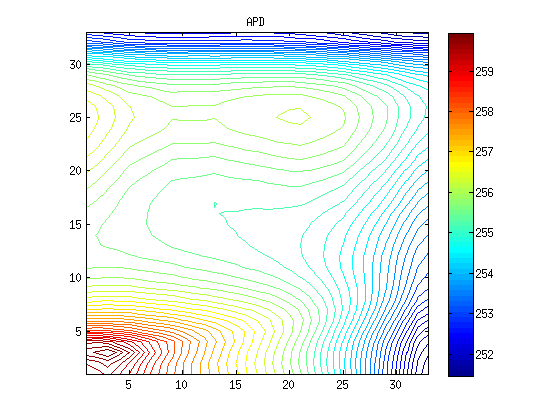}
		\caption{}
	\end{subfigure} \\
	\begin{subfigure}[t]{0.5\textwidth}
		\includegraphics[width=0.7\textwidth]{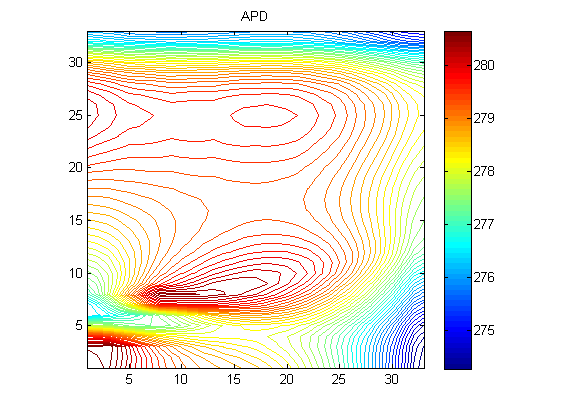}
		\caption{}
	\end{subfigure}\hspace{-14em}
	\begin{subfigure}[t]{0.5\textwidth}
		\centering
		\includegraphics[width=0.7\textwidth]{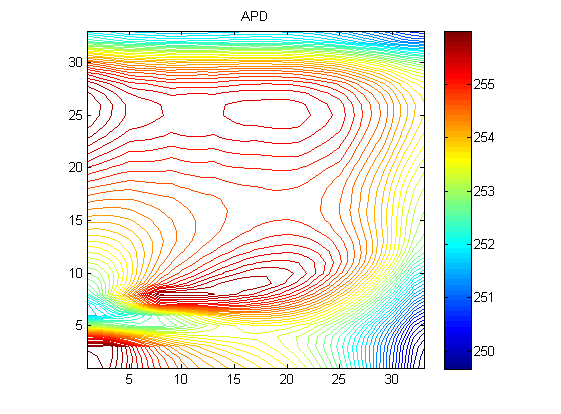}
		\caption{}
	\end{subfigure}\hspace{-10em}
	\begin{subfigure}[t]{0.5\textwidth}
		\centering
		\includegraphics[width=0.7\textwidth]{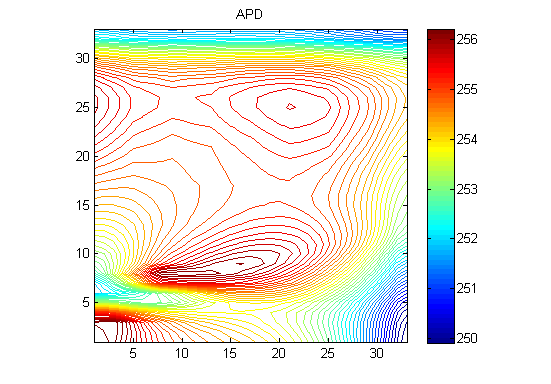}
		\caption{}
	\end{subfigure}\\
	\caption{: \textbf{APD for the cases without $I_{sac}$ (first column), with $I_{sac}$ (second column) and with $I_{sac}+$ CONV (third column) for $[K^+]_o =5.4$mM (first row) and $[K^+]_o =20$mM (second row).}}
	\label{APD_Ko_20}
\end{figure*}

AT, RT and APD contours for S2, S3 and S4 with $[K^+]_o =20mM$ is presented in Fig. \eqref{AT_Ko_20}- \eqref{APD_Ko_20}. From Fig. \eqref{AT_Ko_20}, it is visualized form the isochrone lines that the addition of $I_{sac}$ does not affect the AT pattern. Also, increasing the level of Hyperkalemia leads to a visible hastening in the activation of human cardiac cells. An elliptic area is clearly visible which activated earlier than the surrounding one. This indicates that ischemic region generate a new wave.

Fig. \eqref{RT_Ko_20} and \eqref{APD_Ko_20} indicates that due to the addition of $I_{sac}$ cardiac cells slightly repolarize faster and hence reduces APD. Also the cardiac cells near to the ischemic region repolarize faster or RT of the ischemic cells or neighboring to the ischemic cells reduces and hence these cells go to the resting state earlier and this tendency increases as the Hyperkalemia gets severe.

\begin{figure*}
	\hspace{-6em}
	\begin{subfigure}[t]{0.5\textwidth}
		\centering
		\includegraphics[width=0.6\textwidth]{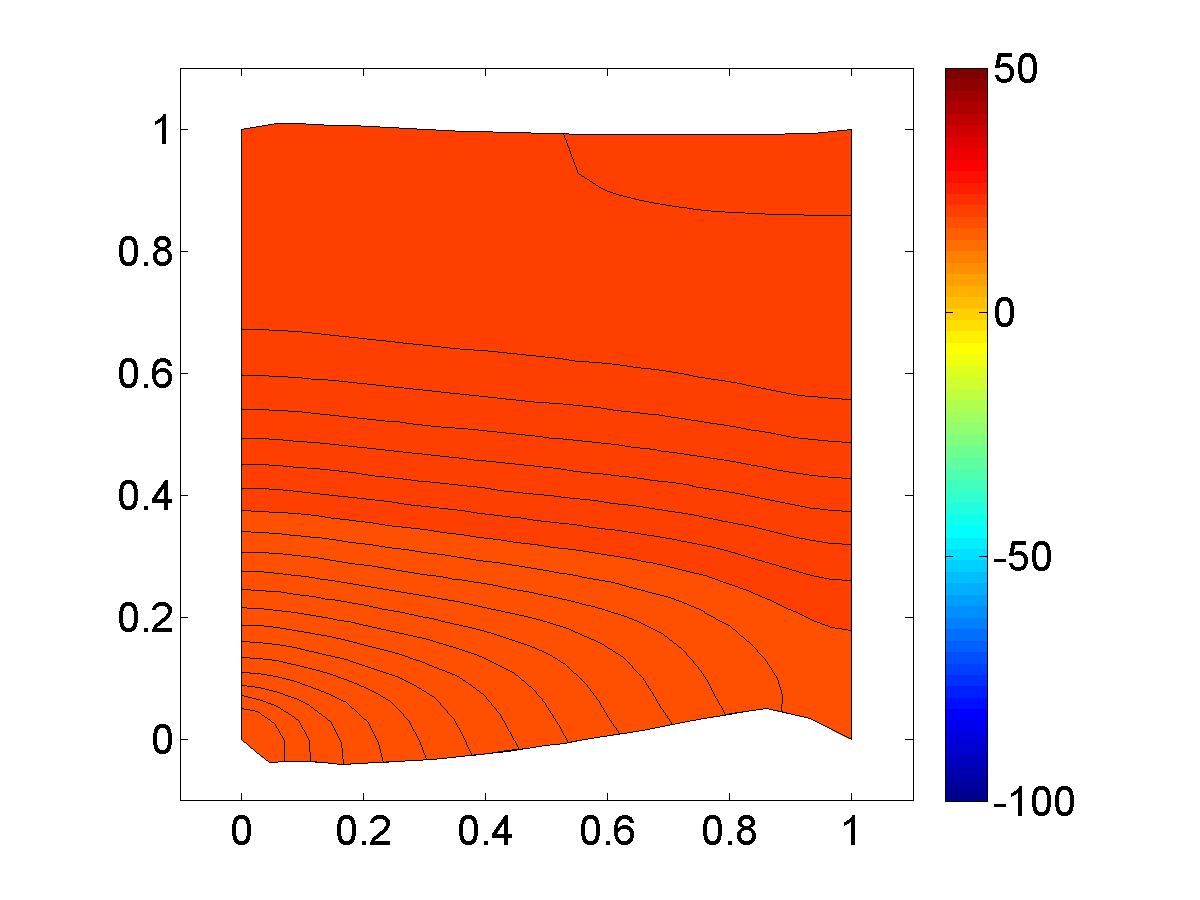}
		\vspace{-1.1em}
		\caption{}
	\end{subfigure} %
	\hspace{-12em}
	\begin{subfigure}[t]{0.5\textwidth}
		\centering
		\includegraphics[width=0.6\textwidth]{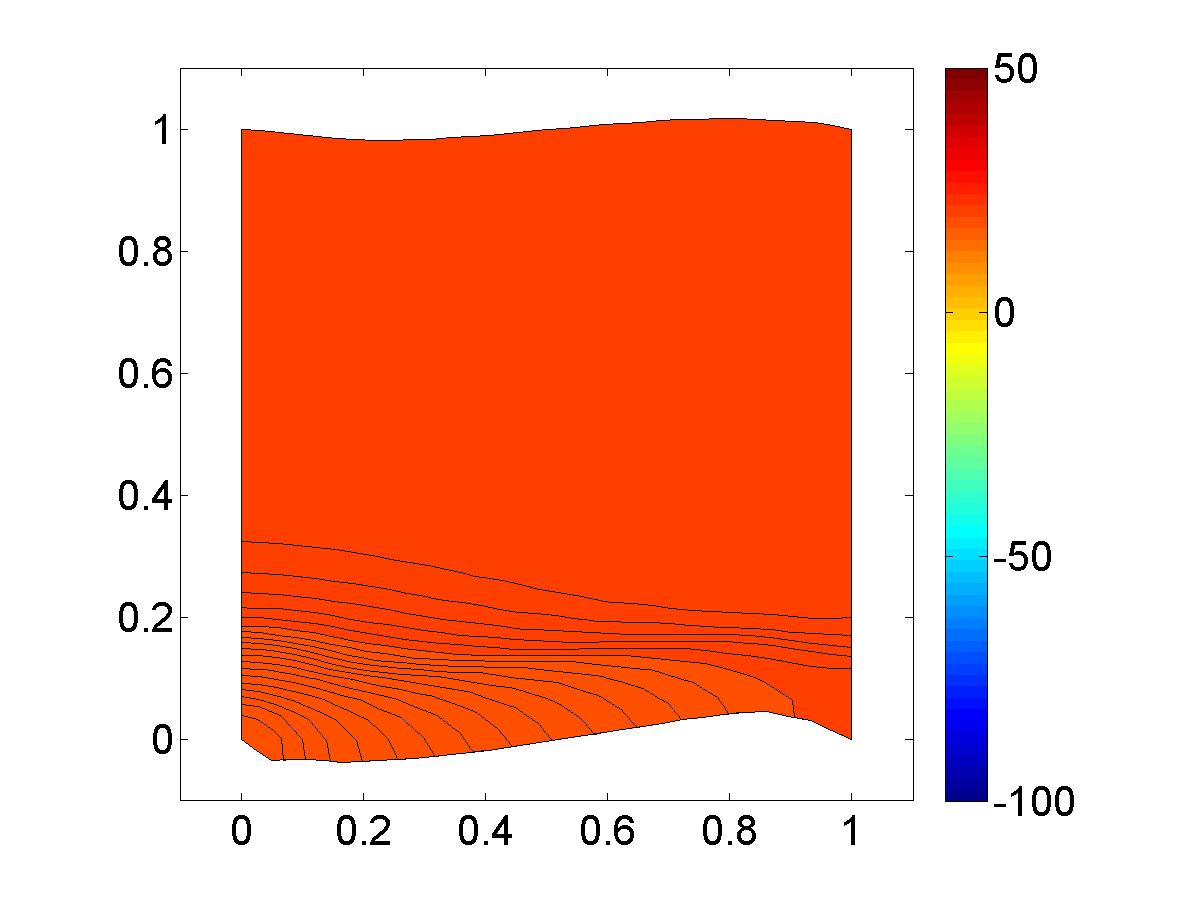}
		\vspace{-1.1em}
		\caption{}
	\end{subfigure}\hspace{-12em}
	\begin{subfigure}[t]{0.5\textwidth}
		\centering
		\includegraphics[width=0.6\textwidth]{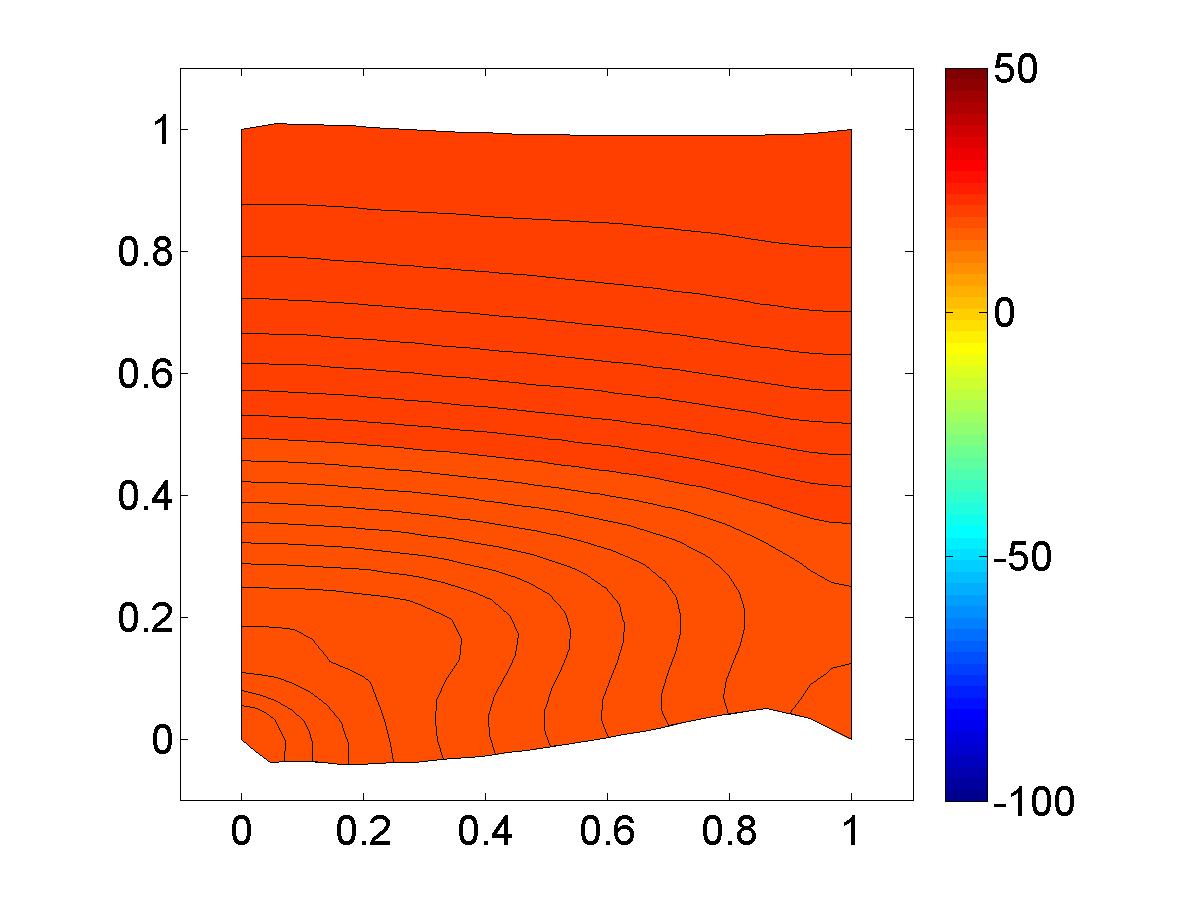}
		\vspace{-1.1em}
		\caption{}
	\end{subfigure} \\
	\vspace{-1em}
	\begin{subfigure}[t]{0.5\textwidth}
		\includegraphics[width=0.6\textwidth]{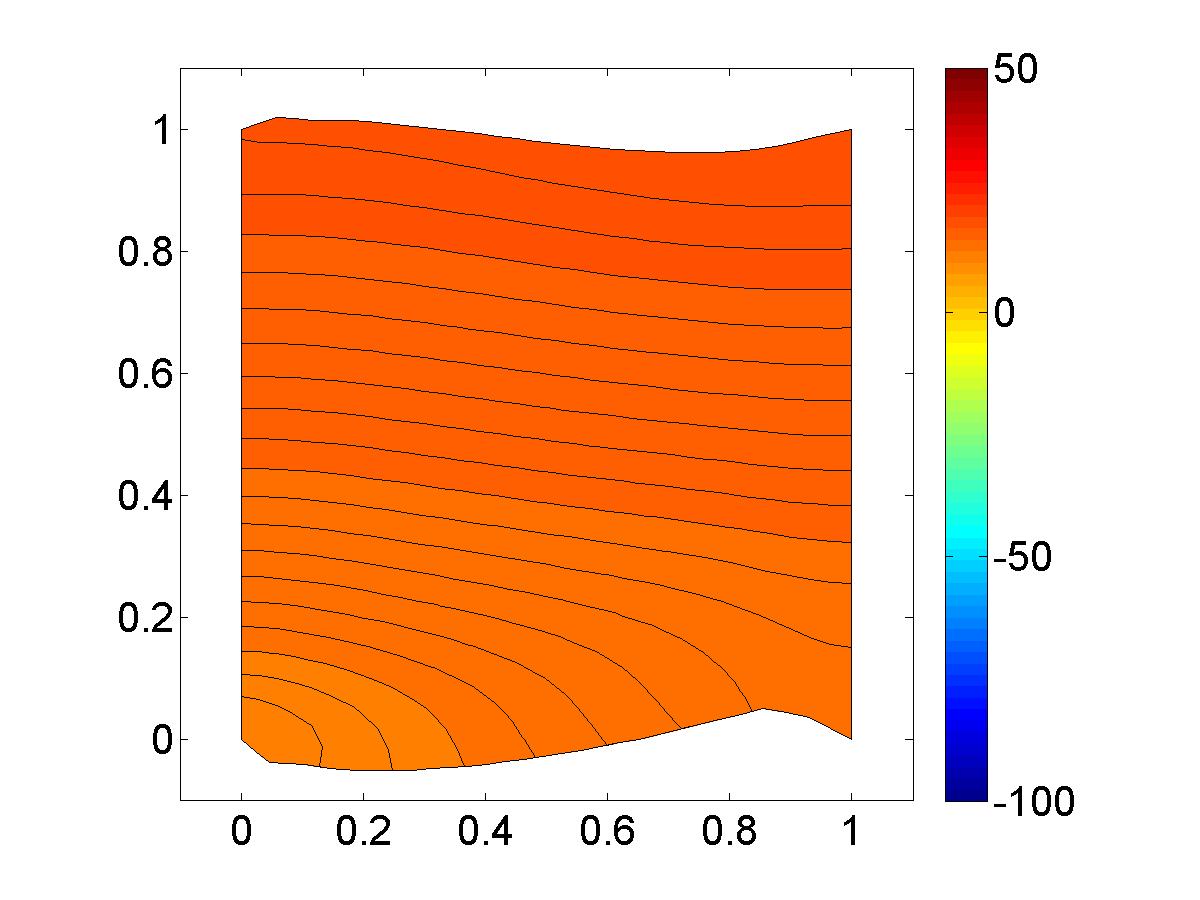}
		\vspace{-1.1em}
		\caption{}
	\end{subfigure}\hspace{-17em}
	\begin{subfigure}[t]{0.5\textwidth}
		\centering
		\includegraphics[width=0.6\textwidth]{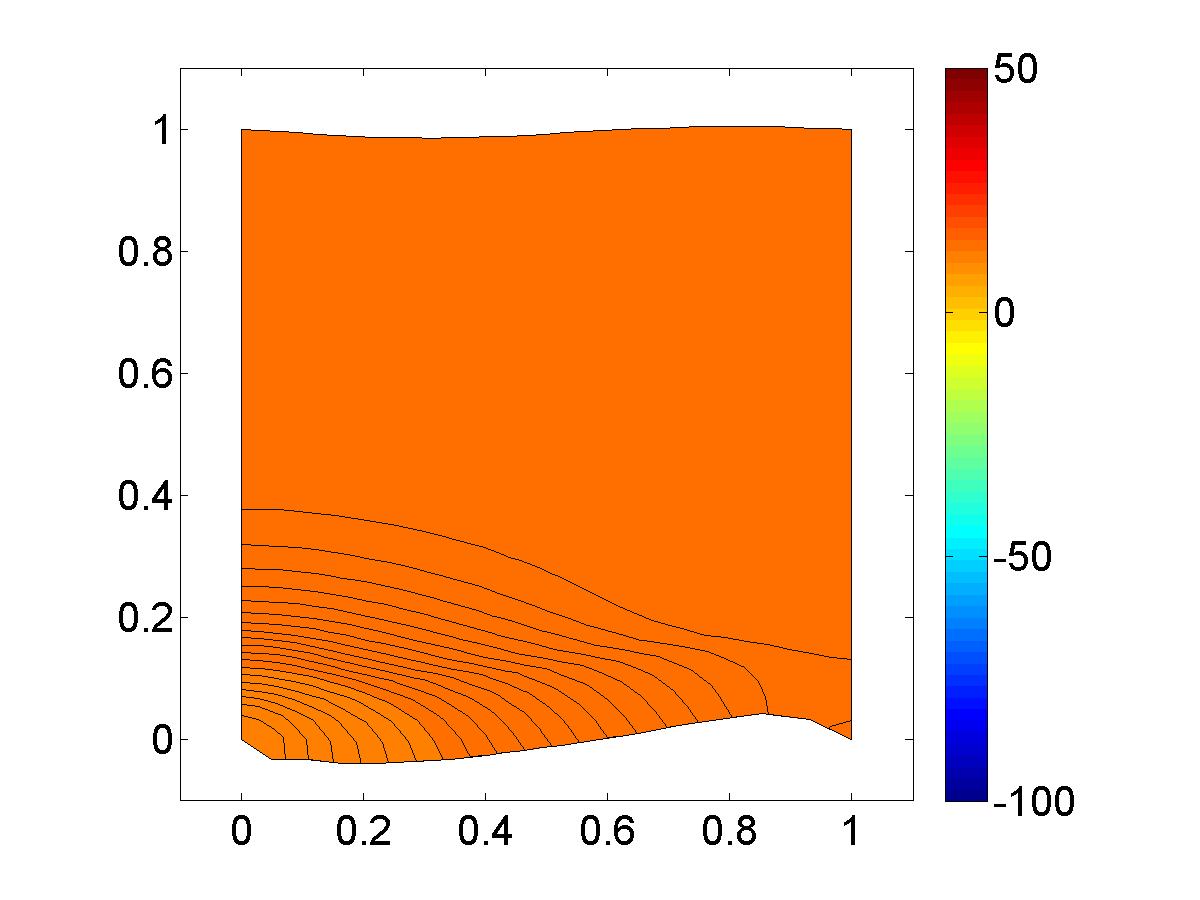}
		\vspace{-1.1em}
		\caption{}
	\end{subfigure}\hspace{-12em}
	\begin{subfigure}[t]{0.5\textwidth}
		\centering
		\includegraphics[width=0.6\textwidth]{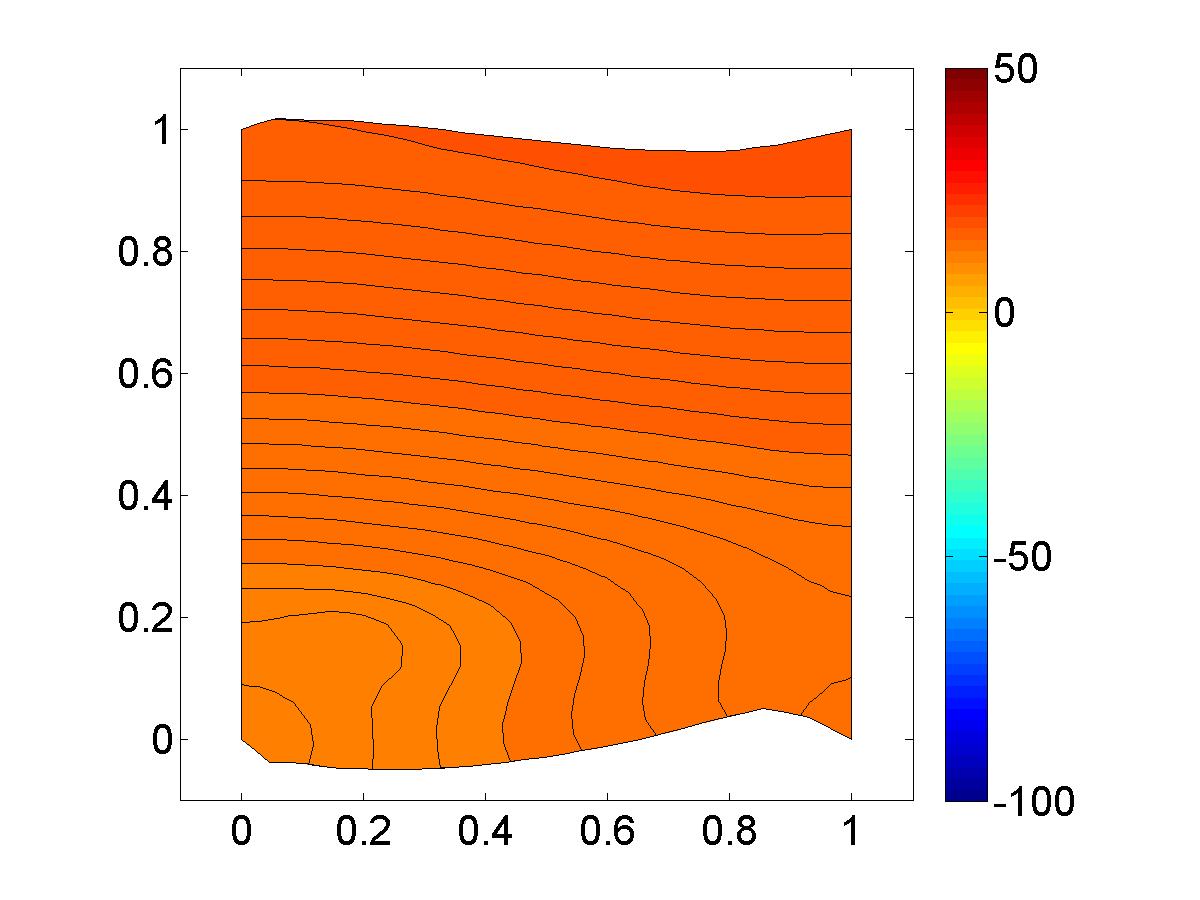}
		\vspace{-1.1em}
		\caption{}
	\end{subfigure} \\
	\begin{subfigure}[t]{0.5\textwidth}
		\includegraphics[width=0.6\textwidth]{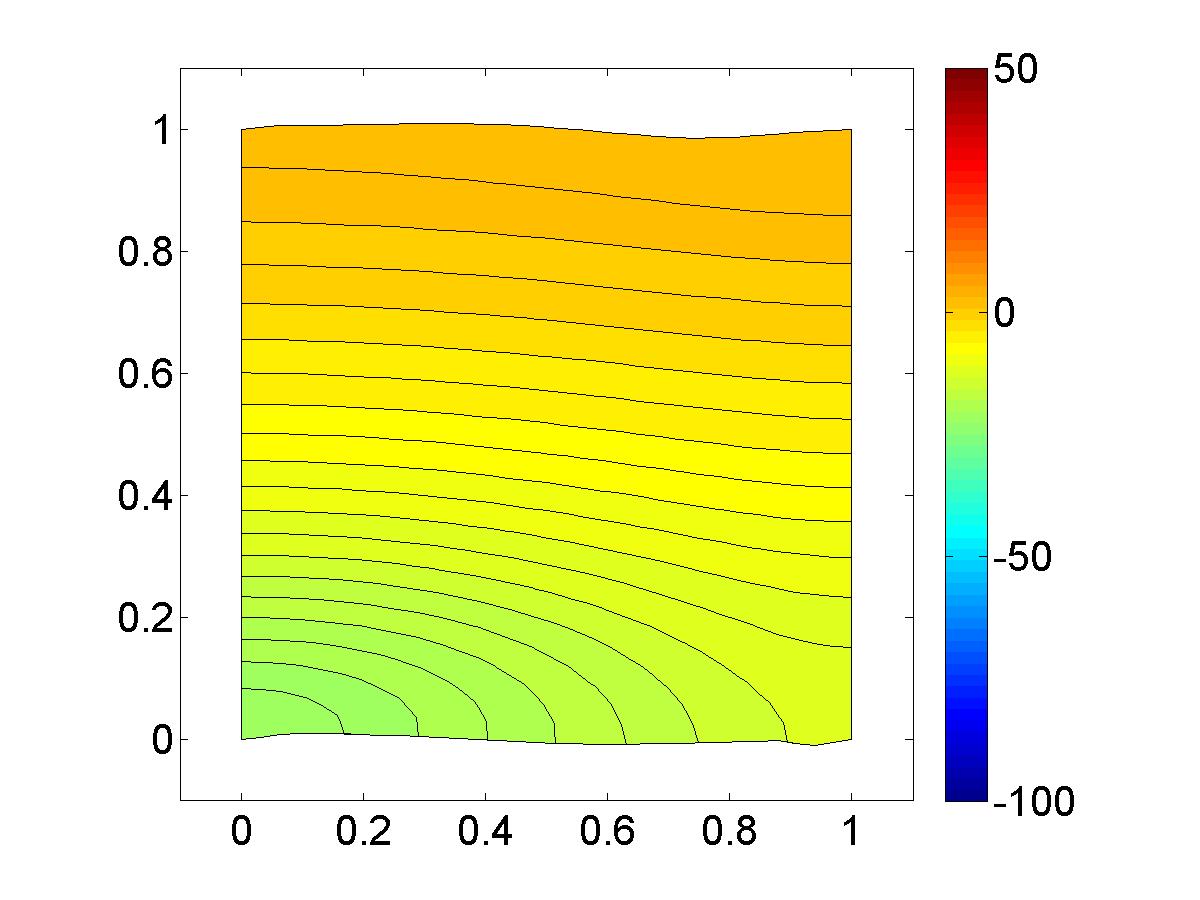}
		\vspace{-1.1em}
		\caption{}
	\end{subfigure}\hspace{-17em}
	\begin{subfigure}[t]{0.5\textwidth}
		\centering
		\includegraphics[width=0.6\textwidth]{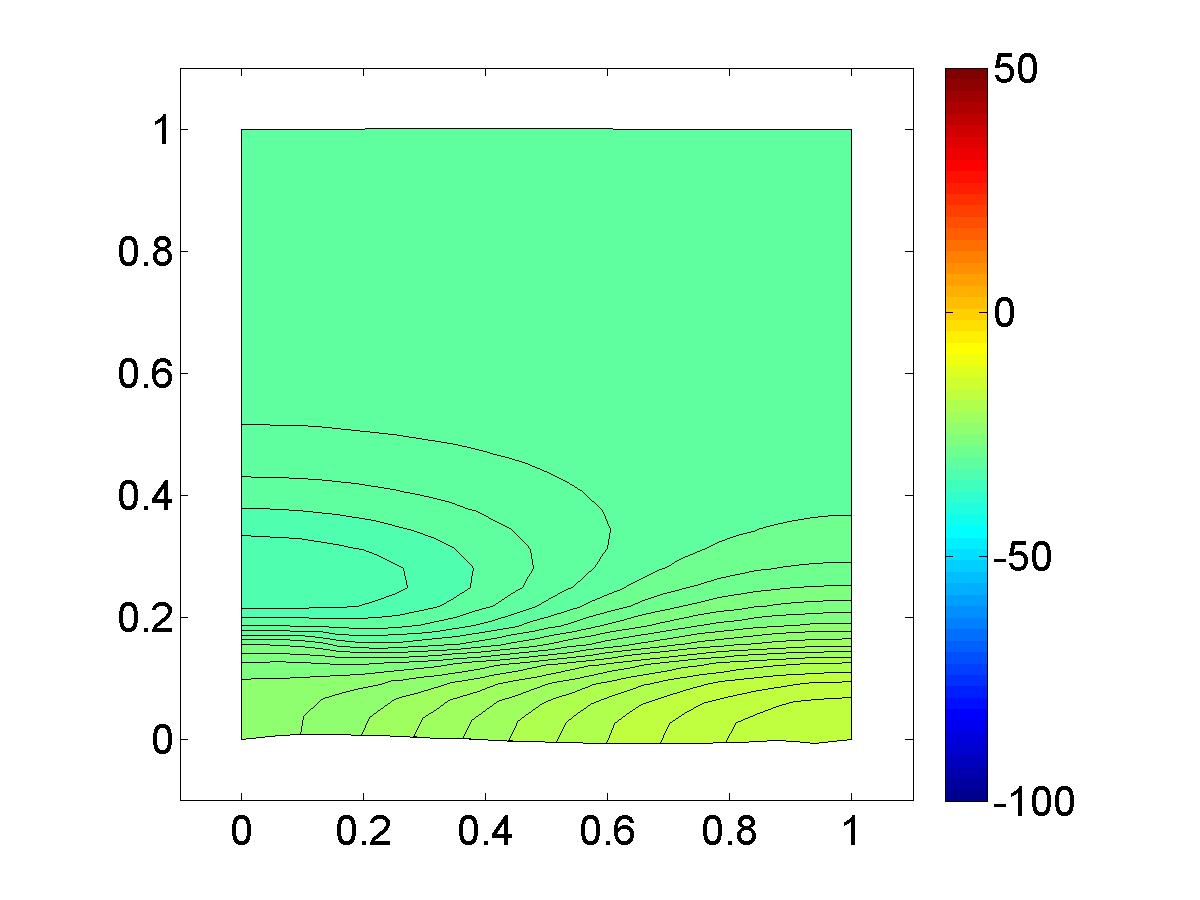}
		\vspace{-1.1em}
		\caption{}
	\end{subfigure}\hspace{-12em}
	\begin{subfigure}[t]{0.5\textwidth}
		\centering
		\includegraphics[width=0.6\textwidth]{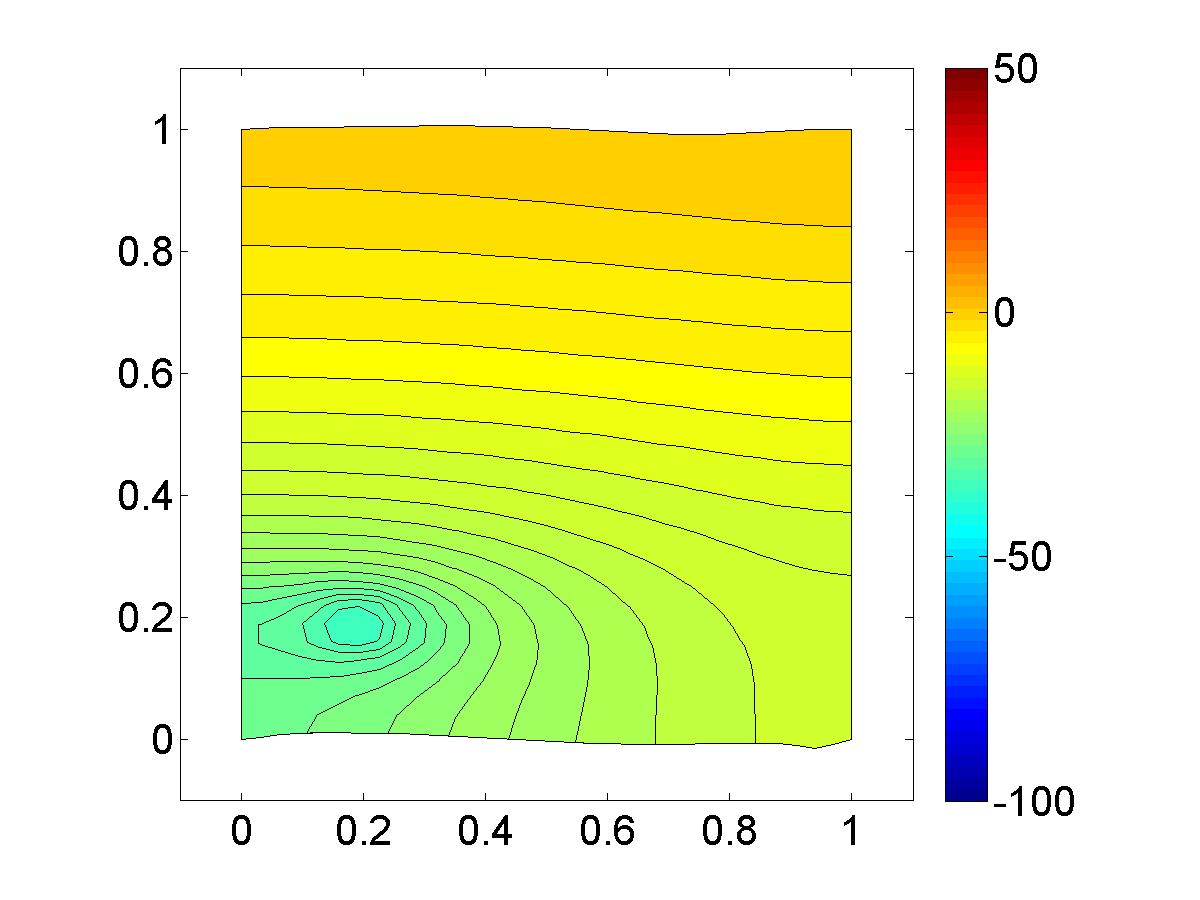}
		\vspace{-1.1em}
		\caption{}
	\end{subfigure}\\
	\vspace{-1.5em}
	\caption{: \textbf{Action potential contours in a deforming domain for the cases without $I_{sac}$ for different values of $[K^+]_o$ = 5.4 mM (first column), $[K^+]_o $=12 mM (second column), $[K^+]_o $ = 20 mM (third column), at time t=90 (first row), t=150 (second row), t=240 (third row).}}
	\label{AP_diffKo_contours_No_ISAC}
\end{figure*}

\begin{figure*}
	\hspace{-0em}		
	\begin{subfigure}[t]{0.5\textwidth}
		\includegraphics[width=0.6\textwidth]{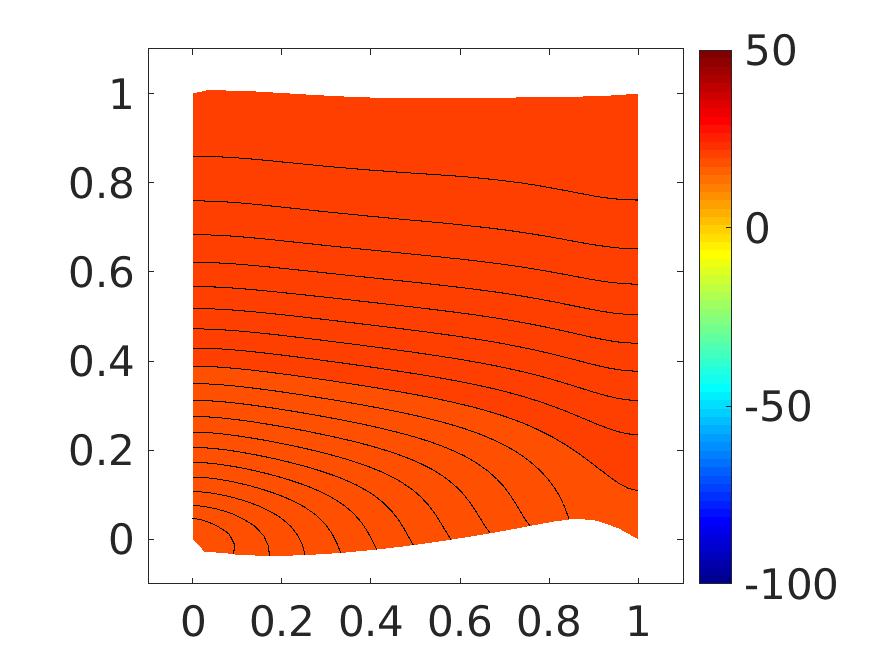}
		\vspace{-1.1em}
		\caption{}
	\end{subfigure} %
	\hspace{-17em}
	\begin{subfigure}[t]{0.5\textwidth}
		\centering
		\includegraphics[width=0.6\textwidth]{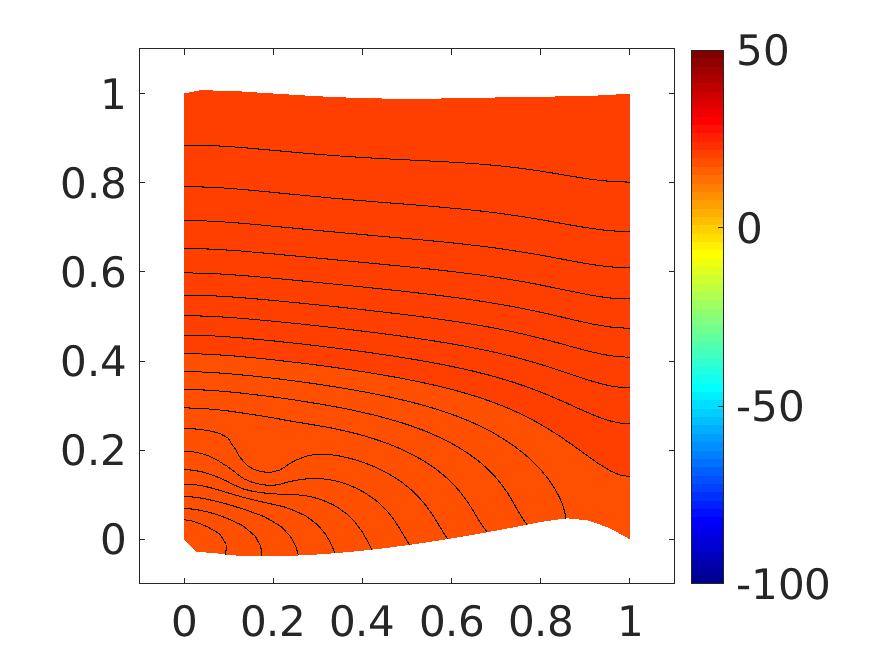}
		\vspace{-1.1em}
		\caption{}
	\end{subfigure}\hspace{-12em}
	\begin{subfigure}[t]{0.5\textwidth}
		\centering
		\includegraphics[width=0.6\textwidth]{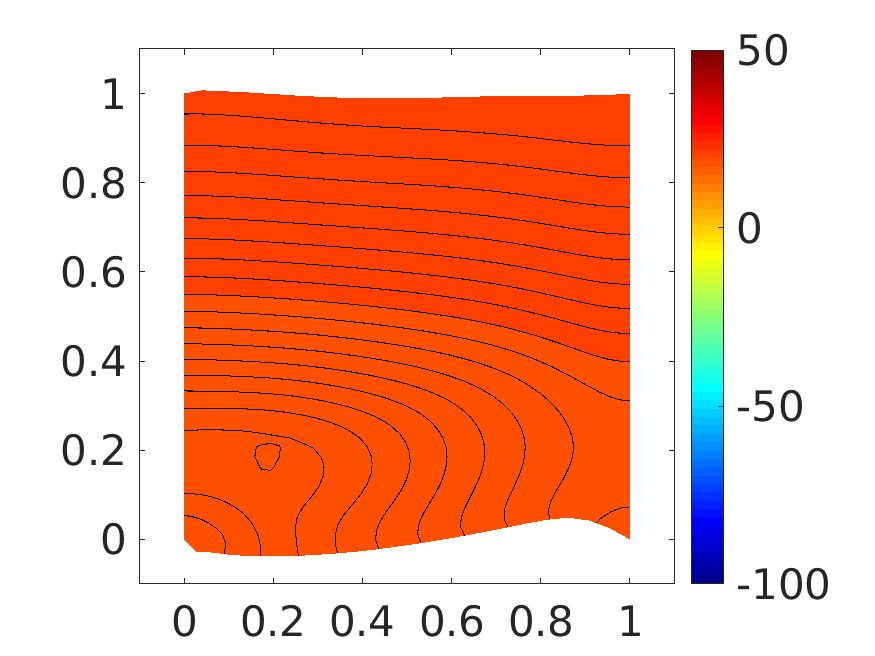}
		\vspace{-1.1em}
		\caption{}
	\end{subfigure} \\
	\vspace{-0.8em}
	\begin{subfigure}[t]{0.5\textwidth}
		\includegraphics[width=0.6\textwidth]{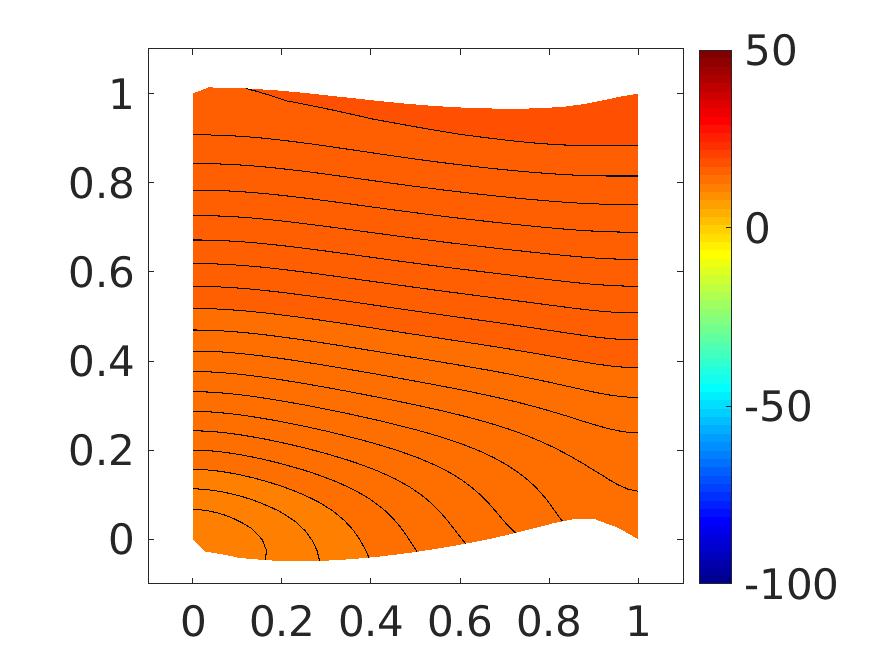}
		\vspace{-1.1em}
		\caption{}
	\end{subfigure}\hspace{-17em}
	\begin{subfigure}[t]{0.5\textwidth}
		\centering
		\includegraphics[width=0.6\textwidth]{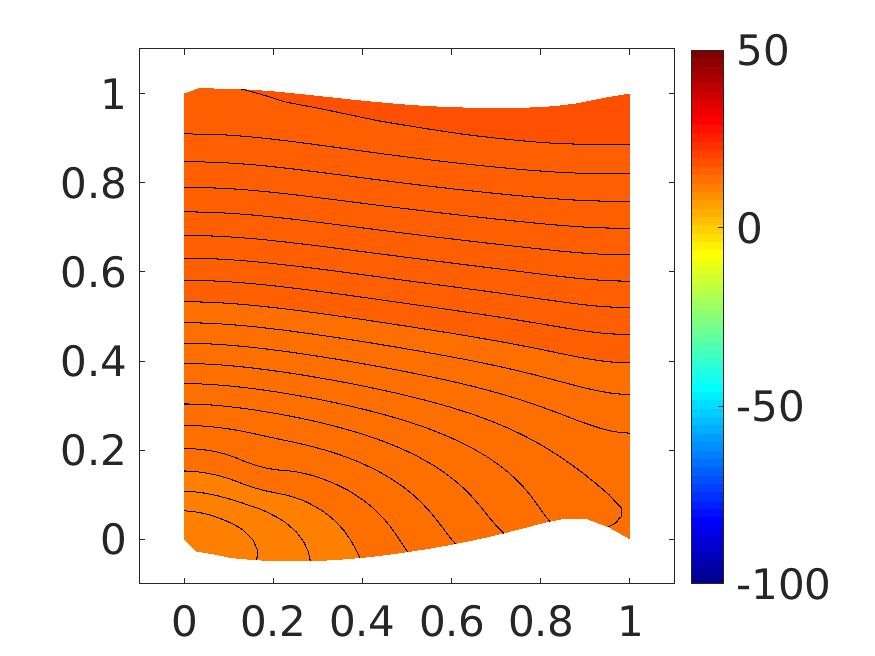}
		\vspace{-1.1em}
		\caption{}
	\end{subfigure}\hspace{-12em}
	\begin{subfigure}[t]{0.5\textwidth}
		\centering
		\includegraphics[width=0.6\textwidth]{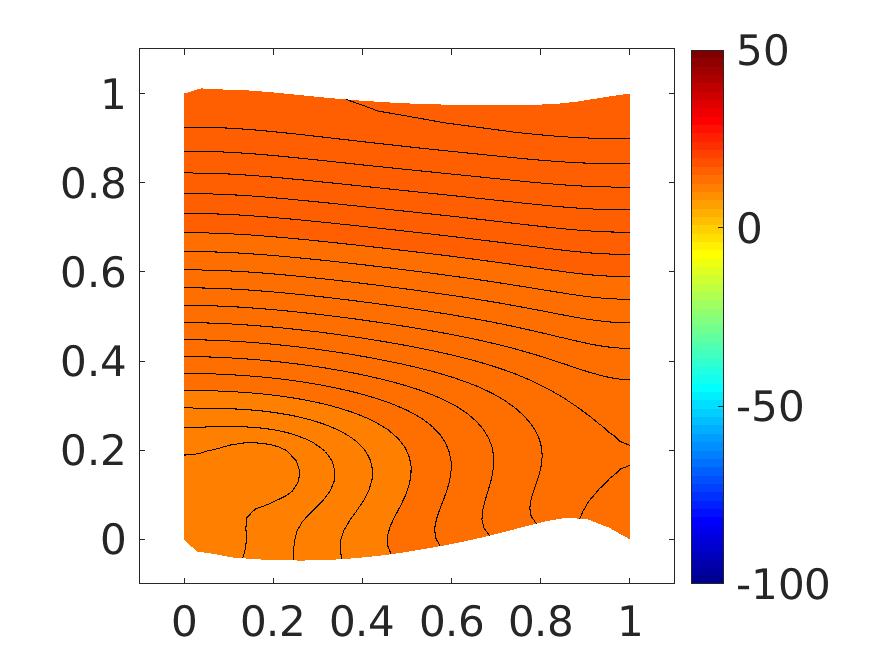}
		\vspace{-1.1em}
		\caption{}
	\end{subfigure} \\
	\vspace{-0.8em}
	\begin{subfigure}[t]{0.5\textwidth}
		\includegraphics[width=0.6\textwidth]{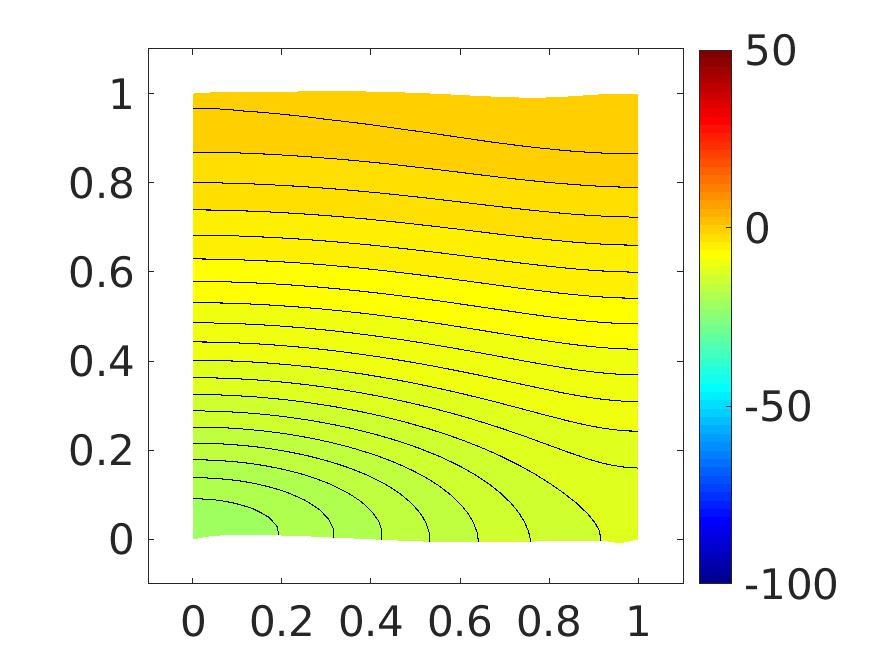}
		\vspace{-1.1em}
		\caption{}
	\end{subfigure}\hspace{-17em}
	\begin{subfigure}[t]{0.5\textwidth}
		\centering
		\includegraphics[width=0.6\textwidth]{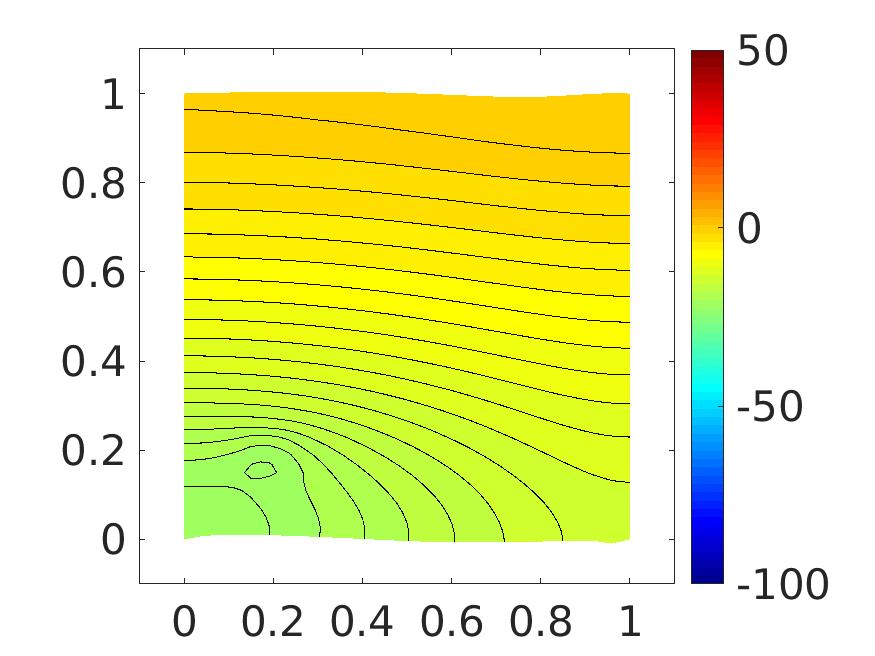}
		\vspace{-1.1em}
		\caption{}
	\end{subfigure}\hspace{-12em}
	\begin{subfigure}[t]{0.5\textwidth}
		\centering
		\includegraphics[width=0.6\textwidth]{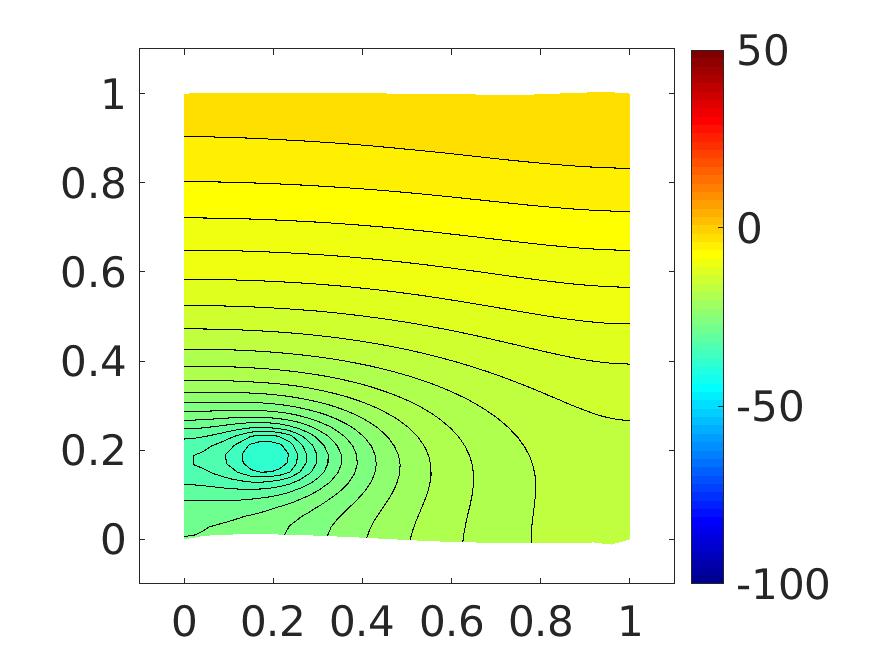}
		\vspace{-1.1em}
		\caption{}
	\end{subfigure}
	\vspace{0em}
	\caption{: \textbf{Action potential contours in a deforming domain for the cases with $I_{sac}$ for different values of $[K^+]_o$ =5.4 (first column), $[K^+]_o $=12 (second column), $[K^+]_o $ = 20 (third column), at time t=90 (first row), t=150 (second row), t=240 (third row).}}
	\label{AP_diffKo_contours_with_ISAC}
\end{figure*}

In Fig. \eqref{AP_diffKo_contours_No_ISAC} and \eqref{AP_diffKo_contours_with_ISAC}, we presented the action potential contours in deformation of the cardiac tissue with the increase in the strength of Hyperkalemia, for S2 and S3, respectively. The change in the deformation of cardiac domain with increasing the Hyperkalemic strength is clearly visible.

\begin{figure*}
	\centering
	\begin{subfigure}[t]{0.5\textwidth}
		\centering
		\includegraphics[width=0.7\textwidth]{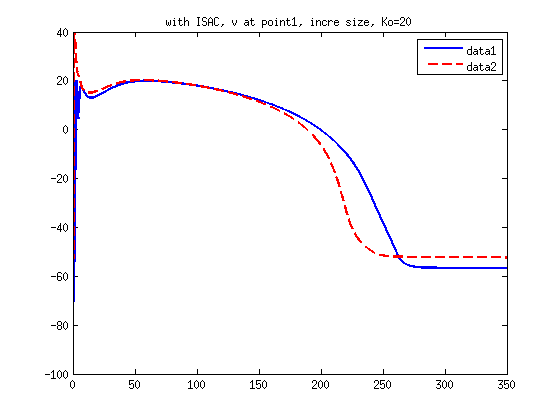}
		\vspace{-1em}
		\caption{}
	\end{subfigure}\hfill	
	\hspace{-10em}
	\begin{subfigure}[t]{0.5\textwidth}
		\centering
		\includegraphics[width=0.7\textwidth]{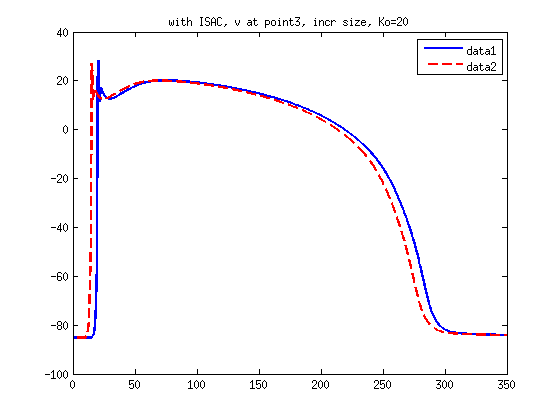}
		\vspace{-1em}
		\caption{}
	\end{subfigure}\\
	\begin{subfigure}[t]{0.5\textwidth}
		\centering
		\includegraphics[width=0.7\textwidth]{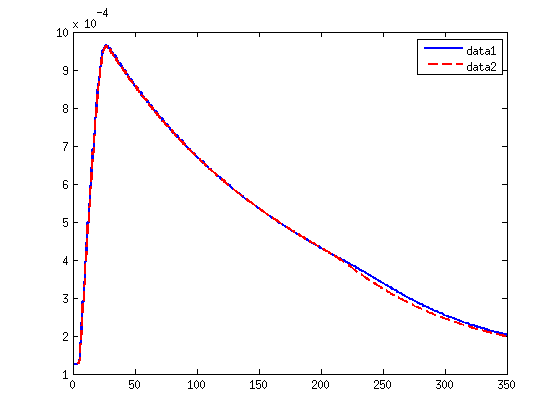}
		\vspace{-1em}
		\caption{}
	\end{subfigure}\hfill
	\hspace{-10em}
	\begin{subfigure}[t]{0.5\textwidth}
		\centering
		\includegraphics[width=0.7\textwidth]{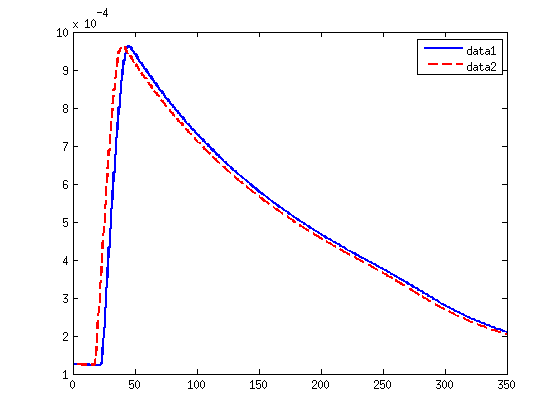}
		\vspace{-1em}
		\caption{}
	\end{subfigure}\\
	\begin{subfigure}[t]{0.5\textwidth}
		\centering
		\includegraphics[width=0.7\textwidth]{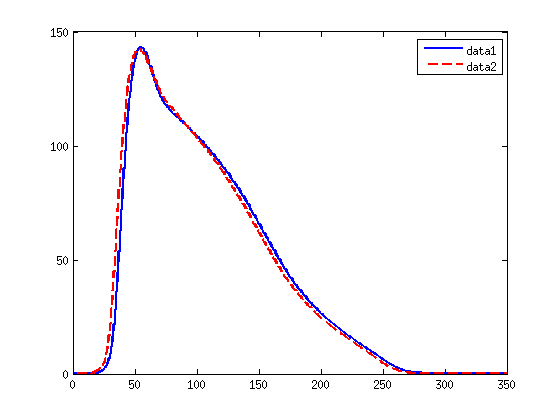}
		\vspace{-1em}
		\caption{}
	\end{subfigure}\hfill
	\hspace{-10em}
	\begin{subfigure}[t]{0.5\textwidth}
		\centering
		\includegraphics[width=0.7\textwidth]{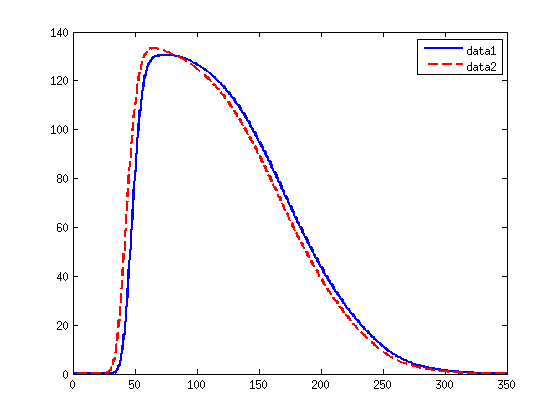}
		\vspace{-1em}
		\caption{}
	\end{subfigure}\\
	\begin{subfigure}[t]{0.5\textwidth}
		\centering
		\includegraphics[width=0.7\textwidth]{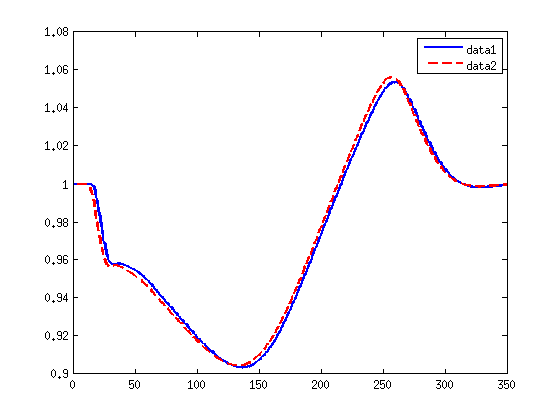}
		\vspace{-1em}
		\caption{}
		\label{}
	\end{subfigure}\hfill
	\hspace{-10em}
	\begin{subfigure}[t]{0.5\textwidth}
		\centering
		\includegraphics[width=0.7\textwidth]{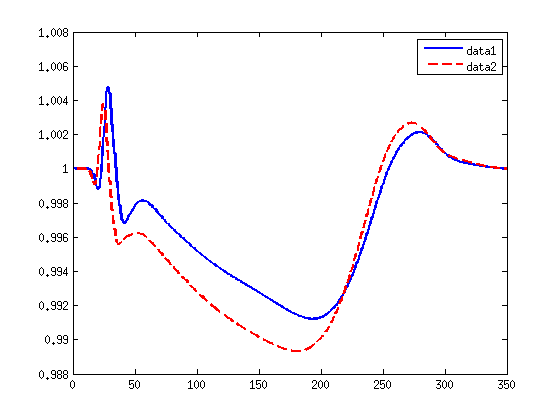}
		\vspace{-1em}
		\caption{}
	\end{subfigure}\\
	\begin{subfigure}[t]{0.5\textwidth}
		\centering
		\includegraphics[width=0.7\textwidth]{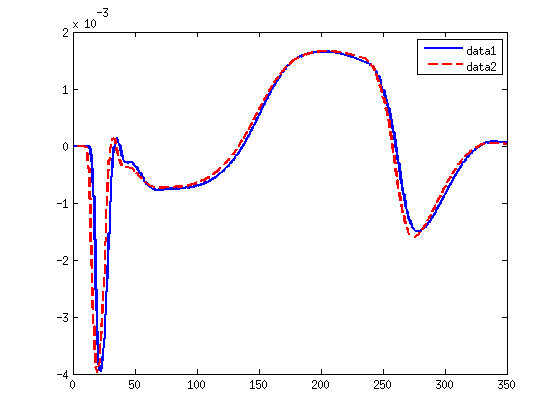}
		\vspace{-1em}
		\caption{}
		\label{}	
	\end{subfigure}\hfill
	\hspace{-10em}
	\begin{subfigure}[t]{0.5\textwidth}
		\centering
		\includegraphics[width=0.7\textwidth]{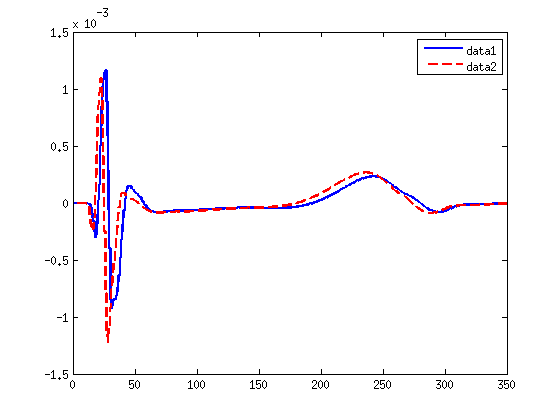}
		\vspace{-1em}
		\caption{}
	\end{subfigure}\\
	\vspace{-1em}
	\caption{: \textbf{$v$(first row), $[Ca^{+2}]_i$(second row), $T_A$(third row), $\lambda$ (fourth row) and $\frac{d \lambda}{dt}$ (fifth row), at M1(first column), M2 (second column), for the case with $I_{sac}$ (second column) for $[K^+]_o =20$, with two different sizes of ischemic regions, data1 ($[0.1563, 0.25]^2$), data2 ($[0.0938, 0.3125]^2$).}}
	\label{incr_size_Ko_Isac}
\end{figure*}

Next, the effect of increase in the size of the ischemic region from $[0.1563, 0.25]^2$ to 
$[0.0938, 0.3125]^2$ on action potential, intracellular calcium ions $[Ca^{+2}]_i$, active tension $T_A$, stretch along the fiber $\lambda$ and stretch rate,  
$\frac{d \lambda}{dt}$ are presented in Fig. \eqref{incr_size_Ko_Isac} for the case $[K^+]_o =20$mM. 
Clearly the spread of severe Hyperkalemic ischemic zone further raises already increased resting potential levels, at $M1$, from  -50mV to -45mV respectively. Also the growth in ischemic size leads to faster repolarization with a reduced APD. There is almost $9\%$ further drop in APD is noticed with a factor of five times increase in the ischemic subregion size. From the action potential plot corresponding to point $M2$, it is visible that there is almost $5\%$ drop in APD for this neighboring healthy cell. Thus, it is amply clear that in the proximity of ischemic subregion the CEA in terms of action potential and hence the ECG of healthy cardiac tissue is affected with the severity of Hyperkalemia. This effect in the ischemic and healthy cells will increase with the spread of ischemic region. 

It is also noticed that there is variation in the waveform of $[Ca^{+2}]_i$ at M1 while it is going to rest, but at M2 the waveform changes during all the phases. The variation in the stretch along the fiber and the stretch rate is more in the neighboring region point M2 with the expansion of the ischemic region while the change at the point M1 inside the ischemic region is negligible. In short, expanding the ischemic region size affect the mechanical properties of the neighboring healthy cells more intensely. 

Hence, it could be concluded that addition of $I_{sac}$ leads to ECG change and mechanical contraction as it reduces the APD and affects the contractility and the length of length of myocytes of healthy and ischemic parts of the human cardiac tissue. It is also observed that the severity of Hyperkalemia leads to reduction in APD, elevation in resting potential, affects the contractile force and contractility and the stretch activated channels, hence affects the QRS complex and QT interval of ECG and the mechanical contraction or can say the electro-mechanical activity of healthy and ischemic regions of human cardiac tissue. All these effects on the electro-mechanical activity of a human heart gets more intense as the ischemic regions expands or more cardiac cells becomes ischemic.

\textbf{Note:} As it is concluded that the addition of $I_{sac}$ plays an important role in the electro-mechanical activity of a human cardiac tissue. So, next results will be with $I_{sac}$ case only.

\subsubsection{{(a) Hypoxia (Effect of ischemia when only $f_{ATP}$ varies)}}
\begin{figure*}
	\begin{subfigure}[t]{0.5\textwidth}
		\centering
		\includegraphics[width=0.7\textwidth]{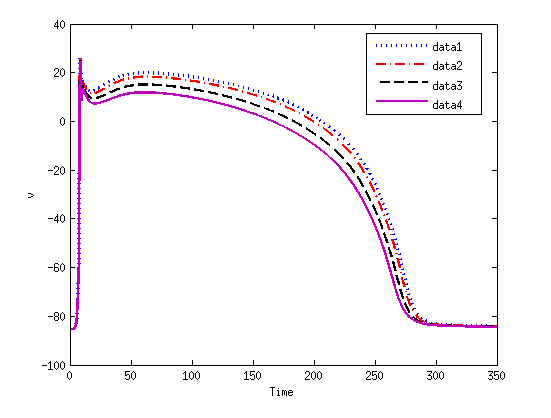}
		\vspace{-1em}
		\caption{}
		\label{AP_Pt1_diff_f}
	\end{subfigure}\hfill	
	\hspace{-10em}
	\begin{subfigure}[t]{0.5\textwidth}
		\centering
		\includegraphics[width=0.7\textwidth]{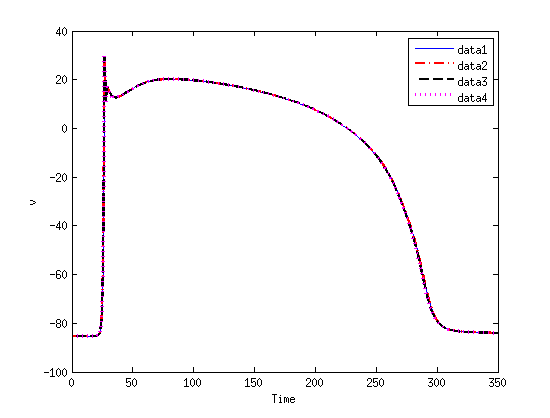}
		\vspace{-1em}
		\caption{}
		\label{AP_Pt2_diff_f}
	\end{subfigure}\\
	\begin{subfigure}[t]{0.5\textwidth}
		\centering
		\includegraphics[width=0.7\textwidth]{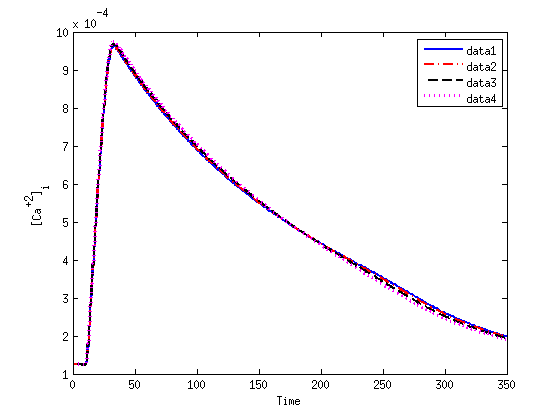}
		\vspace{-1em}
		\caption{}
		\label{Cai_Pt1_diff_f}
	\end{subfigure}\hfill
	\hspace{-10em}
	\begin{subfigure}[t]{0.5\textwidth}
		\centering
		\includegraphics[width=0.7\textwidth]{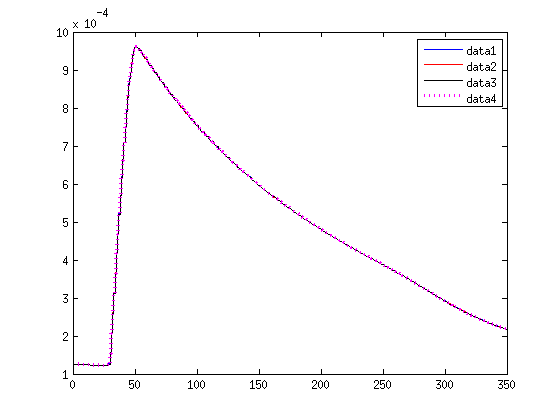}
		\vspace{-1em}
		\caption{}
		\label{Cai_Pt2_diff_f}
	\end{subfigure}\\
	\begin{subfigure}[t]{0.5\textwidth}
		\centering
		\includegraphics[width=0.7\textwidth]{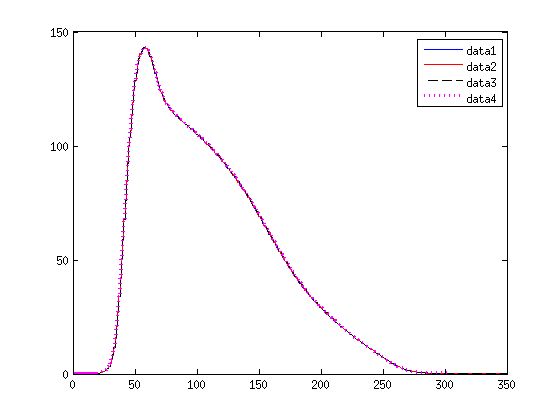}
		\vspace{-1em}
		\caption{}
		\label{Ta_Pt1_diff_f}
	\end{subfigure}\hfill
	\hspace{-10em}
	\begin{subfigure}[t]{0.5\textwidth}
		\centering
		\includegraphics[width=0.7\textwidth]{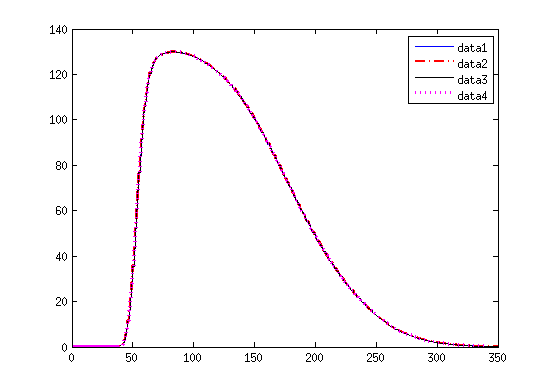}
		\vspace{-1em}
		\caption{}
		\label{Ta_Pt3_diff_f}
	\end{subfigure}\\
	\begin{subfigure}[t]{0.5\textwidth}
		\centering
		\includegraphics[width=0.7\textwidth]{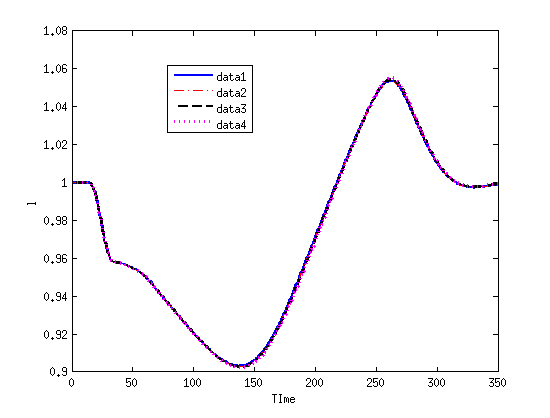}
		\vspace{-1em}
		\caption{}
		\label{lambda_Pt1_diff_f}
		\label{}
	\end{subfigure}\hfill
	\hspace{-10em}
	\begin{subfigure}[t]{0.5\textwidth}
		\centering
		\includegraphics[width=0.7\textwidth]{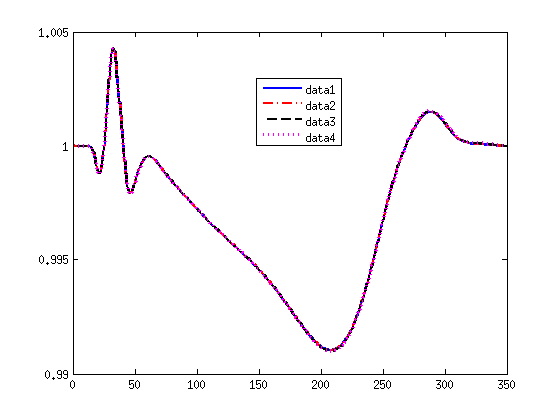}
		\vspace{-1em}
		\caption{}
		\label{lambda_Pt3_diff_f}
	\end{subfigure}\\
	\begin{subfigure}[t]{0.5\textwidth}
		\centering
		\includegraphics[width=0.7\textwidth]{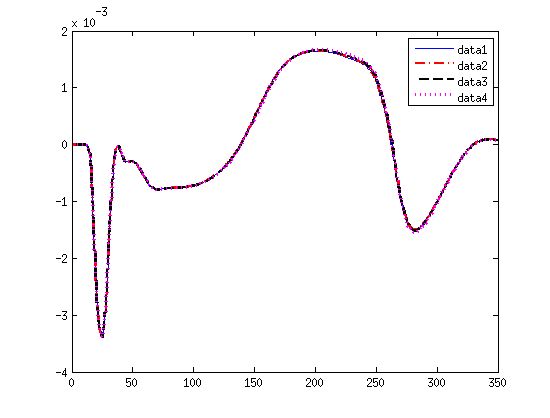}
		\vspace{-2em}
		\caption{}
		\label{dlambda_Pt1_diff_f}	
	\end{subfigure}\hfill
	\hspace{-10em}
	\begin{subfigure}[t]{0.5\textwidth}
		\centering
		\includegraphics[width=0.7\textwidth]{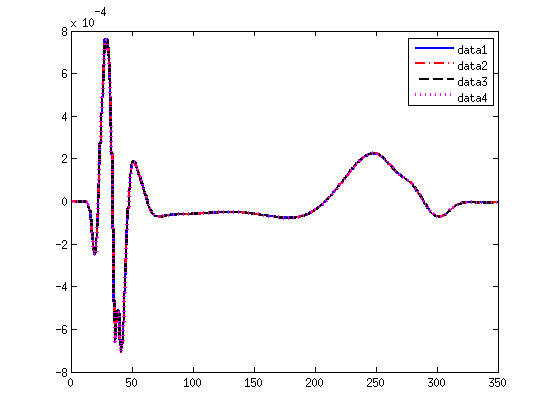}
		\vspace{-2em}
		\caption{}
		\label{dlambda_Pt3_diff_f}
	\end{subfigure}\\
	\vspace{-2em}
	\caption{:\textbf{$v$(first row), $[Ca^{+2}]_i$(second row), $T_A$(third row), $\lambda$ (fourth row) and $\frac{d \lambda}{dt}$ (fifth row), at M1(first column), M2 (second column), for different values of $f_{ATP}$. data1 ($f_{ATP}=0$), data2 ($f_{ATP}=0.1\%$), data3 ($f_{ATP}=0.3\%$), data4 ($f_{ATP}=0.5\%$).}}
	\label{diff_f}
\end{figure*}

$f_{ATP}$ is taken in the range $0.1-0.5\%$. Influence of the $f_{ATP}$ values (oxygen level) on the action potential and the mechanical parameters $[Ca^{+2}]_i$, active tension $(T_A)$, stretch along the fiber $(\lambda_l)$ and stretch rate $(\frac{d\lambda_l}{dt})$ at the two points $M1$ and $M2$ (described above) of the cardiac tissue are estimated and presented in Fig. \eqref{diff_f}. From Fig. \eqref{AP_Pt1_diff_f}, it can be seen that as the oxygen level in the ischemic region of the cardiac tissue decreases (or $f_{ATP}$ value increased), delay in the closing of L-type $Ca^{+2}$ channels is noticed, and further rapid close of the $K^+$ channels takes place, therefore the plateau phase and then the repolarization phase of the cardiac action potential get affected. Therefore, APD also reduces with the increase in the strength of Hypoxia. While from the Fig. \eqref{AP_Pt2_diff_f}, it can be seen that influence of the increase in $f_{ATP}$ values is negligible on the neighborhood point M2 on the non-ischemic region of cardiac tissue.

Further from figures \eqref{Cai_Pt1_diff_f} - \eqref{dlambda_Pt3_diff_f}, it is visulaized that there is negligible change in the waveform of the mechanical parameters, intracellular calcium concentration ($[Ca^{+2}]_i$), active tension $(T_A)$, stretch along the fiber $(\lambda_l)$ and stretch rate $(\frac{d\lambda_l}{dt})$ at the points $M1$ and $M2$. 

\begin{figure*}
	\hspace{-5em}
	\begin{subfigure}[t]{0.5\textwidth}
		\centering
		\includegraphics[width=0.6\textwidth]{ACTI_Ko_norm_with_ISAC_CONV.png}
		\vspace{-1.2em}
		\caption{}
		\label{ACTI_f_0}
	\end{subfigure} 
	\hspace{-12em}
	\begin{subfigure}[t]{0.5\textwidth}
		\centering
		\includegraphics[width=0.6\textwidth]{REPO_Ko_norm_with_ISAC_CONV.png}
		\vspace{-1.2em}
		\caption{}
		\label{REPO_f_0}
	\end{subfigure}\hspace{-12em}
	\begin{subfigure}[t]{0.5\textwidth}
		\centering
		\includegraphics[width=0.6\textwidth]{APD_Ko_norm_with_ISAC_CONV.png}
		\vspace{-1.2em}
		\caption{}
		\label{APD_f_0}
	\end{subfigure}\\
	\begin{subfigure}[t]{0.5\textwidth}
		\includegraphics[width=0.6\textwidth]{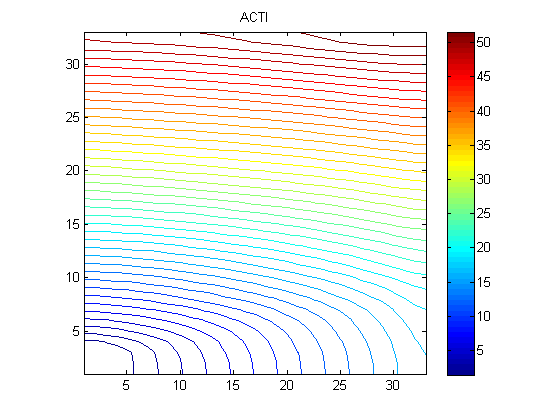}
		\vspace{-1.2em}
		\caption{}
		\label{ACTI_f_30}
	\end{subfigure}\hspace{-17em}
	\begin{subfigure}[t]{0.5\textwidth}
		\centering
		\includegraphics[width=0.6\textwidth]{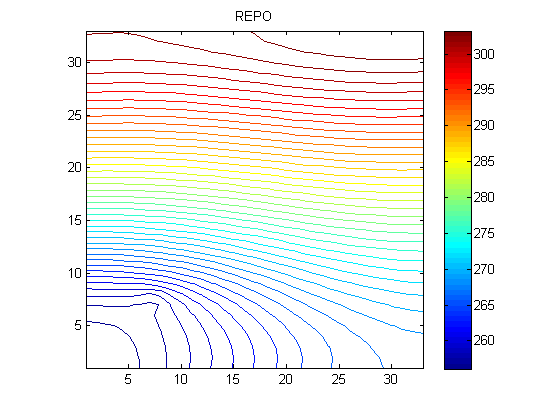}
		\vspace{-1.2em}
		\caption{}
		\label{REPO_f_30}
	\end{subfigure}\hspace{-12em}
	\begin{subfigure}[t]{0.5\textwidth}
		\centering
		\includegraphics[width=0.6\textwidth]{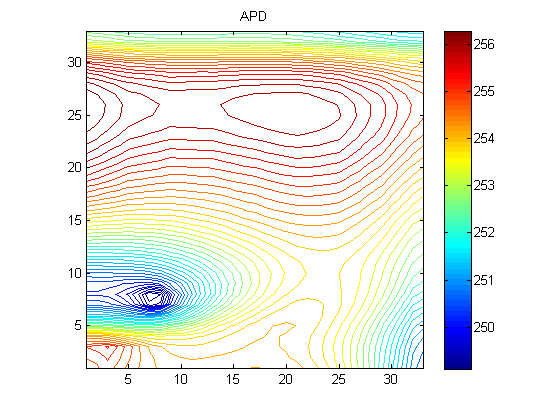}
		\vspace{-1.2em}
		\caption{}
		\label{APD_f_30}
	\end{subfigure} \\
	\begin{subfigure}[t]{0.5\textwidth}
		\includegraphics[width=0.6\textwidth]{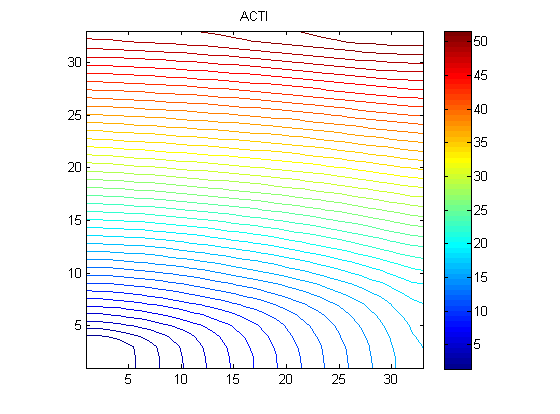}
		\vspace{-1.2em}
		\caption{}
		\label{ACTI_f_50}
	\end{subfigure}\hspace{-17em}
	\begin{subfigure}[t]{0.5\textwidth}
		\centering
		\includegraphics[width=0.6\textwidth]{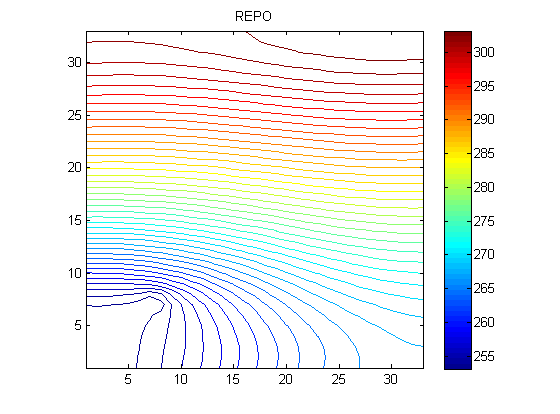}
		\vspace{-1.2em}
		\caption{}
		\label{REPO_f_50}
	\end{subfigure}\hspace{-12em}
	\begin{subfigure}[t]{0.5\textwidth}
		\centering
		\includegraphics[width=0.6\textwidth]{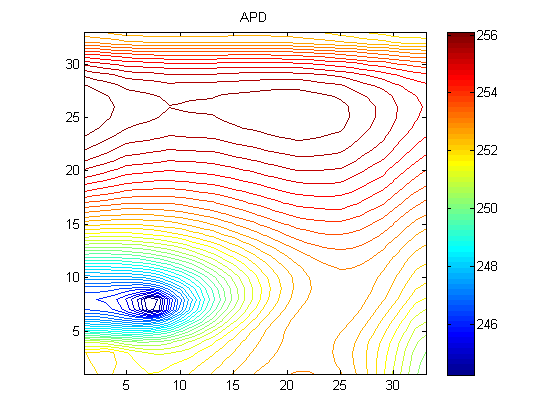}
		\vspace{-1.2em}
		\caption{}
		\label{APD_f_50}
	\end{subfigure}
	\caption{:\textbf{AT (first column), RT (second column), APD (third column) for different values of $f_{ATP}=0$(first row), $f_{ATP}=0.1\%$(second row), $f_{ATP}=0.3\%$(third row), $f_{ATP}=0.5\%$(fourth row).}}
	\label{ACTI_REPO_APD_f}
\end{figure*}

From Fig. \eqref{ACTI_f_0}-\eqref{ACTI_f_50}, it is clear that the activation time remain unaffected by $f_{ATP}$ and hence the isochrones of AT corresponding to these cases turn out to be same in the entire domain including those in ischemic region. The front near the ischemic region is quasi elliptic and becomes flat after that.

RT isochrones show that as we increase the $f_{ATP}$
values in the ischemic region i.e. as the oxygen level in the ischemic region of the cardiac tissue
decreases, the cardiac cells in the ischemic region repolarize earlier i.e. the repolarization 
time of ischemic cell decreases with the increase in $f_{ATP}$ values while the repolarization
time for the non-ischemic cardiac cells remains same. The loss
of quasi-ellipticity and the folding of repolarization time isochrones especially in the ischemic
subregion depict an accelerated repolarization wave propagation in the human cardiac tissue.
The isochrone lines for APD shows that APD in the ischemic region reduces as we are increasing $f_{ATP}$ values. 

\begin{figure*}
	\hspace{-5em}
	\begin{subfigure}[t]{0.5\textwidth}
		\centering
		\includegraphics[width=0.6\textwidth]{v90_k0_5p4_ISAC_CONV.jpeg}
		\vspace{-1.2em}
		\caption{}
		\label{contour_f_Nor_90}
	\end{subfigure} 
	\hspace{-12em}
	\begin{subfigure}[t]{0.5\textwidth}
		\centering
		\includegraphics[width=0.6\textwidth]{v150_k0_5p4_ISAC_CONV.jpeg}
		\vspace{-1.2em}
		\caption{}
		\label{contour_f_Nor_150}
	\end{subfigure}\hspace{-12em}
	\begin{subfigure}[t]{0.5\textwidth}
		\centering
		\includegraphics[width=0.6\textwidth]{v240_k0_5p4_ISAC_CONV.jpeg}
		\vspace{-1.2em}
		\caption{}
		\label{contour_f_Nor_240}
	\end{subfigure} \\
	\begin{subfigure}[t]{0.5\textwidth}
		\includegraphics[width=0.6\textwidth]{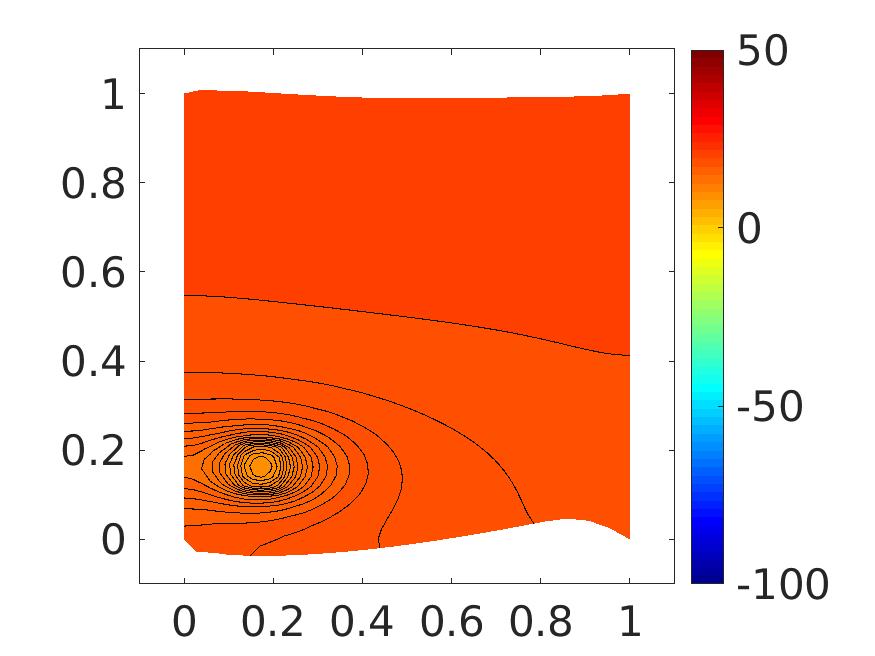}
		\vspace{-1.2em}
		\caption{}
		\label{contour_f_30_90}
	\end{subfigure}\hspace{-17em}
	\begin{subfigure}[t]{0.5\textwidth}
		\centering
		\includegraphics[width=0.6\textwidth]{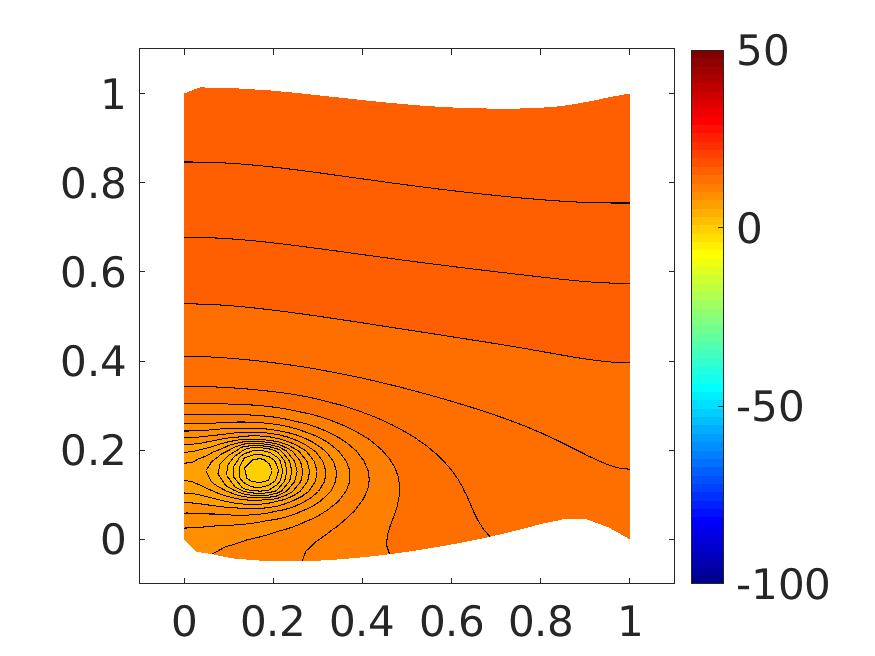}
		\vspace{-1.2em}
		\caption{}
		\label{contour_f_30_150}
	\end{subfigure}\hspace{-12em}
	\begin{subfigure}[t]{0.5\textwidth}
		\centering
		\includegraphics[width=0.6\textwidth]{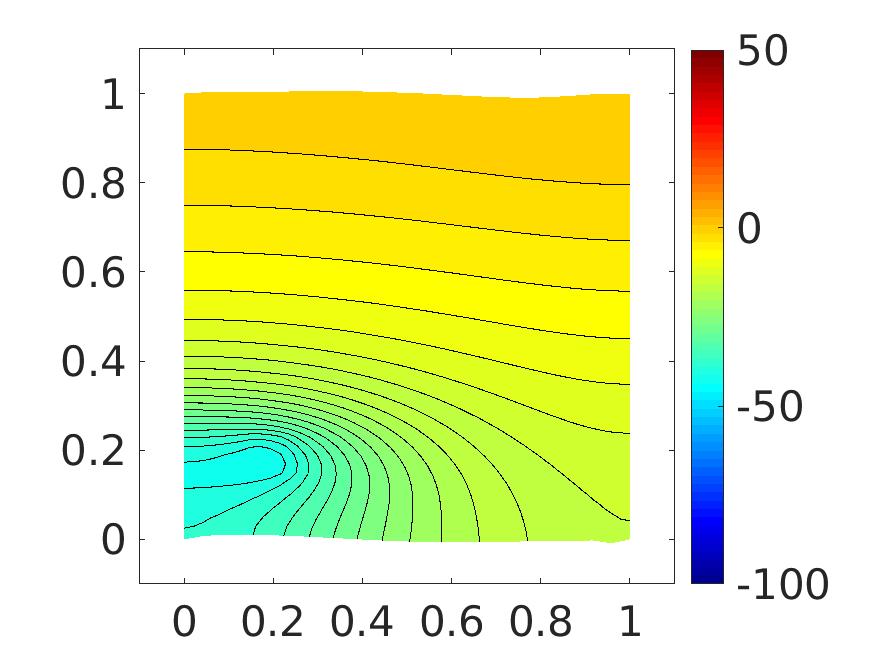}
		\vspace{-1.2em}
		\caption{}
		\label{contour_f_30_240}
	\end{subfigure}\\
	\begin{subfigure}[t]{0.5\textwidth}
		\includegraphics[width=0.6\textwidth]{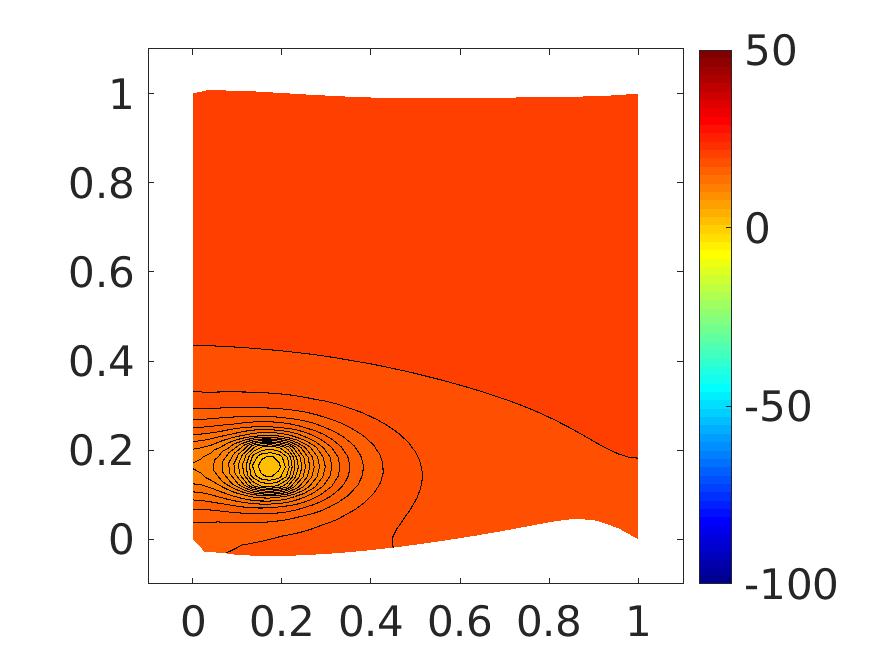}
		\vspace{-1.2em}
		\caption{}
		\label{contour_f_50_90}
	\end{subfigure}\hspace{-17em}
	\begin{subfigure}[t]{0.5\textwidth}
		\centering
		\includegraphics[width=0.6\textwidth]{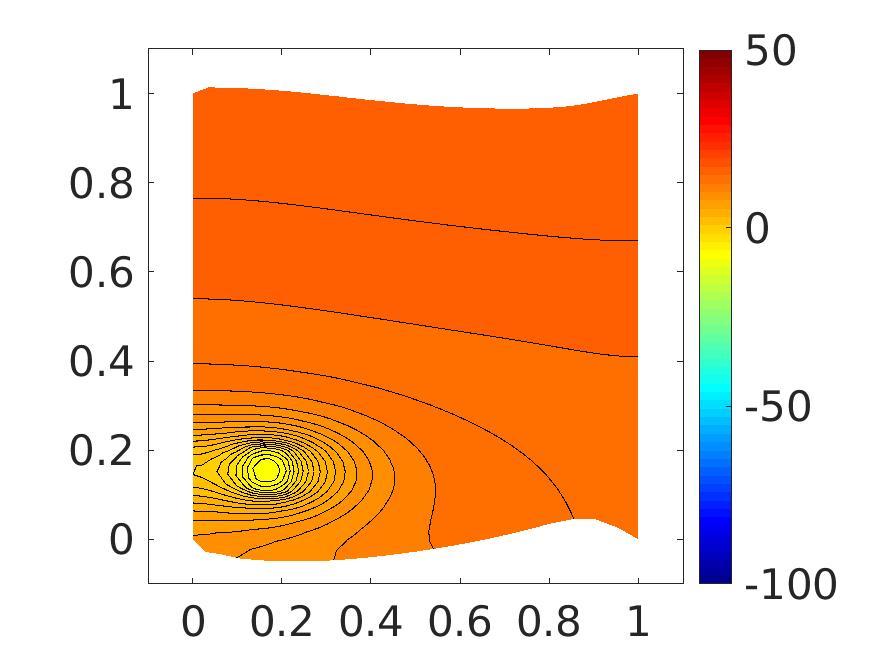}
		\vspace{-1.2em}
		\caption{}
		\label{contour_f_50_150}
	\end{subfigure}\hspace{-12em}
	\begin{subfigure}[t]{0.5\textwidth}
		\centering
		\includegraphics[width=0.6\textwidth]{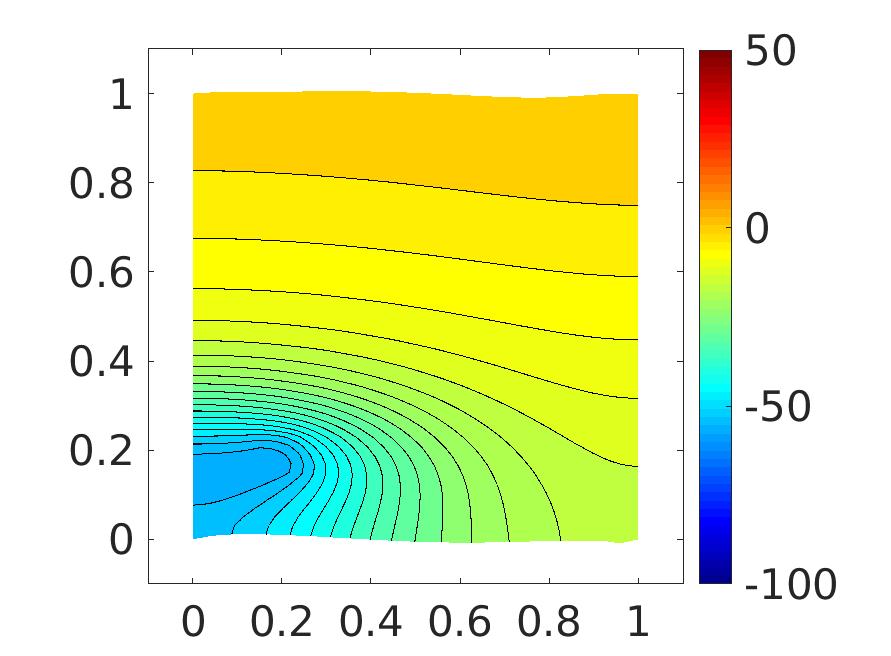}
		\vspace{-1.2em}
		\caption{}
		\label{contour_f_50_240}
	\end{subfigure}
	\vspace{0em}
	\caption{:\textbf{AP contours in a deforming domain for different values of $f_{APT}$ =0 (first row), $f_{APT}=0.3\%$ (second row), $f_{APT}=0.5\%$ (third row), at time t=90 (first column), t=150 (seconf column), t=240 (third column).}}
	\label{AP_contours_f}
\end{figure*}

In Fig. \eqref{contour_f_Nor_90}-\eqref{contour_f_50_240}, we presented the action potential contours in deforming cardiac tissue with 
the increase in the strength of Hypoxia. We can visualize the change in the deformation of cardiac domain with the increasing Hypoxic strength. Thus, change in oxygen level affects the contraction and expansion of ventricle of human heart.

\begin{figure*}
	\centering
	\begin{subfigure}[t]{0.5\textwidth}
		\centering
		\includegraphics[width=0.7\textwidth]{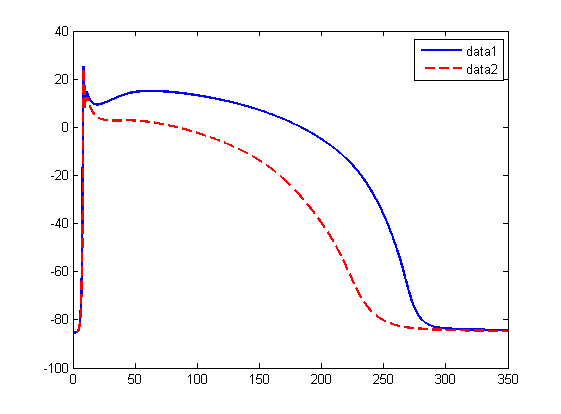}
		\vspace{-1em}
		\caption{}
		\label{AP_Pt1_f_30_incr_size}
	\end{subfigure}\hfill	
	\hspace{-10em}
	\begin{subfigure}[t]{0.5\textwidth}
		\centering
		\includegraphics[width=0.7\textwidth]{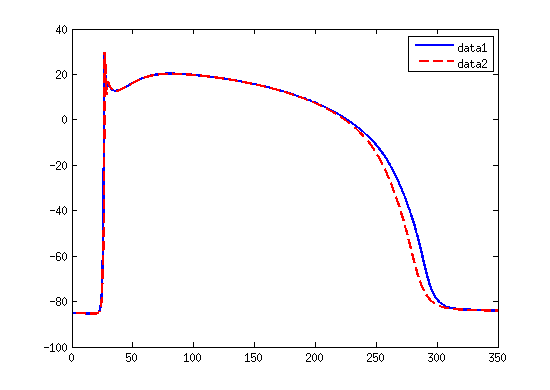}
		\vspace{-1em}
		\caption{}
		\label{AP_Pt2_f_30_incr_size}
	\end{subfigure}\\
	\begin{subfigure}[t]{0.5\textwidth}
		\centering
		\includegraphics[width=0.7\textwidth]{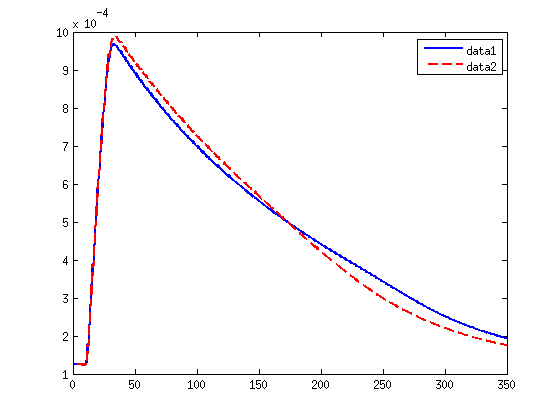}
		\vspace{-1em}
		\caption{}
		\label{Cai_Pt1_f_30_incr_size}
	\end{subfigure}\hfill
	\hspace{-10em}
	\begin{subfigure}[t]{0.5\textwidth}
		\centering
		\includegraphics[width=0.7\textwidth]{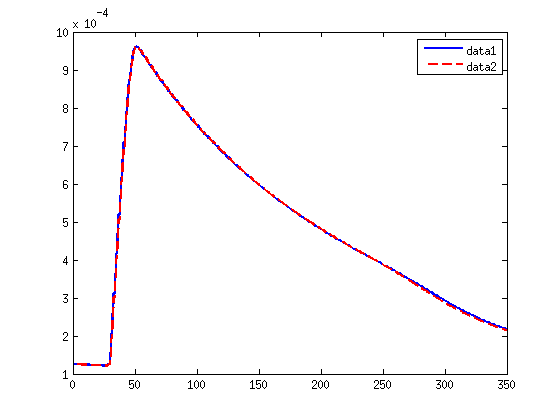}
		\vspace{-1em}
		\caption{}
		\label{Cai_Pt2_f_30_incr_size}
	\end{subfigure}\\
	\begin{subfigure}[t]{0.5\textwidth}
		\centering
		\includegraphics[width=0.7\textwidth]{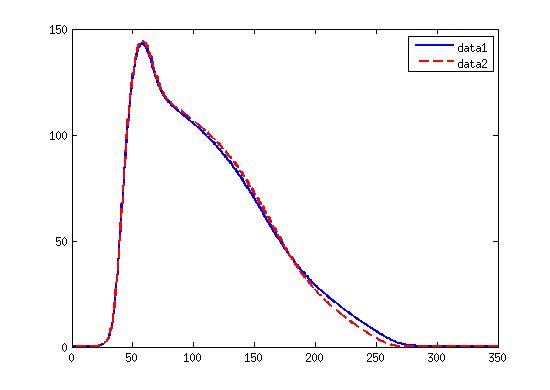}
		\vspace{-1em}
		\caption{}
		\label{Ta_Pt1_f_30_incr_size}
	\end{subfigure}\hfill
	\hspace{-10em}
	\begin{subfigure}[t]{0.5\textwidth}
		\centering
		\includegraphics[width=0.7\textwidth]{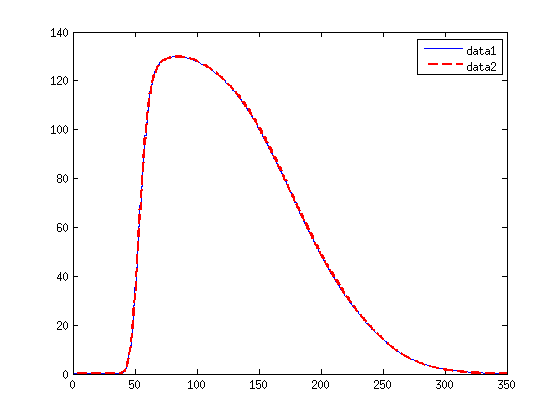}
		\vspace{-1em}
		\caption{}
		\label{Ta_Pt2_f_30_incr_size}
	\end{subfigure}\\
	\begin{subfigure}[t]{0.5\textwidth}
		\centering
		\includegraphics[width=0.7\textwidth]{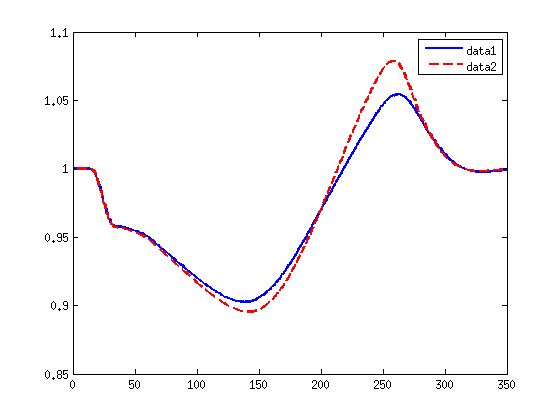}
		\vspace{-1em}
		\caption{}
		\label{lambda_Pt1_f_30_incr_size}
	\end{subfigure}\hfill
	\hspace{-10em}
	\begin{subfigure}[t]{0.5\textwidth}
		\centering
		\includegraphics[width=0.7\textwidth]{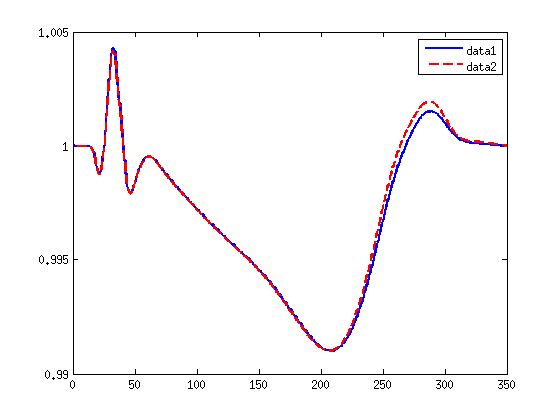}
		\vspace{-1em}
		\caption{}
		\label{lambda_Pt2_f_30_incr_size}
	\end{subfigure}\\
	\begin{subfigure}[t]{0.5\textwidth}
		\centering
		\includegraphics[width=0.7\textwidth]{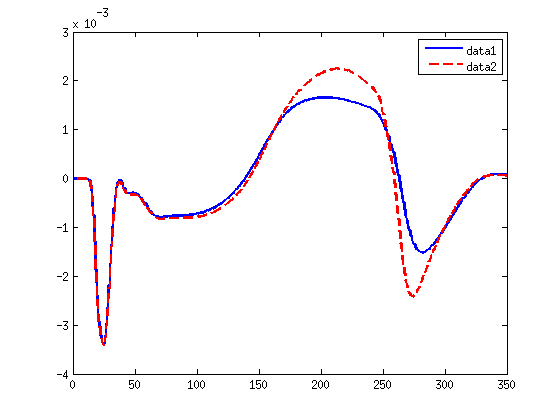}
		\vspace{-1em}
		\caption{}
		\label{dlambda_Pt1_f_30_incr_size}	
	\end{subfigure}\hfill
	\hspace{-10em}
	\begin{subfigure}[t]{0.5\textwidth}
		\centering
		\includegraphics[width=0.7\textwidth]{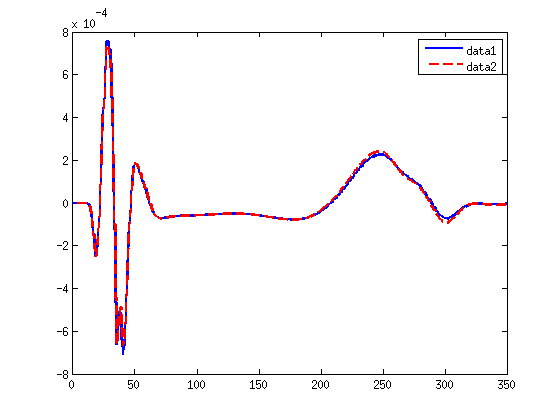}
		\vspace{-1em}
		\caption{}
		\label{dlambda_Pt2_f_30_incr_size}
	\end{subfigure}\\
	\caption{:\textbf{$v$(first row), $[Ca^{+2}]_i$(second row), $T_A$(third row), $\lambda$ (fourth row) and $\frac{d \lambda}{dt}$ (fifth row), at M1(first column), M2 (second column), for the case with $I_{sac}$ (second column) for $f_{ATP} =0.3\%$, with two different sizes of ischemic regions, data1 ($[0.1563, 0.25]^2$), data2 ($[0.0938, 0.3125]^2$).}}
	\label{incr_size_f_Isac}
\end{figure*}

Next, as we increase the size of the ischemic region from $[0.1563, 0.25]^2$ to $[0.0938, 0.3125]^2$, action potential at the points M1 and M2 get influenced as shown in Fig \eqref{incr_size_f_Isac}. From the Fig \eqref{AP_Pt1_f_30_incr_size}, we can conclude that, at $M_1$ (ischemic region point), with the increase in the size of the ischemic region there is further delay in the closing of L-type calcium channels. Thus, a delay in plateau phase and then early repolarization of the action potential takes place. APD also reduces with increase in the size of the ischemic region. From the action potential plot corresponding to point $M_2$ in Fig. \eqref{AP_Pt2_f_30_incr_size}, a delay in plateau phase and early repolarization of the action potential and reduction in the APD of these neighboring cells is also clearly visible. At $M1$ there is $11\%$ drop in APD with a factor of five times increase in the ischemic region size. Depending on the distance of the neighboring points from the ischemic region, more than $5\%$ drop in APD with a factor of five times increase in the ischemic region size is noticed.

Now, from Fig. \eqref{Cai_Pt1_f_30_incr_size} - \eqref{dlambda_Pt2_f_30_incr_size}, it is visible that increasing the size of the ischemic region also considerably affects the waveforms of the mechanical parameters $[Ca^{+2}]_i$, $T_A$, $\lambda$ and $\frac{d \lambda}{dt}$ corresponding to the point M1 than those of the point M2. From the figures \eqref{lambda_Pt1_f_30_incr_size} and \eqref{lambda_Pt2_f_30_incr_size}, it can be noticed that an increase in the size of the ischemic region by a factor of five leads to additional shortening of fiber length during systolic contraction and an enhanced stretch in the fibers post systolic contraction.
At $M1$ $3-4\%$ change in the length of the myocytes relative to the resting length is noticed and at $M_2$ only a $1-2\%$ change in the length of the myocytes relative to the resting length is noticed. Therefore, the stretch rate $\frac{d \lambda}{dt}$ also changes at the points $M_1$ and $M_2$. There is approximately $45\%$ variation in the stretch rate $\frac{d \lambda}{dt}$ at $M1$ and approximately $5-25\%$ variation in the neighboring points is detected. 

Thus, the cardiac electro-mechanical activity mathematically generated by PDE-ODE system leads to tracing of ionic channel dynamics, intracellular calcium ion abnormalities including the tissue structural influences. It may be noted that the ionic channels corresponding to each cardiac cell adds to the generation of cardiac action potential and also forms the basis for calcium induction, calcium release and cardiac mechanical contraction. The change in the action potential phase and hence in ECG pattern affects both the electrical and mechanical function of heart. These numerical results help in tracing the affected action potential phase, leading to a disturbance in ECG pattern. It also helps in identification of corresponding affected ionic channel dynamics. Thus, the current analysis can provide a useful guideline for the treatment of Hyperkalemia and Hypoxia.

Thus, the heart failure due to the  cardiac disease from the abnormal functioning of  the coupled electrical and mechanical activity can be traced through the virtual  reconstruction of ionic channels and intracellular calcium ion abnormalities with cardiac structural reconstruction through the proposed numerical modeling. The ionic channels corresponding to each cardiac cell adds to the generation of cardiac action potential and also forms the basis for calcium induction, calcium release and cardiac mechanical contraction. The change in the action potential phase and hence in ECG pattern affects both the electrical and mechanical function of heart.  

\section{Conclusion}
Local ischemia in the deformed human cardiac tissue is modeled by varying the ischemic parameters value in the ischemic subregion of the cardiac tissue domain. Two types of ischemic effects, namely, Hyperkalemia and Hypoxia, are considered in this work. Monodomain model in a deforming domain is taken with the TT06 human cell level model. The coupled electro-mechanical PDEs-ODEs non-linear system of equations are solved numerically using linear finite elements in space and backward Euler finite difference scheme in time. We examine the cardiac electrical and mechanical activity in terms of the action potential ($v$) and intracellular calcium ion concentration $[Ca^{+2}]_i$, active tension, ($T_A$), stretch ($\lambda$), stretch rate ($ \frac{d \lambda}{dt}$), in several cases for local ischemia. 
We investigated the effect of varying strength of Hyperkalemia and Hypoxia in the ischemic subregion of cardiac tissue, on the electrical and mechanical activity of healthy and ischemic zones in the cardiac muscle. We also 
investigated the impact of increasing the size of the ischemic region on the electrical and mechanical parameters of neighboring cells.

(a) With the severity of Hyperkalemia the concentration of extracellular potassium ions increase leading to, (i) Fall in availability of sodium channels and a decrease in inward sodium current, (ii) AP comes to a resting state significantly prior 
to reaching the normal resting potential level, and, (iii) Increased potassium conductance leading to shortening of RT in turn shortens QT interval is noticed.

These changes in cellular ionic dynamics not only further alters the AT, RP, RT and APD of affected cells with the spread of Hyperkalemic region  but also increasingly alters the AP of healthy cells in its vicinity. Also, this spread of Hyperkalemic region alters the waveform of the mechanical parameters $[Ca^{+2}]_i, T_A, \lambda$ ,$ \frac{d \lambda}{dt}$ of the ischemic and the neighboring healthy cells.

Thus, the severity of Hyperkalemia leads to reduction in APD, elevation in resting potential, affects the contractile force and contractility and the stretch activated channels, hence affects the QRS complex and QT interval of ECG and the mechanical contraction or can say the electro-mechanical activity of healthy and ischemic regions of human cardiac tissue. All these effects on the electro-mechanical activity of a human heart gets more intense as the ischemic regions expands or more cardiac cells becomes ischemic.

(b) Severity of Hypoxic ischemia, with a reduction in intracellular ATP concentration, affect both influx of calcium ions 
and efflux of potassium ions and alters the normal functioning of calcium channels and balanced potassium channels. 
The increase in intracellular $Ca^{+2}$ levels completely disturbs the Phase1, Phase2 and Phase 3 of AP. As the Hypoxically
affected region increases, the Plateau phase, repolarization phase and the APD corresponding to the AP of healthy cells in its
vicinity are affected and the ionic dynamics of Hypoxically already degenerated cells get further disturbed. It is also visible that increasing the size of the ischemic region by a factor of five considerably affects the waveforms of the mechanical parameters $[Ca^{+2}]_i$, $T_A$, $\lambda$ and $\frac{d \lambda}{dt}$. There is approximately $45\%$ variation in the stretch rate $\frac{d \lambda}{dt}$ at $M1$ and approximately $5-25\%$ variation in the neighboring points is detected. 

Thus, the change in the action potential phase and hence in ECG pattern affects both the electrical and mechanical function of heart. These numerical results help in tracing the affected action potential phase, leading to a disturbance in ECG pattern. It also helps in identification of corresponding affected ionic channel dynamics. Thus, the current analysis can provide a useful guideline for the treatment of Hyperkalemia and Hypoxia.

\subsection*{Financial disclosure}

We would like to thank the DST for support through Inspire Fellowship, ID no. is IF130906. We would also like to thanks the ICTP-INdAM.

\subsection*{CONFLICTS OF INTEREST}

The authors declare no potential conflict of interests.

\bibliography{reference}
\bibliographystyle{alpha} 

\end{document}